\documentclass[10pt,a4paper,oneside,leqno]{article}
\usepackage{amsfonts,amssymb}
\usepackage{amscd}
\newtheorem{defn}{Definition}

\newtheorem{observen}{Observation}

\vspace{0.5cm}
 4
\font\ebf=cmbx8
\font\erm=cmr8

\usepackage{graphicx, graphics}
\usepackage{color}

\newcommand{\os}{\oplus\!\!\to}

\newcommand{\fnomial}[2]{ {{#1} \choose {#2}}_F }

\begin{document}
\begin{center}
	\noindent { \textsc{ Natural join construction of graded posets versus ordinal sum and discrete  hyper- boxes}}  \\ 
	\vspace{0.3cm}
	\vspace{0.3cm}
	\noindent Andrzej Krzysztof Kwa\'sniewski \\
	\vspace{0.2cm}
	\noindent {\erm Member of the Institute of Combinatorics and its Applications  }\\
{\erm High School of Mathematics and Applied Informatics} \\
	{\erm  Kamienna 17, PL-15-021 Bia\l ystok, Poland }\noindent\\
	\noindent {\erm e-mail: kwandr@gmail.com}\\
\vspace{0.3cm}
\end{center}

\noindent {\ebf Abstract:}
\vspace{0.1cm}
\noindent {\small One introduces here the natural join $P \os Q$ of graded posets $\left\langle P,\leq_P \right\rangle$ and $\left\langle Q,\leq_Q \right\rangle$ with correspondingly  maximal (in $P$) and minimal (in $Q$) sets  being identical as expressed by ordinal sum  $P\oplus Q$ (see [1]) apart from other definition in [2,3,4,5] and due to that one arrives at a simple proof of the M{\"{o}}bius function formula for cobweb posets.\\
We also quote the other authors  explicit formulas  for the  zeta matrix and the  inverse  of zeta matrix for any graded posets  with the finite set of minimal elements from [2] and [3,4,5]. These formulas are  based on the formulas for cobweb posets and their Hasse diagrams - graphs named KoDAGs, which are interpreted as chains of binary \textbf{complete} (or universal) relations -  joined by the natural join operation.\\
Natural join of two independent sets is therefore the ordinal sum [1] of this trivially ordered posets represented also by directed bi-clique named di-biclique and correspondingly by  their Hasse diagrams - graphs named KoDAGs.\\
Such cobweb posets and equivalently their Hasse diagrams - graphs named KoDAGs -  are also encoded by discrete hyper-boxes and the natural join operation of such discrete hyper-boxes is just Cartesian product of them accompanied with projection out of - sine qua non - common faces. All graded posets with no mute vertices in their Hasse diagrams [2] (i.e. no vertex has  in-degree or out-degree equal zero) are natural join of chain of relations and may be at the same time  interpreted an $n-ary$, $n \in N \cup \left\{\infty \right\}$ relation. The Whitney numbers and characteristic polynomials 
explicit formulas for  cobweb posets are derived.
}

\vspace{0.3cm}


\vspace{0.2cm}

\noindent AMS Classification Numbers: 06A06 ,05B20, 05C75  

\vspace{0.2cm}

\noindent  This is The Internet Gian-Carlo Rota Polish Seminar article, No \textbf{8}, \textbf{Subject 4}, \textbf{2009-07-13}, \emph{http://ii.uwb.edu.pl/akk/sem/sem\_rota.htm}\\


\section{Reference information on ordinal sum and natural join of graded ordered sets and M{\"{o}}bius  function  formula.} 


\noindent \textbf{1.1.} \textbf{Ponderables. [2,3,4,5]} 

\vspace{0.2cm}

\noindent We shall here take for granted the notation and the results of [2,3]. In particular $\left\langle \Pi,\leq \right\rangle$ denotes cobweb partial order set (cobweb poset)  while   $I(\Pi,R)$ denotes its incidence algebra over the ring $R$.  Correspondingly  $\left\langle P,\leq \right\rangle$ denotes arbitrary graded poset while   $I(P,R)$ denotes its incidence algebra over the ring $R$. 
$[n]\equiv \left\{1,2,...,n\right\}$ , for example $[k_F]\equiv \left\{1,2,...,k_F\right\}$.

\vspace{0.1cm}

\noindent $R$ might be taken to be Boolean algebra $2^{\left\{1 \right\}}$ , the field $Z_2=\left\{0,1\right\}$, the ring of integers $Z$ or real or complex or $p$-adic fields. The present article is the next one in a series of papers listed in  order of appearence and these are: [5],[4],[3],[2]. The abbreviation \textbf{DAG} $\equiv$ \textbf{D}irected \textbf{A}cyclic \textbf{G}raph. \\

\vspace{0.1cm}

\noindent Inspired by Gaussian integers sequence notation $\left\{ n_q \right\}_{n\geq 0}$ - the authors   upside down notation  is used throughout this paper i.e. $F_n \equiv  n_F$. 
\textcolor{blue}{\textbf{The Upside Down Notation }} was used since last century effectively (see [2-22] and [29-48]).\\
Through all the paper  $F$ denotes a natural numbers valued sequence $\left\{ n_F \right\}_{n\geq 0} \equiv  \left\{ F_n \right\}_{n\geq 0}$  sometimes specified to be Fibonacci or others - if needed. Among many consequences of this is that \textbf{\textbf{graded posets}} ($\equiv$ their cover relation  digraphs $ \Longleftrightarrow$ Hasse diagrams) \textbf{\textit{are connected}} and sets of their minimal elements are finite.\\

\vspace{0.1cm}

\begin{defn}
\noindent Let  $F = \left\langle k_F \right\rangle_{k=0}^n$ be an arbitrary natural numbers valued sequence, where $n\in N \cup \left\{0\right\}\cup \left\{\infty\right\}$. We say that the graded poset $P = (\Phi,\leq)$ is \textcolor{red}{\textbf{denominated}} (encoded=labeled) by  $F$  iff   $\left|\Phi_k\right| = k_F$ for $k = 0,1,..., n.$ . We shall also use the expression - $F$-graded poset.
\end{defn}

\begin{defn}
\noindent  Let  $n\in N \cup \left\{0\right\}\cup \left\{\infty\right\}$. Let   $r,s \in N \cup \left\{0\right\}$.  Let  $\Pi_n$ be the graded partial ordered set (poset) i.e. $\Pi_n = (\Phi_n,\leq)= ( \bigcup_{k=0}^n \Phi_k ,\leq)$ and $\left\langle \Phi_k \right\rangle_{k=0}^n$ constitutes ordered partition of $\Pi_n$. A graded poset   $\Pi_n$  with finite set of minimal 
elements is called \textbf{cobweb poset} \textsl{iff}  
$$\forall x,y \in \Phi \  i.e. \  x \in \Phi_r \ and \  y \in \Phi_s \   r \neq s\ \Rightarrow \   x < y   \ or \ y < x  , $$ 
 $\Pi_\infty \equiv \Pi. $
\end{defn}

\noindent \textbf{Note}. By definition of $\Pi_n$ being graded its  levels $\Phi_r \in \left\{\Phi_k\right\}_k^n$ are independent sets, $n\in N \cup \left\{0\right\}\cup \left\{\infty\right\}$.

\vspace{0.2cm}

\noindent The Definition 2  is the reason for calling Hasse digraph $D = \left\langle \Phi, \leq \cdot \right\rangle $ of the poset $\Pi = (\Phi,\leq))$ a \textbf{\textcolor{red}{Ko}}DAG as in  Professor \textbf{\textcolor{red}{K}}azimierz   \textbf{\textcolor{red}{K}}uratowski native language one word \textbf{\textcolor{red}{Ko}mplet} means \textbf{complete ensemble}- see more in  [2,3]
and for the history of this name see:  The Internet Gian-Carlo Polish Seminar Subject 1. \textit{ oDAGs and KoDAGs in Company }(Dec. 2008). Examples - see Fig.1 and Fig.3. and consult also 
[2,3,4,5] and references therein.

\begin{defn}

$$\langle\Phi_k \to \Phi_n \rangle \equiv \oplus_{n\geq s\geq k} \Phi_s $$
is called the layer of the cobweb poset [2,3,4,5.6.7.8.9.10.11] where $\oplus$ denotes ordinal sum of posets while $\Phi_k$ stay for independent sets (see below) hence 
$$\Pi = \equiv \oplus_{n\geq 0} \Phi_s $$.
\end{defn}

\begin{figure}[ht]
\begin{center}
	\includegraphics[width=100mm]{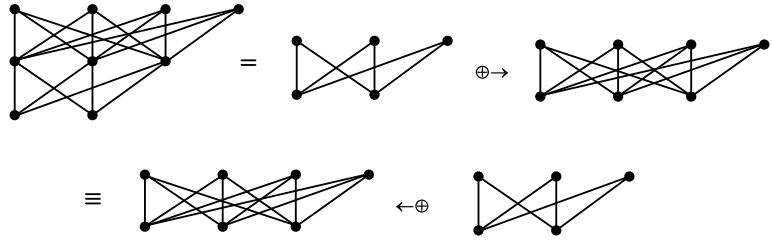}
	\caption 
	{Display of the natural join of bipartite layers $\left\langle \Phi_k \rightarrow \Phi_{k+1} \right\rangle$ $F = N $.} \label{fig:representation}
\end{center}
\end{figure}
	

\noindent \textbf{Comment 1.} Colligate and make \textbf{identifications} of graded DAGs\textbf{ with  $n$-ary relations} digraph representation as in [2,3,4,5]:

$$ \leq = \Phi_0\times\Phi_1\times ... \times\Phi_n \Longleftrightarrow cobweb \  poset  \Longleftrightarrow  KoDAG ,$$
for the natural join of di-bicliques and similarly for $\leq$ being  natural join  of any sequence binary relations [see Fig.1,3]

$$ \leq \subseteq \Phi_0\times\Phi_1\times ... \times\Phi_n \Longleftrightarrow cobweb \  poset  \Longleftrightarrow  KoDAG .$$

\noindent \textcolor{red}{\textbf{Warning.}} Note that \textcolor{red}{\textbf{not for all}} $F$-graded posets their partial orders may be consequently identified with 
$n$-ary relations, where  $F = \left\langle k_F \right\rangle_{k=1}^n $ while $n \in N \cup \left\{\infty\right\}$. This is possible iff
no  biadjacency matrices entering the natural join for $\leq$  has a zero column or a zero row. If  a vertex $m \in \Phi_k$ has not 
either incoming or outgoing arcs then we shall call it the \textcolor{red}{\textbf{mute}} node [2]. This naming being adopted we may say now:\\

\vspace{0.1cm}

\noindent $F$-\textit{graded poset may be identified with} $n$-\textit{ary relation  as above iff it is} $F$-\textit{graded poset with no mute nodes.}\\ 
Equivalently - zero columns or rows in bi-adjacency matrices are forbidden.  See and compare with figures below.




\begin{figure}[ht]
\begin{center}
	\includegraphics[width=100mm]{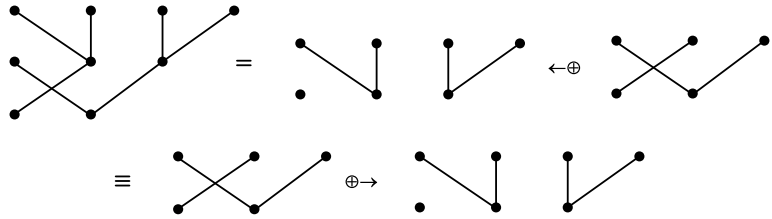}
	\caption {Display of the  natural join $\os$ of  bipartite digraphs with one mute node. }\end{center}\label{fig:representation}
\end{figure}

\begin{figure}[ht]
\begin{center}
	\includegraphics[width=100mm]{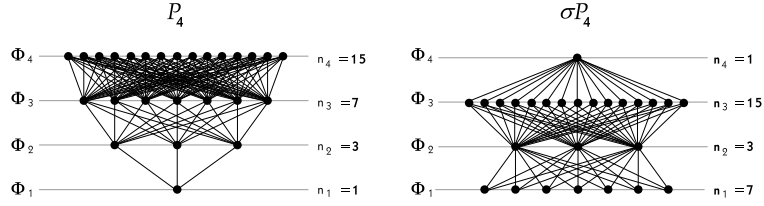}
	\caption {Display of the layer  $\langle\Phi_1 \to \Phi_4 \rangle$ = the subposet $\Pi_4$ of the  $F$ = Gaussian integers sequence $(q=2)$ $F$-cobweb poset and $\sigma \Pi_4$ subposet of the $\sigma$ permuted Gaussian $(q=2)$ $F$-cobweb poset.}\label{fig:representation}
\end{center}
\end{figure}



\noindent \textbf{1.2.} \textbf{Ordinal sum and natural join. } 

\vspace{0.2cm}

\noindent  Let us recall that the ordinal sum [linear sum] of two disjoint ordered sets $P$ and $Q$, denoted by  $P\oplus Q$, is the union of $P$ and $Q$, with $P$'s elements ordered as in $P$ while $Q$'s elements are correspondingly ordered as in $Q$, and  for each $x\in P$ and $y \in Q$ we put  $x \leq y$ .  
The Hasse diagram of  $P\oplus Q$ we construct placing $Q$'s diagram just above $P$'s diagram and with an edge between each minimal element of $Q$ and each maximal element of $P$.  

\vspace{0.2cm}

\noindent Here we propose to add to  the standard operations   $1,2,3$  on ordered sets below from [1] (definitions and properties see: [1] Chapter  III)i.e. operations:\\

\vspace{0.1cm}

\noindent 1.     dual  $P^*$  of  $P$ ;\\

\vspace{0.1cm}

\noindent 2.     the disjoint union  $P+Q =$  cardinal sum ;\\          

\vspace{0.1cm}

\noindent 3.     the ordinal sum   $P\oplus Q$ .\\

\vspace{0.1cm}

\noindent Namely - we propose to \textbf{add} the following binary operation on \textbf{graded} posets:

\vspace{0.1cm}

\noindent \textbf{4.}  this is  \textcolor{red} {\textbf{natural join}} [5,4,3,2] :    $P \os Q$, which here now is   expressed by the ordinal sum  $P\oplus Q$  below.

\vspace{0.1cm}

\noindent \textbf{5.}   Let   $P =\Pi_1 \oplus \Pi_2$  and  $Q = \Pi_2 \oplus \Pi_3$ then  we define  $\os$ via identity

\vspace{0.1cm}

\begin{center}
$P \os Q  \equiv \Pi_1 \oplus \Pi_2 \oplus \Pi_3 $
 \end{center}

\vspace{0.2cm}

\noindent 6.   Therefore   the  cobweb  poset  identified till now with  the natural join of chain of bipartite complete digraphs  ["di-biqliques"] $B_k = \Phi_k \oplus \Phi_{k+1}$ - see Fig.1. -
(contact [2,3,4,5] and references therein) is now defined equivalently as the  ordinal sum  [linear sum ] of chain of trivially ordered sets.\\
Cobweb  poset  $\Pi$  is a  linear sum of trivially ordered sets   $ \left\{\Phi_k \right\}_{k\geq0}$. Nethertheless, the definition of  cobweb poset as natural join of relations {binary, ternary, etc.) as in [5,4,3,2] - by no means - provides additionally   advantages  - visual, based on sight interpretation included. \textit{All graded posets with no mute vertices} in their Hasse diagrams
(i.e. no vertex has  in-degree or out-degree equal zero - consult Fig.2 and compare with Fig.1) are natural join of chain of relations and may be at the same time  interpreted an $n-ary$, $n \in N \cup \left\{\infty \right\}$ relation. As for cobweb posets - these have also discrete hyper-boxes representation [11,2,49,16]. 

\vspace{0.1cm}

\noindent Illustration  -  see Fig.4, Fig.8, Fig.7, Fig.6.

\vspace{0.1cm}

\begin{center}

 $\Pi = \oplus_{k\geq 0} \Phi_k$
 \end{center}
 
 \vspace{0.1cm}
 
\noindent $\left\{\Phi_k \right\}_{k\geq0}$ are then independent sets of $\Pi$. See Fig.4. and Fig.5.

\vspace{0.1cm}



\begin{figure}[ht]
\begin{center}
	\includegraphics[width=100mm]{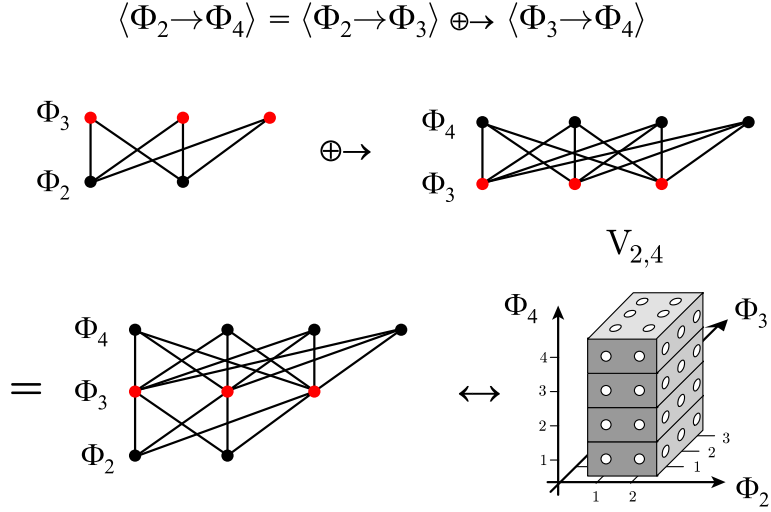}
	\caption 	{Display of the natural join of bipartite layers $\left\langle \Phi_k \rightarrow \Phi_{k+1} \right\rangle$ $F = N $, resulting in  $2\cdot3\cdot4$ maximal chains and equivalent hyper-box $V_{2,4}$ with $2\cdot3\cdot4$ white circle-dots.}\label{fig:representation} 
\end{center}
\end{figure}

\begin{figure}[ht]
\begin{center}
	\includegraphics[width=100mm]{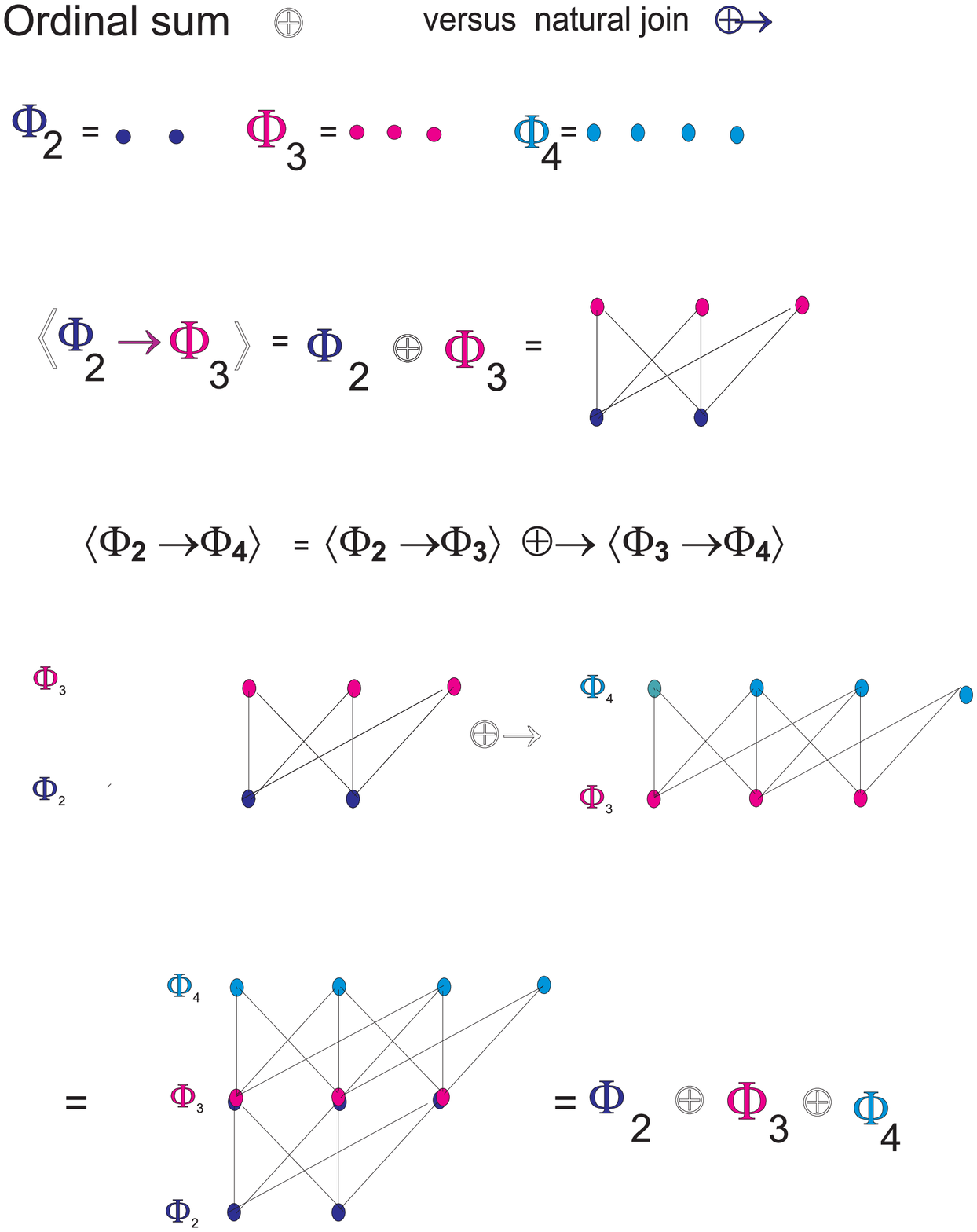}
	\caption {Display of the ordinal sum versus natural join for $F=N$.}\label{fig:representation} 
\end{center}
\end{figure}

\begin{figure}[ht]
\begin{center}
	\includegraphics[width=70mm]{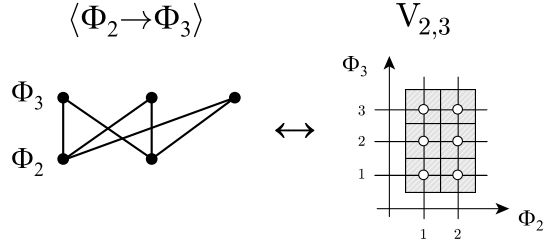}
	\caption{Bipartite  layer $\langle\Phi_3 \rightarrow \Phi_4 \rangle$ with six maximal chains and equivalent hyper-box $V_{2,3}$  with six white circle-dots} \label{fig:representation}
\end{center}
\end{figure}

\begin{figure}[ht]
\begin{center}
	\includegraphics[width=70mm]{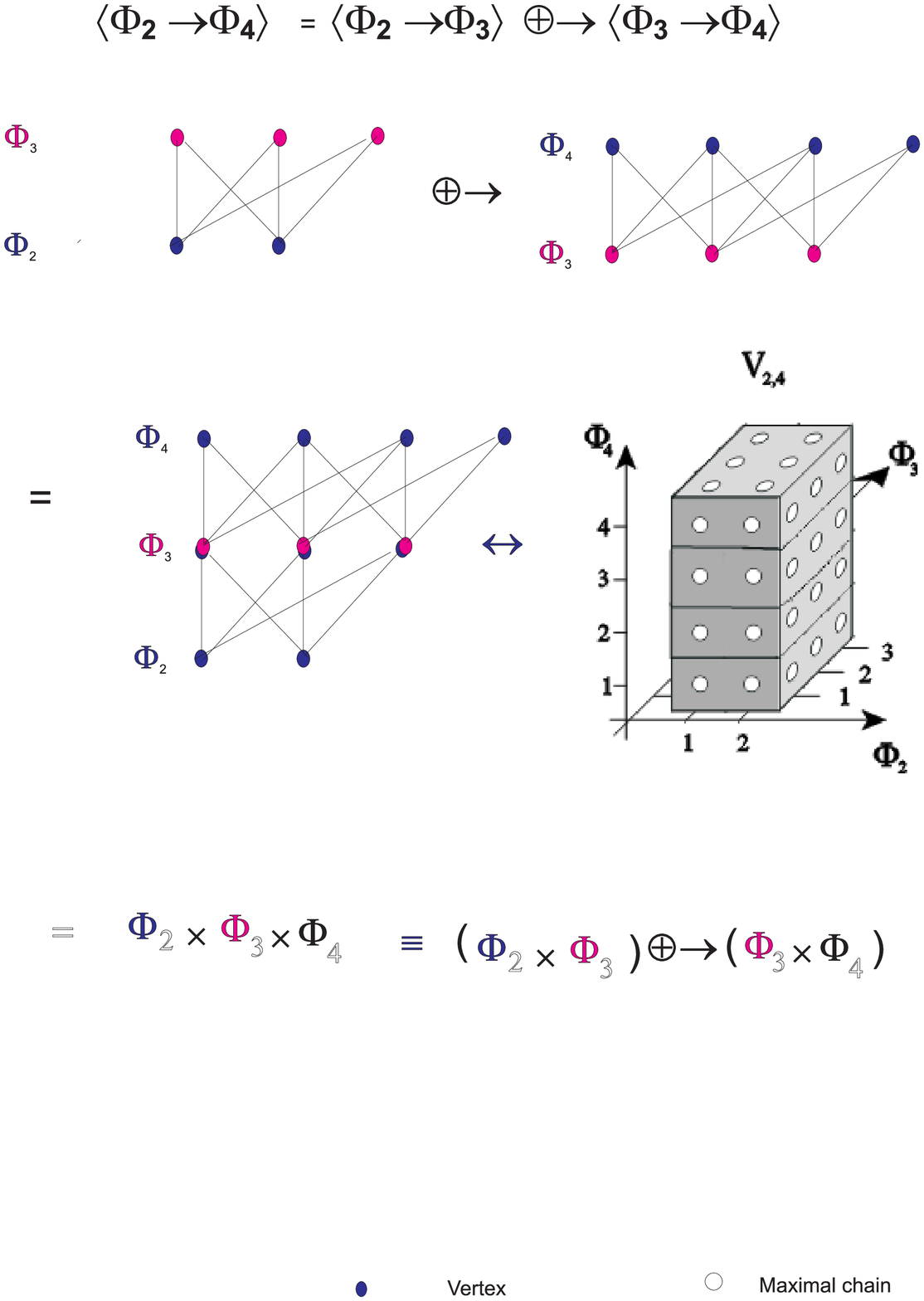}
	\caption{Natural join of layers versus Cartesian product constituting the hyper-box $V_{2,4}$.} \label{fig:representation}
\end{center}
\end{figure}

\begin{figure}[ht]
\begin{center}
	\includegraphics[width=70mm]{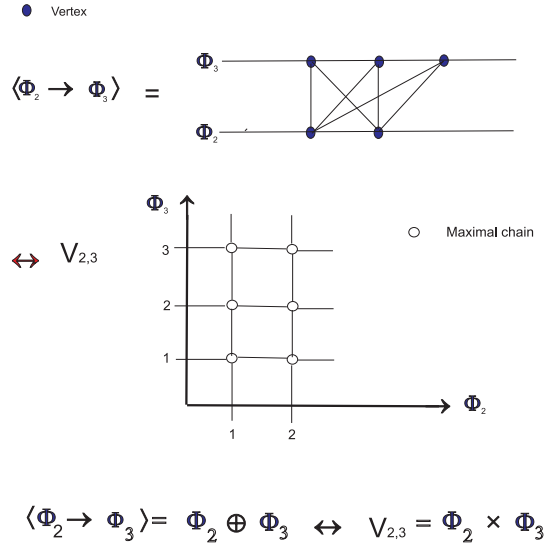}
	\caption{Natural join as ordinal sum of levels $\Phi_2\os \Phi_3 \equiv  \Phi_2\oplus \Phi_3 \equiv V_{2,3}$.}\label{fig:representation}
\end{center}
\end{figure}

\begin{figure}[ht]
\begin{center}
	\includegraphics[width=70mm]{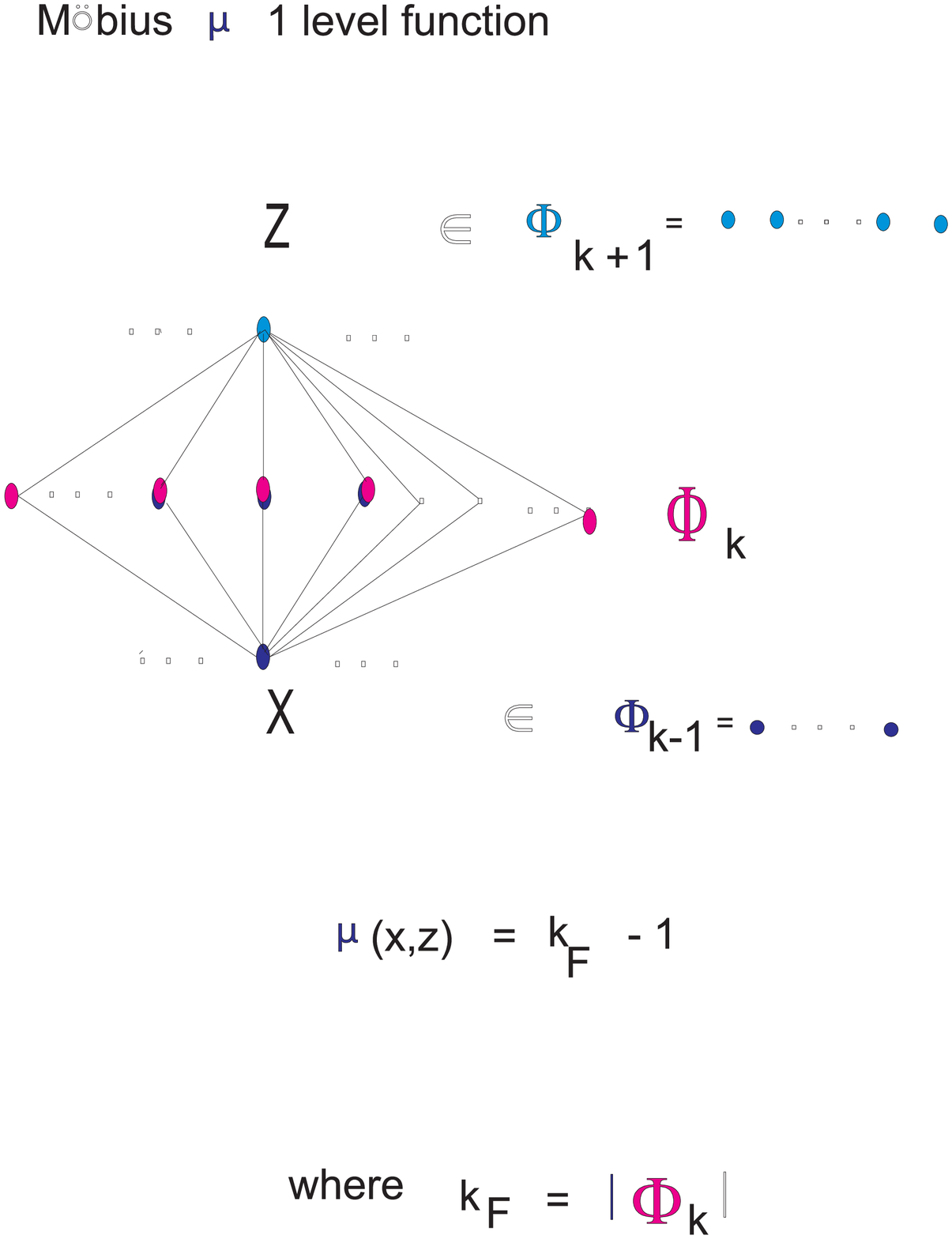}
	\caption{M\"{o}bius 1 level function for  $[x,z] = x\oplus \Phi_k\oplus z $ , $x\in \Phi_{k-1}, \ z \in \Phi_{k+1}$.}\label{fig:representation}
\end{center}
\end{figure}

\begin{figure}[ht]
\begin{center}
	\includegraphics[width=70mm]{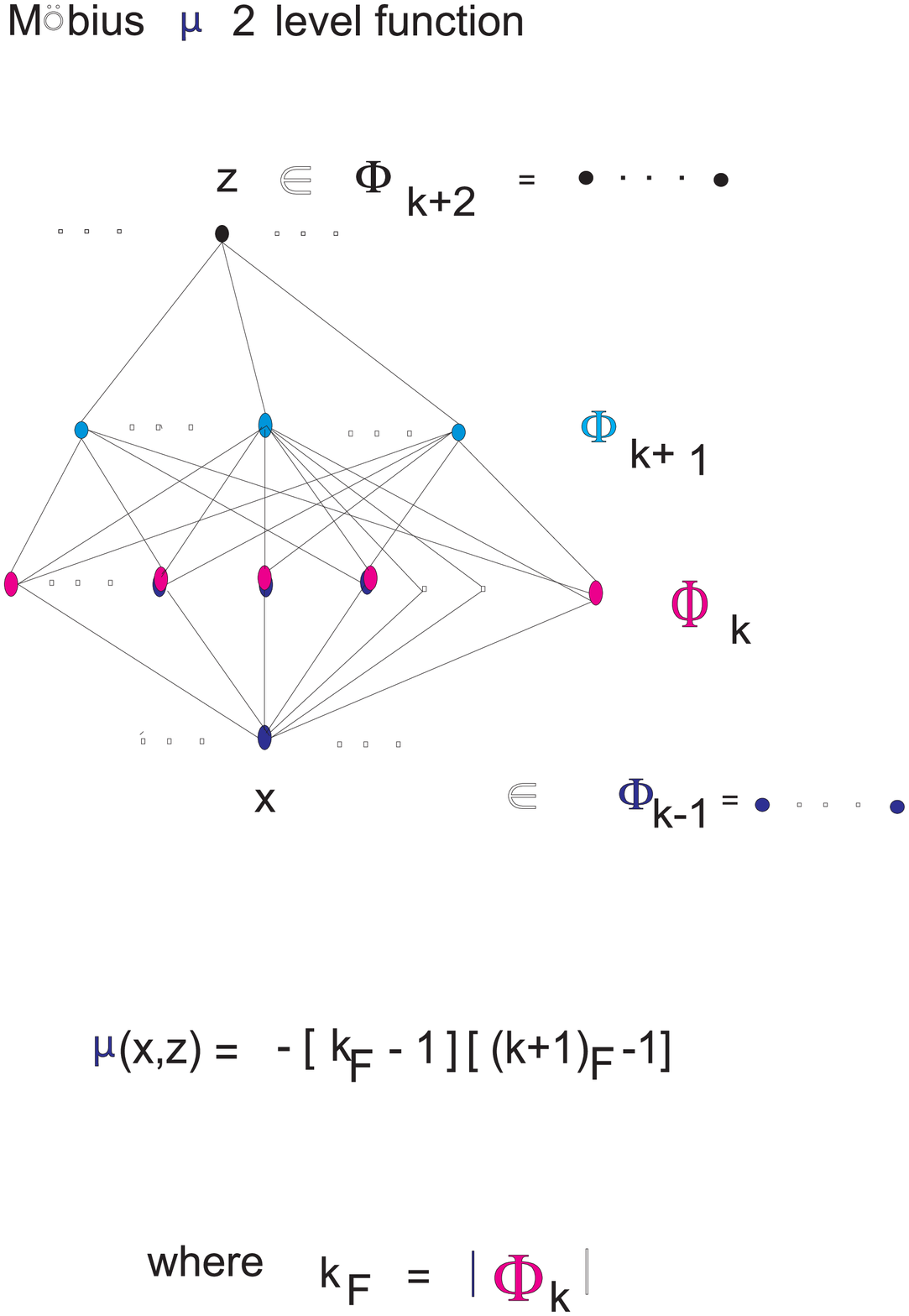}
	\caption{M\"{o}bius 2 level function  for  $[x,z] = x\oplus \Phi_k\oplus \Phi_{k+1}\oplus z $ , $x\in \Phi_{k-1}, \ z \in \Phi_{k+2}$.}\label{fig:representation}
\end{center}
\end{figure}



\noindent \textbf{1.3. Colligation of the above with hyper-boxes from [11]} \\

\noindent  Recall [2]: $C_{max}(\Pi_n)$ is the set of all maximal chains of $\Pi_n$. Recall [2]: 

$$C^{k,n}_{max} = \big\{ \mathrm{maximal\ chains\ in\ } \langle \Phi_k \rightarrow \Phi_n \rangle \big\}.$$

\noindent Consult  Section 3. in  [11] in order to view $C_{max}(\Pi_n)$ or  $C^{k,n}_{max}$ as the hyper-box of points.\\ 
Namely [11,2] denoting with $V_{k,n}$ the discrete finite rectangular $F$-hyper-box or $(k,n)-F$-hyper-box or in everyday parlance just $(k,n)$-box

$$
	V_{k,n} = [k_F]\times [(k+1)_F]\times ... \times[n_F]
$$
\noindent we identify (see Figure 7.) the following two just by agreement according to the $F$-natural identification:
$$
	C^{k,n}_{max} \equiv V_{k,n}
$$
 
i.e.

$$
C^{k,n}_{max} = \big\{ \mathrm{maximal\ chains\ in\ } \langle \Phi_k \rightarrow \Phi_n \rangle \big\} \equiv V_{k,n}.
$$


\vspace{0.1cm}

\noindent Illustration  -  see Fig.8, Fig.4, Fig.7, Fig.6 and for more  see [11], [47], [12] and  [2]  with specific indication on [16].

\vspace{0.2cm}

\noindent \textcolor{red}{\textbf{Important}}.  Accordingly  the \textcolor{blue}{\textbf{natural join operation}} of  discrete hyper-boxes - ( cobweb posets are  encoded by discrete hyper-boxes! [11])  is just \textcolor{blue}{\textbf{Cartesian product of them accompanied with projection out}} of sine qua non common faces (see \textcolor{blue}{\textbf{Fig.7}}) which is schematically represented by an sample case below  [note - cobweb poset might be defined also as the ordinal sum of its independence  sets $\left\{\Phi_k\right\}_{k \geq )}$ ],  

\vspace{0.1cm}

$$\left(\Phi_k \times \Phi_{k+1}\right) \os \left(\Phi_{k+1} \times \Phi_{k+2}\right) \equiv  \Phi_k \times \Phi_{k+1} \times \Phi_{k+2} $$

\noindent see Fig.5 and Fig.7.

\vspace{0.2cm}

\noindent \textbf{1.4. Cobweb posets' M\"{o}bius function}\\

\noindent  In order to find out the  M\"{o}bius function $\mu$ of the cobweb poset $\Pi$ note the following obvious statements  ($\left|\Phi_k\right|= k_F$)
  
\vspace{0.1cm}
  
\noindent Obvious statement: $\mu(x,x) = 1, \   \mu(x,y) = - 1$ for $x\prec\cdot y , \ $and$  \   \mu(x,z) = k_F - 1 $  for   $[x,z] = x\oplus \Phi_k\oplus z $.

\vspace{0.1cm}

\noindent Obvious statement: ($\left|\Phi_k\right|= k_F$)  for  $[x,z] = x\oplus \Phi_k\oplus z $ , $x\in \Phi_{k-1}, \ z \in \Phi_{k+1}$

\vspace{0.1cm}

\noindent Obvious statement: $$\mu(x,z) = [ k_F - 1].$$ 

\vspace{0.1cm}

\noindent Obvious statement ($\left|\Phi_k\right|= k_F$)  :  for  $[x,z] = x\oplus \Phi_k\oplus \Phi_{k+1}\oplus z \equiv x\oplus B_k \oplus z$ , $x\in \Phi_{k-1}, \ z \in \Phi_{k+2}$ 

\vspace{0.1cm}
                                                                     
$$\mu(x,z) = - [ k_F - 1][ (k+1)_F - 1].$$

\vspace{0.1cm}

\noindent Illustration  ( see Fig.9 and Fig.10).
                                                                     
\vspace{0.1cm}

\noindent Hence via induction  for  $x \in \Phi_r$ ,   $z \in \Phi_s$  and  for     
$$s>r , \ [x,z] = x\oplus \Phi_{r+1}\oplus...\oplus \Phi_{s-1}\oplus z$$ 
we have the final obvious statement:

  $$\mu(x,y) = (-1)^{s-r}\prod_{k=r+1}^{s-1}[k_F-1]. $$


\vspace{0.1cm} 

 \noindent Compare the above  formula as derived and discussed in [2] and declared in [46, 24] for Fibonacci sequence. Compare also with [25-28].

\vspace{0.2cm}

\noindent Naturally the values of $\mu(x,y)$ depend only on the rank of its arguments - here $r(x)=r$ and $r(y)=s$ - which is the reason of coding matrix existence 
for M\"{o}bius function in a natural labelling representation ( see: examples in Section 4 and Definition 7). The rank function is here defined as follows:
$r(x)= r$  if  $x\in \Phi_r$.



\section{Combinatorial interpretation.}

\noindent For \textbf{combinatorial interpretation of cobweb posets} via their cover relation digraphs (Hasse diagrams) called KoDAGs see [7,6,2,3,5]. The recent equivalent formulation 
of this combinatorial interpretation is to be found in [6] (Feb 2009) or [8] from which we quote it here down.

\begin{defn}
$F$-\textbf{nomial} \textbf{coefficients} are defined as follows
$$
	\fnomial{n}{k} = \frac{n_F!}{k_F!(n-k)_F!} 
	= \frac{n_F\cdot(n-1)_F\cdot ...\cdot(n-k+1)_F}{1_F\cdot 2_F\cdot ... \cdot k_F}
	= \frac{n^{\underline{k}}_F}{k_F!}
$$
\noindent while $n,k\in \mathbb{N}$ and $0_F! = n^{\underline{0}}_F = 1$  with $n^{\underline{k}}_F \equiv \frac{n_F!}{k_F!}$ staying for falling factorial.
$F$ is called  $F$-graded poset  \textcolor{blue}{\textbf{admissible}} sequence iff  $\fnomial{n}{k} \in N \cup\left\{0\right\}$ ( In particular we shall use the expression - $F$-cobweb admissible sequence).
\end{defn}

\begin{defn}

$$
C_{max}(\Pi_n) \equiv  \left\{c=<x_0,x_1,...,x_n>, \: x_s \in \Phi_s, \:s=0,...,n \right\} 
$$  
i.e. $C_{max}(\Pi_n)$ is the set of all maximal chains of $\Pi_n$

\end{defn}

\noindent and consequently (see Section 2 in [11]  on Cobweb posets' coding via $N^\infty$ lattice boxes)

\begin{defn} ($C^{k,n}_{max}$) Let  

$$
C_{max}\langle\Phi_k \to \Phi_n \rangle \equiv \left\{c=<x_k,x_{k+1},...,x_n>, \: x_s \in \Phi_s, \:s=k,...,n \right\}\equiv
$$

$$
	\equiv \big\{ \mathrm{maximal\ chains\ in\ } \langle \Phi_k \rightarrow \Phi_n \rangle \big\} \equiv
	C_{max}\big( \langle \Phi_k \rightarrow \Phi_n \rangle \big) \equiv
	C^{k,n}_{max}.
$$
\end{defn}

\noindent \textbf{Note.} The $C_{max}\langle\Phi_k \to \Phi_n \rangle \equiv C^{k,n}_{max}$
is the hyper-box points'  set [11] of  Hasse sub-diagram corresponding maximal chains and it defines biunivoquely 
the layer $\langle\Phi_k \to \Phi_n \rangle = \bigcup_{s=k}^n\Phi_s$  as the set of maximal chains' nodes (and vice versa) -
for  these arbitrary $F$-denominated \textbf{graded} DAGs (KoDAGs included).

\vspace{0.1cm}

\noindent The equivalent to that of [7,6,8] formulation of the fractals reminiscent combinatorial interpretation of cobweb posets via their cover relation digraphs (Hasse diagrams) is the following.

\vspace{0.2cm}

\noindent \textbf{Theorem 1} [10,8,7,2] \\
\noindent(Kwa\'sniewski) \textit{For $F$-cobweb admissible sequences $F$-nomial coefficient $\fnomial{n}{k}$ is the cardinality of the family of \emph{equipotent} to  $C_{max}(P_m)$ mutually disjoint maximal chains sets, all together \textbf{partitioning } the set of maximal chains  $C_{max}\langle\Phi_{k+1} \to \Phi_n \rangle$  of the layer   $\langle\Phi_{k+1} \to \Phi_n \rangle$, where $m=n-k$.}

\vspace{0.1cm} \noindent For environment needed and then  simple combinatorial proof see [7,8,3,4,5]  easily accessible via Arxiv.\\

\vspace{0.1cm}

\noindent \textbf{Comment 2}. For the  above Kwa\'sniewski combinatorial  interpretation of  $F$-nomials' array \textit{it does not matter}  of course whether the diagram is being directed  or not, as this combinatorial interpretation is  equally valid for partitions  of the family of  $SimplePath_{max}(\Phi_k - \Phi_n)$ in  comparability graph of the Hasse  digraph with self-explanatory notation used on the way. The other insight into this irrelevance for combinatoric interpretation is [9]: colligate the coding of $C^{k,n}_{max}$ by hyper-boxes. (More on that soon).  And to this end recall  what really also matters here : a poset is graded if and only if every connected component of its \textbf{comparability graph} is graded. We are concerned here with connected graded graphs and digraphs.

\vspace{0.2cm}

\noindent For the relevant recent developments see [9]  while \textbf{[10]} is their all \textbf{source paper} as well as those reporting on the broader affiliated research (see [11-22,24-28,46-48] and references therein). The inspiration for  "`philosophy"'  of notation in mathematics as that in Knuth's from [23] - in the case of "`upside-downs"'  has been  driven by Gauss "`$q$-Natural numbers"'$\equiv N_q = \left\{n_q=q^0+q^1+...+q^{n-1}\right\}_{n\geq 0}$ from finite geometries of  linear subspaces lattices over Galois fields. As for the earlier use and origins of the use of this author's upside down notation see [29-45].

\vspace{0.1cm}

\noindent In discrete hyper-boxes language the combinatorial interpretation reads: 

\vspace{0.2cm}

\noindent \textbf{Theorem 2}\\
\noindent(Kwa\'sniewski) \textit{For $F$-cobweb admissible sequences $F$-nomial coefficient $\fnomial{n}{k}$ is the cardinality of the family of \emph{equipotent} to  $V_{0,m}$ mutually disjoint 
discrete hyper-boxes, all together \textbf{partitioning } the discrete hyper-box   $V_{k+1,n}$    $\equiv$  to the layer   $\langle\Phi_{k+1} \to \Phi_n \rangle$, where $m=n-k$.}

\vspace{0.2cm}

\noindent \textbf{Comment 3.} Recall:  all graded posets with no mute vertices in their Hasse diagrams [2] (i.e. no vertex has  in-degree or out-degree equal zero) are natural join of chain of relations and may be at the same time  interpreted an $n-ary$, $n \in N \cup \left\{\infty \right\}$ relation.\\
\noindent Colligate any binary relation $R$ with Hasse digraph cover relation  $\prec\cdot$  and identify as in [3,2] $\zeta({\bf R})\equiv {\bf
R}^{*}$ with incidence algebra zeta function and with zeta matrix of the poset associated to its Hasse digraph, where the \textbf{reflexive} reachability  relation $\zeta({\bf R})\equiv {\bf
R}^{*}$ is defined as

$${\bf R^*}=R^0\cup R^1\cup R^2\cup\ldots\cup R^n\cup\ldots  
\bigcup_{k>0}{R}^k={\bf R}^{\infty}\cup {\bf I}_A =$$
 
\begin{center}
=  {\bf transitive} and {\bf reflexive} closure of $\bf{R}$  $\Leftrightarrow$
\end{center}
 
$$
\Leftrightarrow\;  A(R^{\infty})=A({ R})^{\copyright 0}\vee A({
R})^{\copyright 1}\vee A({ R})^{\copyright 2}\vee \ldots \vee A(R)^{\copyright n}\vee \ldots ,
$$

\noindent where $A({\bf R})$ is the Boolean  adjacency matrix of the
relation ${\bf R}$ simple digraph  and  $\copyright$ stays for Boolean product.

\vspace{0.2cm}

\noindent Then colligate and/or recall from [3] the resulting schemes.
\noindent {\bf Schemes:}
 $$ <=\prec\cdot ^{\infty}=\;connectivity\;of\;\prec\cdot$$
 $$ \leq=\prec\cdot ^{*}=\;reflexive \; reachability\;of\;\prec\cdot$$
 $$\prec\cdot ^{*}=\zeta(\prec\cdot).$$



\vspace{0.2cm}

\noindent \textcolor{red}{\textbf{Remark 1. Obvious.}} Needed also for the next Section. Compare with the Observation 3. below.\\ 
\vspace{0.1cm}

\noindent The  $\zeta$ matrix ($\equiv$ the algebra structure coding element of the incidence algebra $I(P,\leq)$ is the characteristic function $ \chi$ of a partial order relation $\leq$\textbf{ for any} given $F$ - graded poset including  $F$- cobweb posets  $\Pi$ :

$$  \zeta  =    \chi\left( \leq \right). $$

\noindent The consequent  (customary-like notation included) notation of other  algebra  $I(P,\leq)$ important  elements then - for the any fixed order $\leq $   -  is the following [3,2,1]:  
$$  \zeta  =    \chi\left( \leq \right)= \zeta_< + \delta, $$
$$  \zeta_<  =    \chi\left( <  \right) =  \zeta - \delta \equiv \rho ,\: (reachability=connectivity),$$

$$  \zeta_{\prec\cdot}  =    \chi\left( \prec\cdot \right) \equiv \kappa, \:(cover), $$

$$  \zeta_{\leq \cdot}  =    \chi\left( \leq \cdot \right)= \kappa +\delta \equiv \eta, \:(reflexive "`cover"'). $$

$$
\eta	= \kappa +\delta  = \left[ \begin{array}{llllll}
I_1 &	B_1 & zeros \\
& I_2 &	B_2 & zeros \\
& & I_3 &	B_3  & zeros \\
	& & ... & \\
	 & & & I_n & B_n & zeros
	\end{array} \right]
$$

$$
\eta^{-1}	= \left(\delta + \kappa\right)^{-1} = \sum_{k\geq 0} (-\kappa)^k  = \left[ \begin{array}{llllll}
I_1 &	-B_1 & B_1B_2 & ... \\
&  I_2 & -B_2  & ... \\
& & I_3 & -B_3 & ... \\
	& & ... & \\
	& & &  I_n & -B_n & ... 
	\end{array} \right]
$$
\noindent Recall from [3] : $B(A)$ is the   biadjacency i.e cover relation $\prec\cdot$  matrix of the adjacency matrix $A$.\\
\textcolor{red}{\textbf{Note}}:  biadjacency and  cover relation $\prec\cdot$  matrix for bipartite digraphs coincide. By extension - we call  \textbf{cover relation} $\prec\cdot$ matrix $\kappa$ the bi-adjacency matrix too in order to keep reminiscent convocations going on.

\vspace{0.1cm}
\noindent As a consequence  - quoting [3] - we have:

$$
	B\left(\os_{i=1}^n G_i \right) \equiv  B [\os_{i=1}^n  A(G_i)]  =  \oplus_{i=1}^n  B[A(G_i) ]  \equiv  \mathrm{diag} (B_1 , B_2 , ..., B_n) = 
$$
$$
	= \left[ \begin{array}{lllll}
	B_1 \\
	& B_2 \\
	& & B_3 \\
	& ... & ... & ...\\
	& & & & B_n
	\end{array} \right],
$$

\noindent or equivalently 

$$
	\kappa \equiv \chi \left(\os_{i=1}^n \prec\cdot_i \right)  \equiv 
$$
$$
	\equiv \left[ \begin{array}{llllll}
	0 & B_1 \\
	& 0 & B_2 \\
	& & 0 & B_3 \\
	& ... & ... & ...\\
	& & & & 0 & B_n\\
	 & & & & & 0 

	\end{array} \right],
$$

\noindent $n \in N \cup \{\infty\}$

\vspace{0.2cm}

\noindent  In view of the all above the following is obvious;

$$\left(A \os B \right)^{-1} \neq A^{-1}\os B^{-1} $$
except for the trivial case.


\vspace{0.2cm}

\noindent  Anticipating considerations of Section III and   customarily  allowing for the identifications:  $ \chi \left(\prec\cdot\right) \equiv \prec\cdot  \equiv  \kappa $ 
- consider $[Max] \in I(P,R)$: 
$$[Max] = {(I - \prec\cdot)}^{-1}= \prec\cdot^0 + \prec\cdot^1 + \prec\cdot^2 +...+ \prec\cdot^k + ... = \sum_{k\geq 0} \kappa^k$$
in order to note that ( $x \in \Phi_t  \equiv  x = x_t \in \Phi_t$)

\begin{center}
$[Max]_{s,t}$ = the number of all maximal chains in the poset interval $[x_s,x_t] \equiv  [s,t].$
\end{center}

\vspace{0.1cm}

\noindent where $x_s \in \Phi_s$  and $x_t \in \Phi_t$ for , say , $s \leq t$ with the reflexivity (loop) convention adopted i.e. $[Max]_{t,t}=1$.

\vspace{0.2cm}

\noindent \textcolor{blue}{\textbf{Sub-Remark 1.1.}} It is now a good - prepared for - place to note further relevant properties of constructs as to be used in the sequel. 
These are the following.

$$C_{max}\left(\langle\Phi_r \to \Phi_k \rangle \os \langle\Phi_k \to \Phi_s \rangle \right) = C_{max}\langle\Phi_r\to \Phi_s \rangle,$$
for $r\leq k \leq s$ while  $\left| \Phi_n \right| \equiv n_F$.

\vspace{0.2cm}

\noindent  Let $ \left|C^{k,n}_{max}\right|\equiv C^{k,n}$.  Then for $F$-cobweb posets (what about just $F$-graded?) we note that

$$ C^{r,k}C^{k,s} = k_F C^{r,s},$$ hence  $$C^{r,k}C^{k,s} =  C^{r,s} \;iff \; k_F =1 $$ 
for $r\leq k \leq s$ while 

$$
C^{r,k}C^{k+1,s} = C^{r,s}
$$
for $r\leq k < s$ while  $\left| \Phi_n \right| \equiv n_F$. Let us now see in more detail  how this  kind (Q.M.?) of  mimics of \textbf{\textit{Markov}} \textbf{property} is intrinsic for natural joins of digraphs. For that to do  consider levels i.e. independent (stable)  sets $\Phi_k = \left\{x_{k,i}\right\}^{k_F}_{i=1} $ and extend the notation accordingly so as to encompass 

$$
\langle \Phi_r \to x_{k,i} \rangle = \left\{c=<x_r,x_{r+1},...,x_{k-1},x_{k,i} >, \: x_s \in \Phi_s, \:s=r,...,k-1 \right\}.
$$
Let

$$
\left|C_{max}\langle \Phi_r \to x_{k,i} \rangle\right| \equiv C^{r,k,i}.
$$ 
Then
$$
\sum_{i=1}^{k_F} C^{r,k,i}C^{s,k,i}= C^{r,s}
$$
for $r\leq k < s$, ...  (for $r\leq k \leq s$ ?). In the case of cobweb posets (what about just $F$-graded?) the numbers $C^{s,k,i}$ are the same for each $ i = 1,...,k_F$ therefore we   
have for cobwebs

$$k_F C^{r,k,i}C^{s,k,i}= C^{r,s}$$
which in view of $k_F C^{r,k,i}= C^{r,k}$ is of course consistent with $ C^{r,k}C^{k,s} = k_F C^{r,s}.$ We consequently notice  that  - with self-evident extension of notation:

$$
\langle\Phi_k\to \Phi_n \rangle =   \bigcup_{i,j=1}^{k_F,n_F} \langle x_{k,i}\to x_{n,j} \rangle.
$$

\vspace{0.1cm}

\noindent The frequently used block matrices are: 1) $I (s\times k)$ which denotes $(s\times k)$  matrix  of  ones  i.e.  $[ I (s\times k) ]_{ij} = 1$;  $1 \leq i \leq  s,  1\leq j  \leq k.$  and  $n \in N \cup \{\infty\}$, 2) and   $B(s\times k)$ which stays for  $(s\times k)$  matrix  of  ones  and zeros accordingly to the $F$-graded poset has been fixed - see Observation 2.

\noindent In the block matrices language the above \textbf{\textit{Markov}} \textbf{property} for cobweb posets (what about just $F$-graded?) reads as follows (to be used in Section 3)  for example :
$$
I\left(r_F\times (r+1)_F\right) I\left( (r+1)_F\times (r+2)_F\right) = (r+1)_F I\left(r_F\times (r+2)_F\right).
$$
Well, what about then just $F$-graded? - See  Comment 3 and its Warning.



\section{ Zeta  and  M{\"{o}}bius  funcions formulas and Matrices in  natural labeling with examples}

\noindent \textbf{Examples} of $\zeta({\bf \leq}) \in I(\Pi,Z)$\\

\vspace{0.1cm}
  
\noindent Let $F$ denotes arbitrary natural numbers valued sequence.  Let  $A_N$ be the Hasse matrix i.e. adjacency matrix of cover relation $\prec\cdot$ digraph denominated by sequence $N$ [1]. Then the zeta matrix  $\zeta = (1-\mathbf{A}_N)^{-1\copyright}$ for the denominated by  $F = N$ cobweb poset is of the form [3]  (see also [17-22,6]):

\vspace{1mm}

$$ \left[\begin{array}{ccccccccccccccccc}
\textbf{1 }& 1 & 1 & 1 & 1 & 1 & 1 & 1 & 1 & 1 & 1 & 1 & 1 & 1 & 1 & 1 & \cdots\\
0 & \textbf{1} & \textbf{\textcolor{blue}{0}} & 1 & 1 & 1 & 1 & 1 & 1 & 1 & 1 & 1 & 1 & 1 & 1 & 1 & \cdots\\
0 & \textbf{0} & \textbf{1} & 1 & 1 & 1 & 1 & 1 & 1 & 1 & 1 & 1 & 1 & 1 & 1 & 1 & \cdots\\
0 & 0 & 0 & \textbf{1} & \textbf{\textcolor{blue}{0}} & \textbf{\textcolor{blue}{0}} & 1 & 1 & 1 & 1 & 1 & 1 & 1 & 1 & 1 & 1 & \cdots\\
0 & 0 & 0 & \textbf{0} & \textbf{1} & \textbf{\textcolor{blue}{0}} & 1 & 1 & 1 & 1 & 1 & 1 & 1 & 1 & 1 & 1 & \cdots\\
0 & 0 & 0 & \textbf{0} & \textbf{0} &\textbf{ 1 }& 1 & 1 & 1 & 1 & 1 & 1 & 1 & 1 & 1 & 1 & \cdots\\
0 & 0 & 0 & 0 & 0 & 0 & \textbf{1} & \textbf{\textcolor{blue}{0}}& \textbf{\textcolor{blue}{0}}&\textbf{\textcolor{blue}{0}} & 1 & 1 & 1 & 1 & 1 & 1 & \cdots\\
0 & 0 & 0 & 0 & 0 & 0 & \textbf{0} & \textbf{1} & \textbf{\textcolor{blue}{0}} & \textbf{\textcolor{blue}{0}} & 1 & 1 & 1 & 1 & 1 & 1 & \cdots\\
0 & 0 & 0 & 0 & 0 & 0 & \textbf{0} & \textbf{0} & \textbf{1} & \textbf{\textcolor{blue}{0}}& 1 & 1 & 1& 1 & 1 & 1 & \cdots\\
0 & 0 & 0 & 0 & 0 & 0 & \textbf{0} & \textbf{0} & \textbf{0}& \textbf{1} & 1 & 1 & 1 & 1 & 1 & 1 & \cdots\\
0 & 0 & 0 & 0 & 0 & 0 & 0 & 0 & 0 & 0 & \textbf{1} & \textbf{\textcolor{blue}{0}} & \textbf{\textcolor{blue}{0}} & \textbf{\textcolor{blue}{0}}& \textbf{\textcolor{blue}{0}}& 1 & \cdots\\
0 & 0 & 0 & 0 & 0 & 0 & 0 & 0 & 0 & 0 & \textbf{0} & \textbf{1} & \textbf{\textcolor{blue}{0}}& \textbf{\textcolor{blue}{0}} & \textbf{\textcolor{blue}{0}} & 1 & \cdots\\
0 & 0 & 0 & 0 & 0 & 0 & 0 & 0 & 0 & 0 & \textbf{0} & \textbf{0} & \textbf{1} & \textbf{\textcolor{blue}{0}} & \textbf{\textcolor{blue}{0}} & 1 & \cdots\\
0 & 0 & 0 & 0 & 0 & 0 & 0 & 0 & 0 & 0 & \textbf{0} & \textbf{0} & \textbf{0} & \textbf{1} & \textbf{\textcolor{blue}{0}}& 1 & \cdots\\
0 & 0 & 0 & 0 & 0 & 0 & 0 & 0 & 0 & 0 & \textbf{0} & \textbf{0} & \textbf{0} & \textbf{0} & \textbf{1} & 1 & \cdots\\
0 & 0 & 0 & 0 & 0 & 0 & 0 & 0 & 0 & 0 & 0 & 0 & 0 & 0 & 0 & \textbf{1} & \cdots\\
. & . & . & . & . & . & . & . & . & . & . & . & . & . & . & . & . \cdots\\
 \end{array}\right]$$
   
\noindent \textbf{Example.1 $\zeta_N$.  The incidence matrix for the N-cobweb poset.} 

\vspace{0.2cm}

\noindent Note that the  matrix $\zeta$ representing uniquely its corresponding  cobweb poset does  exhibits  a staircase structure of zeros above the diagonal (see above, see below) which is characteristic to Hasse diagrams of \textbf{all} cobweb posets and for graded posets it is characteristic too.

\vspace{2mm}

$$ \left[\begin{array}{ccccccccccccccccc}
\textbf{1} & 1 & 1 & 1 & 1 & 1 & 1 & 1 & 1 & 1 & 1 & 1 & 1 & 1 & 1 & 1 & \cdots\\
0 & \textbf{1} & 1 & 1 & 1 & 1 & 1 & 1 & 1 & 1 & 1 & 1 & 1 & 1 & 1 & 1 & \cdots\\
0 & 0 & \textbf{1} & 1 & 1 & 1 & 1 & 1 & 1 & 1 & 1 & 1 & 1 & 1 & 1 & 1 & \cdots\\
0 & 0 & 0 & \textbf{1} & \textbf{\textcolor{red}{0}} & 1 & 1 & 1 & 1 & 1 & 1 & 1 & 1 & 1 & 1 & 1 & \cdots\\
0 & 0 & 0 & \textbf{0} & \textbf{1} & 1 & 1 & 1 & 1 & 1 & 1 & 1 & 1 & 1 & 1 & 1 & \cdots\\
0 & 0 & 0 & 0 & 0 & \textbf{1} & \textbf{\textcolor{red}{0}} & \textbf{\textcolor{red}{0}} & 1 & 1 & 1 & 1 & 1 & 1 & 1 & 1 & \cdots\\
0 & 0 & 0 & 0 & 0 & \textbf{0} & \textbf{1} & \textbf{\textcolor{red}{0}} & 1 & 1 & 1 & 1 & 1 & 1 & 1 & 1 & \cdots\\
0 & 0 & 0 & 0 & 0 & \textbf{0} & \textbf{0} & \textbf{1} & 1 & 1 & 1 & 1 & 1 & 1 & 1 & 1 & \cdots\\
0 & 0 & 0 & 0 & 0 & 0 & 0 & 0 & \textbf{1} & \textbf{\textcolor{red}{0}} & \textbf{\textcolor{red}{0}} & \textbf{\textcolor{red}{0}} & \textbf{\textcolor{red}{0}} & 1 & 1 & 1 & \cdots\\
0 & 0 & 0 & 0 & 0 & 0 & 0 & 0 & \textbf{0} & \textbf{1} & \textbf{\textcolor{red}{0}}& \textbf{\textcolor{red}{0}} & \textbf{\textcolor{red}{0}} & 1 & 1 & 1 & \cdots\\
0 & 0 & 0 & 0 & 0 & 0 & 0 & 0 & \textbf{0} & \textbf{0} & \textbf{1} & \textbf{\textcolor{red}{0}}& \textbf{\textcolor{red}{0}} & 1 & 1 & 1 & \cdots\\
0 & 0 & 0 & 0 & 0 & 0 & 0 & 0 & \textbf{0} & \textbf{0} & \textbf{0} & \textbf{1} & \textbf{\textcolor{red}{0}}& 1 & 1 & 1 & \cdots\\
0 & 0 & 0 & 0 & 0 & 0 & 0 & 0 & \textbf{0} & \textbf{0} & \textbf{0} & \textbf{0} & \textbf{1} & 1 & 1 & 1 & \cdots\\
0 & 0 & 0 & 0 & 0 & 0 & 0 & 0 & 0 & 0 & 0 & 0 & 0 & \textbf{1} & \textbf{\textcolor{red}{0}} & \textbf{\textcolor{red}{0}} & \cdots\\
0 & 0 & 0 & 0 & 0 & 0 & 0 & 0 & 0 & 0 & 0 & 0 & 0 & \textbf{0} & \textbf{1} & \textbf{\textcolor{red}{0}} & \cdots\\
0 & 0 & 0 & 0 & 0 & 0 & 0 & 0 & 0 & 0 & 0 & 0 & 0 & \textbf{0} & \textbf{0} & \textbf{1} & \cdots\\
. & . & . & . & . & . & . & . & . & . & . & . & . & . & . & . & . \cdots\\
 \end{array}\right]$$

\vspace{1mm} \noindent \textbf{Example.2  $\zeta_F$.  The matrix
$\zeta$ for the Fibonacci cobweb poset associated to \textbf{$F$-KoDAG} Hasse digraph. }

\vspace{2mm}   

\noindent The above remarks are visualized as below  [17-22,6]. Namely  - apart from $F$ - label, the another label and simultaneously visual code of cobweb graded poset is  its "`La scala"'  descending down there to infinity with  picture which looks like that below.

\vspace{0.2cm}

\begin{flushright}

1 - - - - - - - - - - - - - - - - - - - - - - - - - - - - - - - - - - - - - - - - - - - - - - - - - - -\\
1 - - - - - - - - - - - - - - - - - - - - - - - - - - - - - - - - - - - - - - - - - - - - - - - - - -\\
1\textbf{ \textcolor{red}{0}} - - - - - - - - - - - - - - - - - - - - - - - - - - - - - - - - - - - - - - - - - - - - - - -\\
1 - - - - - - - - - - - - - - - - - - - - - - - - - - - - - - - - - - - - - - - - - - - - - - -\\
1 \textbf{\textcolor{red}{0 0}} - - - - - - - - - - - - - - - - - - - - - - - - - - - - - - - - - - - - - - - - - - -\\
1 \textbf{\textcolor{red}{0}} - - - - - - - - - - - - - - - - - - - - - - - - - - - - - - - - - - - - - - - - - - -\\
1 - - - - - - - - - - - - - - - - - - - - - - - - - - - - - - - - - - - - - - - - - - -\\
$F_{5}-1\; \textbf{\textcolor{red}{0}}'s \;\; $1 \textbf{\textcolor{red}{0 0 0 0}} - - - - - - - - - - - - - - - - - - - - -  - - - - - - - - - - - - - - -\\
1 \textbf{\textcolor{red}{0 0 0}} - - - - - - - - - - - - - - - - - - - - - - - - - - - - - - - - - - - -\\
1 \textbf{\textcolor{red}{0 0}} - - - - - - - - - - - - - - - - - - - - - - - - - - - - - - - - - - - -\\
1 \textbf{\textcolor{red}{0}} - - - - - - - - - - - - - - - - - - - - - - - - - - - - - - - - - - - -\\
1 - - - - - - - - - - - - - - - - - - - - - - - - - - - - - - - - - - - -\\
$F_{6}-1\;zeros \quad \quad \quad \quad $1 \textbf{\textcolor{red}{0 0 0 0 0 0 0}} - - - - - - - - - - - - - - - - - - - - - - - - - -\\
1 \textbf{\textcolor{red}{0 0 0 0 0 0}} - - - - - - - - - - - - - - - - - - - - - - - - -\\
1 \textbf{\textcolor{red}{0 0 0 0 0}} - - - - - - - - - - - - - - - - - - - - - - - - -\\
1 \textbf{\textcolor{red}{0 0 0 0}} - - - - - - - - - - - - - - - - - - - - - - - - -\\
1 \textbf{\textcolor{red}{0 0 0}} - - - - - - - - - - - - - - - - - - - - - - - - -\\
1 \textbf{\textcolor{red}{0 0}} - - - - - - - - - - - - - - - - - - - - - - - - -\\
1 \textbf{\textcolor{red}{0}} - - - - - - - - - - - - - - - - - - - - - - - - -\\
1 - - - - - - - - - - - - - - - - - - - - - - - - -\\
$F_{7}-1\;zeros \quad \quad \quad \quad \quad \quad \quad \quad
\quad \quad   \;$1 \textbf{\textcolor{red}{\textcolor{red}{0 0 0 0 0 0 0 0 0 0 0 0}}} - - - - - - - \\
1\textbf{ \textcolor{red}{0 0 0 0 0 0 0 0 0 0 0}} - - - - - - -\\
1\textbf{ \textcolor{red}{0 0 0 0 0 0 0 0 0 0}} - - - - - - -\\
1 \textbf{\textcolor{red}{0 0 0 0 0 0 0 0 0}} - - - - - - -\\
1\textbf{ \textcolor{red}{0 0 0 0 0 0 0 0}} - - - - - - - \\
1\textbf{ \textcolor{red}{0 0 0 0 0 0 0}} - - - - - - - \\
1 \textbf{\textcolor{red}{0 0 0 0 0 0}} - - - - - - - \\
1 \textbf{\textcolor{red}{0 0 0 0 0}} - - - - - - - \\
1 \textbf{\textcolor{red}{0 0 0 0}} - - - - - - - \\
1 \textbf{\textcolor{red}{0 0 0}} - - - - - - - \\
1 \textbf{\textcolor{red}{0 0}} - - - - - - -  \\
1
\textbf{\textcolor{red}{0}} - - - - - - - \\
1 - - - - - - - \\
$F_{8}-1\;zeros \quad \quad \quad \quad \quad \quad \quad \quad
\quad \quad \quad \quad \quad \quad \quad \quad \quad \quad \quad \quad \;\,$.........................\\
{\em and so on}
\end{flushright}

\vspace{0.1cm}
\begin{center}
\noindent \textbf{ Example.3  La Scala di Fibonacci . The  staircase structure  of  incidence
matrix $\zeta_F$}  for \textbf{$F$=Fibonacci sequence} \end{center}

\vspace{0.2cm}

\noindent \textit{Note.} The picture above is drawn for the sequence $F= \left\langle F_\textbf{ \textcolor{red}{1}}, F_2,F_3,...,F_n,... \right\rangle$, where 
$F_k$ are Fibonacci numbers.

\vspace{0.1cm}

\noindent \textbf{Description }of the Figure "`La Scala di Fibonacci"' following [17-22,6].
\noindent If one defines (see:  [17-22] and for earlier references therein as well as in all [3-10]) the  Fibonacci poset  $\Pi=\left\langle P,\leq \right\rangle$ 
with help of  its incidence matrix $\zeta$ representing $P$ uniquely   then one arrives at $\zeta$ with easily recognizable staircase-like  structure - of
zeros in the upper part of this upper triangle matrix   $\zeta$.  This structure is depicted by the Figure "La Scala di Fibonacci"'  where:  empty places
mean zero values (under diagonal)  and    filled with  --   places mean values one (above the diagonal).

\vspace{0.2cm}

\noindent \textbf{Advice.} \textbf{Simultaneous} perpetual \textbf{Exercises}. How the all above and coming figures , formulas and expressions change (simplify) 
in the case of $2^{\left\{1 \right\}}$ replacing the ring  $Z$ of integers in $I(\Pi,Z)$.

\vspace{0.1cm}

\noindent \textbf{Comment 4.} 
The given  $F$-denominated staircase zeros structure above the diagonal of zeta matrix $zeta$ is the \textbf{unique characteristics} of  its corresponding  \textbf{$F$-KoDAG} Hasse digraphs,
where  $F$ denotes \textbf{any} natural numbers valued sequence as shown below.  

\vspace{0.1cm}

\noindent  For that to deliver we use the Gaussian coefficients inherited  \textbf{upside down notation} i.e.  $F_n \equiv  n_F$ (see [2-18], [29-32],and  the Appendix in [9] extracted from [34]) and  recall the Upside Down Notation Principle.\\
Let us also easier the portraying task putting  \textcolor{red} {$n_F = 1 $}. Then - apart from $F$ - label, the another label and simultaneously visual code of cobweb graded poset is  its "`La scala"'  descending down there to infinity with  picture which looks like that below , where

\vspace{0.1cm}
\noindent \textbf{recall}  the $F = \left\langle k_F \right\rangle_{k=0}^n$ is an arbitrary natural numbers valued sequence finite or infinite as $n\in N \cup \left\{0\right\}\cup \left\{\infty\right\}$.

\vspace{0.1cm}

\begin{flushright}
1  ($1_{F}-1)\;zeros$ - - - - - - - - - - - - - - - - - - - - - - - - - -
- - - - - - - - - - - - - - - -\\
............................................................................................................................
0 ... 0 1 0 - - - - - - - - - - - - - - - - - - - - - - - - - - - - - - - - - - - - - - - - - - - - - - 
0 ... 0 0 1 - - - - - - - - - - - - - - - - - - - - - - - - - - - - - - -
- - - - - - - - - - - - - - -\\
0 ... 0 0 0 1 ($2_{F}-1)\;zeros$ - - - - - - - - - - - - - - - - - - - - - -
- - - - - - - - - - - - - \\
..............................................................................................................................                             

0 ... 0 0 0 1 0 - - - - - - - - - - - - - - - - - - - - - - - - - - - - - -  - - - - - - - - - - - - - -\\  
0 ... 0 0 0 0 1 - - - - - - - - - - - - - - - - - - - - - - - - - - - -  - - - - - - - - - - - - - - - -\\
0 ... 0 0 0 0 0 1 ($3_{F}-1) \;zeros$ - - - - - - - - - - - - - - - - - - - - - - - - - - - - - - - -\\
.............................................................................................................................                              
0 ...0 0 0 ...0 1 0 0 - - - - - - - - - - - - - - - - - - - - - - - - - - - - - - - - - - - - - - - - -  
0 ...0 0 0 ...0 0 1 0 - - - - - - - - - - - - - - - - - - - - - - - - - - - - - - - - - - - - - - - - -\\
0 ...0 0 0 ...0 0 0 1 - - - - - - - - - - - - - - - - - - - - - - - - - - - - - - - - - - - - - - -  - - -\\ 
0 ...0 0 0 ...0 0 0 0 1 ($4_{F}-1)\;zeros \;$  - - - - - - - - - - - - - - - - - - - - - - - - - - - - - -\\

{\em and so on}
\end{flushright}
\begin{center}
\textbf{Example.4  La scala  $\textbf{F}$-Generale. The  assumptive, perspicacious  staircase structure  of the incidence
matrix $\zeta_F$ } for \textbf{\textcolor{red}{any}}  $\textbf{F}$ \textbf{natural numbers valued sequence }\end{center}

\vspace{0.2cm} 
 
\noindent Another special case \textbf{Example} is delivered below.

\vspace{0.1cm}

\begin{flushright}

1 - - - - - - - - - - - - - - - - - - - - - - - - - - - - - - - - -
- - - - - - - - - - - - - - - -\\
  1 - - - - - - - - - - - - - - - - - - - - - - - - - - - - - - - - -
- - - - - - - - - - - - - - -\\
    1 \textbf{\textcolor{green}{0 0}} - - - - - - - - - - - - - - - - - - - - - - - - - - - - - - -
- - - - - - - - - - - - -\\
     1 \textbf{\textcolor{green}{0}} - - - - - - - - - - - - - - - - - - - - - - - - - - - - -
- - - - - - - - - - - - - - -\\
       1 - - - - - - - - - - - - - - - - - - - - - - - - - - - -
- - - - - - - - - - - - - - - -\\
         1 \textbf{\textcolor{green}{0 0}} - - - - - - - - - - - - - - - - - - - - - - - - - - - - - -
- - - - -  - - - - -\\
           1 \textbf{\textcolor{green}{0}} - - - - - - - - - - - - - - - - - - - - - - - - - - - - - -
- - - - - -  - - - -\\
                      1 - - - - - - - - - - - - - - - - - - - - - - - - - - - - - -
- - - - - -  - - - -\\    
                                     1 \textbf{\textcolor{green}{0 0}} - - - - - - - - - - - - - - - - - - - - -
- - - - - - - - - - - - - - -\\
                                       1 \textbf{\textcolor{green}{0}} - - - - - - - - - - - - - - - - - - - - -
- - - - - - - - - - - - - - -\\
                                                  1 - - - - - - - - - - - - - - - - - - - - - - - - - - - - - - - - - - - -\\
..........................................\\

{\em and so on}\\
..........................................\\

                        1 \textbf{\textcolor{green}{0 0}} - - - - - - - - - - - - - - - - - - - - - - - - - - - - \\
                           1 \textbf{\textcolor{green}{0}} -  - - - - - - - - - - - - - - - - - - - - - -  - - - - -\\
                              1 - - - -  - - - - - - - - - - - - - - - - - - - - - - - -\\

{\em and so on}
\end{flushright}
\begin{center}
\textbf{Example.5  $\zeta_F$.  The matrix $\zeta$ for  ($0_F=1_F=1$ and $n_F=3$ for $n \geq 2$) the special sequence \textcolor{green}{F} constituting the label sequence denominating  cobweb poset associated to $F$-KoDAG Hasse digraph. }
\end{center}

\vspace{0.2cm}
\noindent \textbf{Advice.} \textbf{Simultaneous} perpetual \textbf{Exercises}. How the all above and coming picture Examples, figures, formulas and expressions change (simplify) 
in the case of $2^{\left\{1 \right\}}$ replacing the ring  $Z$ of integers in $I(\Pi,Z)$.

\vspace{0.3cm}



\noindent \textbf{Graded Posets' $\zeta$ matrix formula.}\\ 
Recall now following [3] that any graded poset  with the finite set of minimal elements is an $F$- sequence denominated sub-poset of  its corresponding cobweb poset.
\noindent The Observation 2 in SNACK supplies the simple recipe for the biadjacency (reduced adjacency) matrix of  Hasse digraph coding  any given  graded poset  with the finite set of minimal elements. The recipe for zeta matrix is then standard. We illustrate this by the [3] source example; the source example as the adjacency  matrices  i.e  zeta matrices of any given  graded poset  with the finite set of minimal elements are sub-matrices of their corresponding cobweb posets and as such have the same block matrix structure and differ "`only"' by eventual additional zeros in upper triangle matrix part while staying to be of the same  cobweb poset block type. 
 
\vspace{0.2cm}

\noindent The explicit  expression for zeta matrix $\zeta_F$ of cobweb posets  via known blocks of zeros and ones for arbitrary natural numbers valued $F$- sequence  was given in [3]  due to more than  mnemonic  efficiency  of the up-side-down notation being applied (see [1] and references therein). With this notation inspired by Gauss  and replacing  $k$ - natural numbers with   "$k_F$"  numbers (Note. The Upside Down Notation Principle has been used in [1]) one gets  :
$$
	\mathbf{A}_F = \left[\begin{array}{llllll}
	0_{1_F\times 1_F} & I(1_F \times 2_F) & 0_{1_F \times \infty} \\
	0_{2_F\times 1_F} & 0_{2_F\times 2_F} & I(2_F \times 3_F) & 0_{2_F \times \infty} \\
	0_{3_F\times 1_F} & 0_{3_F\times 2_F} & 0_{3_F\times 3_F} & I(3_F \times 4_F) & 0_{3_F \times \infty} \\
	0_{4_F\times 1_F} & 0_{4_F\times 2_F} & 0_{4_F\times 3_F} & 0_{4_F\times 4_F} & I(4_F \times 5_F) & 0_{4_F \times \infty} \\
	... & etc & ... & and\ so\ on & ...
	\end{array}\right]
$$

\noindent and

$$
	\zeta_F = exp_\copyright[\mathbf{A}_F] \equiv (1 - \mathbf{A}_F)^{-1\copyright} \equiv I_{\infty\times\infty} + \mathbf{A}_F + \mathbf{A}_F^{\copyright 2} + ... =
$$
$$
	= \left[\begin{array}{lllll}
	I_{1_F\times 1_F} & I(1_F\times\infty) \\
	O_{2_F\times 1_F} & I_{2_F\times 2_F} & I(2_F\times\infty) \\
	O_{3_F\times 1_F} & O_{3_F\times 2_F} & I_{3_F\times 3_F} & I(3_F\times\infty) \\
	O_{4_F\times 1_F} & O_{4_F\times 2_F} & O_{4_F\times 3_F} & I_{4_F\times 4_F} & I(4_F\times\infty) \\
	... & etc & ... & and\ so\ on & ...
	\end{array}\right]
$$
\noindent where  $I (s\times k)$  stays for $(s\times k)$  matrix  of  ones  i.e.  $[ I (s\times k) ]_{ij} = 1$;  $1 \leq i \leq  s,  1\leq j  \leq k.$  and  $n \in N \cup \{\infty\}$

\vspace{0.3cm}

\noindent Particular examples of the above block structure of $\zeta$ matrix (resulting from $\zeta$ being  a result of natural join  operations on the way)  are supplied by Examples 1 ,2,3,4,5 above and Examples 6, 7, 8 represented by  Fig.4, Fig.5, Fig.6 below.  As a matter of fact - \textbf{ all} elements $\sigma$ of the incidence algebra $I(P,R)$ including $\zeta$ i.e. characteristic function of the partial order  $\leq$  or M{\"{o}}bius function $\mu = \zeta^{-1}$ (as exemplified with Examples 9,10,11,12 below) \textbf{have the same block structure} encoded by $F$ sequence chosen. Recall that $R$ from $I(P,R)$ denotes commutative ring and  for example $R$ might be taken to be Boolean algebra $2^{\left\{1 \right\}}$ , the field $Z_2=\left\{0,1\right\}$ the ring $Z$ of integers or real or complex or $p$-adic numbers.

\vspace{0.1cm}

\noindent Namely, arbitrary $\sigma \in I(P,R)$ is of the form

$$
	\sigma = \left[\begin{array}{lllll}
	D_{1_F\times 1_F} & M(1_F\times\infty) \\
	O_{2_F\times 1_F} & D_{2_F\times 2_F} & M(2_F\times\infty) \\
	O_{3_F\times 1_F} & O_{3_F\times 2_F} & D_{3_F\times 3_F} & M(3_F\times\infty) \\
	O_{4_F\times 1_F} & O_{4_F\times 2_F} & O_{4_F\times 3_F} & D_{4_F\times 4_F} & M(4_F\times\infty) \\
	... & etc & ... & and\ so\ on & ...
	\end{array}\right]
$$
\noindent where $D_{k_F\times k_F}$ denotes diagonal  $k_F\times k_F$ matrix  while $M(n_F\times\infty)$  stays for arbitrary  $n_F\times\infty$ matrix and both with matrix elements from the ring  $R$= $2^{\left\{1 \right\}}$ , $Z_2=\left\{0,1\right\}$, $Z$ etc.

\vspace{0.2cm}

\noindent In more detail: it is trivial to note that  all elements  $\sigma \in I(P,R)$  - including  $\zeta^{-1}$  for which   $D_{k_F\times k_F} = I_{k_F\times k_F}$ - are  of  matrix block form resulting from $\os$ of the subsequent bipartite digraphs  $\left\langle \Phi_k \cup \Phi_{k+1},R \right\rangle,\;\;R \subseteq \Phi_k \times \Phi_{k+1},\ \ \left|\Phi_k \right|=k_F $ i.e.

$$
	\sigma  = \left[\begin{array}{llllll}
	D_{1_F\times 1_F} & M(1_F \times 2_F) & M(1_F \times 3_F) & M(1_F \times 4_F) & M(1_F \times 5_F) & M(1_F \times 6_F)\\
	0_{2_F\times 1_F} & D_{2_F\times 2_F} & M(2_F \times 3_F) & M(2_F \times 4_F) & M(2_F \times 5_F) & M(2_F \times 6_F)\\
	0_{3_F\times 1_F} & 0_{3_F\times 2_F} & D_{3_F\times 3_F} & M(3_F \times 4_F) & M(3_F \times 5_F) & M(3_F \times 6_F)\\
	0_{4_F\times 1_F} & 0_{4_F\times 2_F} & 0_{4_F\times 3_F} & D_{4_F\times 4_F} & M(4_F \times 5_F) & M(4_F \times 6_F) \\
	... & etc & ... & and\ so\ on & ...
	\end{array}\right]
$$
where  $M(k_F \times (k+1)_F)$  denote corresponding  $k_F \times (k+1)_F$ matrices with matrix elements from the ring  $R$= $2^{\left\{1 \right\}}$ , $Z_2=\left\{0,1\right\}$, $Z$ etc. \textcolor{red}{\textbf{However...} } for some seemingly most useful of them ...

\vspace{0.2cm}
\noindent  \textcolor{blue}{\textbf{The New Name: $\os$-natural}} $\Leftrightarrow$   $M(k_F \times (k+1)_F)_{r,s} =  c_{r,s} B(k_F \times (k+1)_F)_{r,s}$. \\ 
In the case of $\zeta$  or  August Ferdinand M{\"{o}}bius matrices  motivating examples of specifically natural elements $\sigma \in I(P,R)$ (i.e.  \textcolor{blue}{\textbf{$\os$-natural}}  including those obtained via \textcolor{red}{\textbf{the ruling formula}}) - so in the case of such type elements $\sigma \in I(P,R)$ we ascertain - and may prove via just see it - that
 
$$M(k_F \times (k+1)_F)_{r,s} =  c_{r,s}B(k_F \times (k+1)_F)_{r,s},$$ 
where the rectangular "`zero-one"' $B(k_F \times (k+1)_F)$ matrices from Observation 2. are obtained from the $F$-cobweb poset matrices $I(k_F \times (k+1)_F)$ by replacing some ones by zeros.\\
Moreover (see Observation 3)  - in the case of M{\"{o}}bius  $\mu = \zeta^{-1}$ matrix as it is obligatory  \textcolor{red}{$\textbf{c}_{r,r+1}= \textbf{-1}$}.

\vspace{0.2cm}
\noindent The motivating example of  \textcolor{blue}{\textbf{$\os$-natural}} element of the incidence algebra is $\zeta_F$ due to  the algorithm of \textcolor{red}{\textbf{the ruling formula}} considered over the $R= 2^{\left\{1\right\}}$ ring in this particular case element:

$$
	\zeta_F = I_{\infty\times\infty} + \mathbf{A}_F + \mathbf{A}_F^{\copyright 2} + ...= (1 - \mathbf{A}_F)^{-1\copyright} 
$$
where 

$$
	\mathbf{A}_F = \left[\begin{array}{llllll}
	0_{1_F\times 1_F} & B(1_F \times 2_F) & 0_{1_F \times \infty} \\
	0_{2_F\times 1_F} & 0_{2_F\times 2_F} & B(2_F \times 3_F) & 0_{2_F \times \infty} \\
	0_{3_F\times 1_F} & 0_{3_F\times 2_F} & 0_{3_F\times 3_F} & B(3_F \times 4_F) & 0_{3_F \times \infty} \\
	0_{4_F\times 1_F} & 0_{4_F\times 2_F} & 0_{4_F\times 3_F} & 0_{4_F\times 4_F} & B(4_F \times 5_F) & 0_{4_F \times \infty} \\
	... & etc & ... & and\ so\ on & ...
	\end{array}\right]
$$
and where  $B(k_F \times (k+1)_F)$  are introduced by the Observation 2.

\vspace{0.2cm}
\noindent For the other example of \textcolor{blue}{\textbf{$\os$-natural}} element is $[Max]$  given by the algorithm  of \textcolor{red}{\textbf{the ruling formula}} over the 
$R= Z$ ring see further on in below.

\vspace{0.2cm}
\noindent  For the sake of the forthcoming Observation 1 we introduce  the  set of corresponding Hasse diagram maximal chains called the layer of the graded DAG called KoDAG to be just this [8,6,5,4]:
$$\langle\Phi_k \to \Phi_n \rangle = \left\{c=<x_k,x_{k+1},...,x_n>, \: x_s \in \Phi_s, \:s=k,...,n \right\}.$$

\vspace{0.2cm}

\begin{observen}
(see [3] - and consult the Remark 1). Let us denote by  $\langle\Phi_k\to\Phi_{k+1}\rangle$ the di-bicliques  denominated by subsequent levels $\Phi_k, \Phi_{k+1}$ of the graded  $F$-poset $P(D) = (\Phi, \leq)$  i.e. levels $\Phi_k , \Phi_{k+1}$ of  its cover relation graded digraph  $D = (\Phi,\prec\!\!\cdot$)  [Hasse diagram].   Then

$$
	B\left(\os_{k=1}^n \langle\Phi_k\to\Phi_{k+1}\rangle \right) = \mathrm{diag}(I_1,I_2,...,I_n) = 
$$
$$
	= \left[ \begin{array}{lllll}
	I(1_F\times 2_F) \\
	& I(2_F\times 3_F) \\
	& & I(3_F\times 4_F) \\
	& & ... \\
	& & & & I(n_F\times (n+1)_F)
	\end{array} \right]
$$
\vspace{0.1cm}
\noindent where $I_k \equiv I(k_F \times (k+1)_F)$, $k = 1,...,n$ and where  - recall - $I (s\times k)$  stays for $(s\times k)$  matrix  of  ones  i.e.  $[ I (s\times k) ]_{ij} = 1$;  $1 \leq i \leq  s,  1\leq j  \leq k.$  and  $n \in N \cup \{\infty\}$.  
\end{observen}

\vspace{0.2cm}
\noindent  The binary  natural join operation  $\os$ being defined for such pairs of arguments (matrices, digraphs, graphs,  relations of varying arity,..) which do satisfy the natural join condition 
(see [3] and [5,4,2]) is associative of course iff performable and obviously $\os$ is noncommutative.

\vspace{0.2cm}

\noindent The recipe for any connected - \textbf{hence} $F$-\textbf{denominated} - the recipe for \textbf{any } given graded poset with a finite minimal elements set is supplied via the following observation.

\vspace{0.2cm}

\begin{observen}
(see [3]- and consult the Remark 1). Consider bigraphs' chain obtained from the above di-bicliques' chain via  deleting or no  arcs making thus [if deleting arcs] some or all of the di-bicliques $ \langle\Phi_k\to\Phi_{k+1}\rangle$  not di-bicliques; denote  them as  $G_k$. Let $B_k = B(G_k)$ denotes their biadjacency matrices correspondingly.  Then for any such $F$-denominated chain [hence any chain ] of bipartite digraphs  $G_k$  the general formula is:

$$
 B\left( \os_{i=1}^n G_i \right) \equiv  B [\os_{i=1}^n  A(G_i)] =  \oplus_{i=1}^n  B[A(G_i) ]  \equiv  \mathrm{diag} (B_1 , B_2 , ..., B_n) =
$$
$$
	= \left[ \begin{array}{lllll}
	B_1 \\
	& B_2 \\
	& & B_3 \\
	& & ... \\
	& & & & B_n
	\end{array} \right]
$$

\noindent $n \in N \cup \{\infty\}$.
\end{observen}

\vspace{2mm}

\noindent \textbf{!} Let us recall that $\zeta$ is defined for any poset as follows ($p,q \in P$):
$$
\zeta (p,q)=\left\{ \begin{array}{cl} 1&for\ p \leq q,\\0&otherwise.
\end{array} \right. 
$$
This is the reason why in the above \textcolor{red}{\textbf{ruling formula}}:

$$
	\zeta_F = I_{\infty\times\infty} + \mathbf{A}_F + \mathbf{A}_F^{\copyright 2} + ...= (1 - \mathbf{A}_F)^{-1\copyright} 
$$
the Boolean powers are used. If this rule is applied with $Z$-ring or other ring $R, Z\subseteq R$ powers then we get

$$
	[Max]_F = \mathbf{A}_F^0 + \mathbf{A}_F^1 + \mathbf{A}_F^2 + ...= (1 - \mathbf{A}_F)^{-1}=
$$
$$
= \left[\begin{array}{llllll}
	I_{1_F\times 1_F} & B(1_F \times 2_F) & B(1_F \times 3_F) & B(1_F \times 4_F) & B(1_F \times 5_F) & ... \\
	0_{2_F\times 1_F} & I_{2_F\times 2_F} & B(2_F \times 3_F) & B(2_F \times 4_F) & B(2_F \times 5_F) & ... \\
	0_{3_F\times 1_F} & 0_{3_F\times 2_F} & I_{3_F\times 3_F} & B(3_F \times 4_F) & B(3_F \times 5_F)  & ...\\
	0_{4_F\times 1_F} & 0_{4_F\times 2_F} & 0_{4_F\times 3_F} & I_{4_F\times 4_F} & B(4_F \times 5_F) & ... \\
	... & etc & ... & and\ so\ on & ...
	\end{array}\right]
$$
where  $B(k_F \times (k+1)_F)$  are introduced by the Observation 2.\\ 
It is a matter of simple observation and induction to see that 

$$
B(r_F \times (r+2)_F)= B(r_F \times (r+1)_F)B((r+1)_F \times (r+2)_F)
$$
and consequently for $s>r+2$
$$
B(r_F \times s_F)= B(r_F \times (r+1)_F)B((r+1)_F \times (r+2)_F)...B((s-2)_F \times (s-1)_F)B((s-1)_F \times s_F).
$$
In the case of $F$-cobweb posets - replace  $B(r_F \times s_F)$  by   $I(r_F \times s_F)$ and then one may use the "`Markov"' property.\\
What about then just $F$-graded posets case ? - See  Comment 3 and its Warning.

\vspace{3mm}

\noindent \textbf{Remark \textbf{2}. $F$-graded poset construction - summary.} The knowledge of $\zeta$  matrix explicit form enables one  to construct (calculate) via standard algorithms the M{\"{o}}bius matrix $\mu =\zeta^{-1}$ and other typical elements of incidence algebra perfectly suitable for calculating number of chains, of maximal chains etc. in finite sub-posets of $P$. Right from the definition of $P$ via its Hasse diagram. 
\noindent The way the $\zeta$  is written above underlines the fact that this is the staircase structure encoding formula for \textbf{any} natural numbers valued sequence $F$. Recall: this \textbf{sequence } $\textbf{F}$  \textbf{serves as} \textbf{the label} encoding all resulting digraphs and combinatorial objects.\\
The subsequent di-biclique of bipartite digraph adjoining via natural join $\os$ is in one to one correspondence with  adjoining another subsequent one step down of La Scala. In another words : 
one more step down La Scala - one more di-biclique $\os$-adjoint. 

\vspace{2mm}

\noindent To this end define the $L$-\textbf{L}ogic function as follows:  

$$
L\left([Max]_F\right) = \zeta_F,\;\;\zeta_{r,s}=\Big\{\begin{array}{l}1\;\;\;\;\; [Max]_{r,s}>0 \\0\;\;\;\;\;[Max]_{r,s}=0\end{array}.
$$ 
This completes the natural join $\os$ structural description of $\zeta_F$ matrix construction for any $F$-graded poset and will be of use as a guide while looking 
for the similar form of M{\"{o}}bius matrix $\mu =\zeta^{-1}$ bearing in mind that for $s>r$
$$
B(r_F \times s_F)= \prod_{i=r}^{s-1} B(i_F \times (i+1)_F).
$$

\vspace{4mm}
\noindent \textcolor{blue}{\textbf{Remark 3}. The \textbf{choice of} $F$-poset  $\Pi$  \textbf{labeling} and then \textbf{Knuth notation}.}\\ 
\noindent If one defines \textbf{\textcolor{blue}{any}} graded $F$ poset $P$  with help of  its incidence matrix $\zeta$ representing $P$ uniquely   then \textbf{in case of cobweb posets} one arrives at
$\zeta$ with  \textbf{Type characterization La Scala \textcolor{blue}{code}} of zeros in the upper part of this upper triangle matrix   $\zeta$ due to implicit natural for right-handed oriented choice of nodes labeling. See all figures above. In the case of arbitrary $F$-graded poset $P$ apart from La Scala additional zeros appear. These are the fixed zeros of  $B(i_F \times (i+1)_F)$ yielding all the other zeros from $B(r_F \times s_F)$ in the upper block triangle of $\zeta$ matrix via the product formula above.
\noindent Let us make now this choice of labeling \textbf{ explicit}. For that to do it is enough to focus on any cobweb poset $\Pi$ as a sample case.

\vspace{0.2cm}

\noindent \textcolor{red}{\textbf{Remark 3.1.}}

\vspace{0.2cm}

\noindent A bit of history. The matrix elements of $\zeta(x,y)$ matrix for Fibonacci cobweb poset were given in 2003 ([18,22] Kwa\'sniewski) using   $x,y \in N \cup \left\{\textcolor{blue}{\textbf{0}}\right\}$ labels of  vertices in their "`natural"' order i.e.  applying the natural labeling as in [52] - see Fig.11. Namely:

\begin{figure}[ht]
\begin{center}
	\includegraphics[width=70mm]{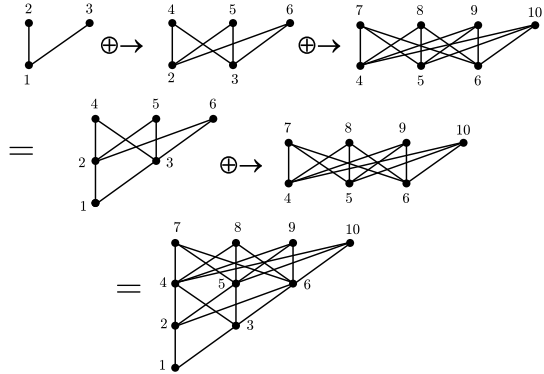}
	\caption{Natural join in natural labeling.} \label{fig:representation}
\end{center}
\end{figure}

\noindent \textbf{1.} set  \textcolor{blue}{\textbf{$k=0$}},\\ 
\textbf{2.} then label subsequent vertices - from the left to the right - along the level  $k$,\\
\textbf{3.} repeat  2.  for $k \rightarrow  k+1$  until  $k=n+1$ ;  $n  \in N\cup\left\{\infty\right\}$ 

\vspace{0.2cm}
\noindent As the result we obtain  the $\zeta$ matrix for Fibonacci sequence as presented by the the Fig. \textit{La Scala di Fibonacci} dating back to 2003 [18,22].

\vspace{0.2cm}

\noindent  \textcolor{blue}{\textbf{The origin - of effectiveness}}. Inspired [29-44] by Gauss $n_q = q^0 + q^1 + ...+ q^{n-1}$ finite geometries numbers and in the spirit of Knuth "`notationlogy"' [23] we shall refer here also to the  upside down notation effectiveness as in [2-6] or earlier in [29-44]. 
\textbf{As for} that upside down attitude  $F_n \equiv n_F$ being much more than "`just a convention"' to be used substantially in what follows as well as  for  the reader's convenience - \textbf{let us} recall it just here quoting it as The Principle according to  Kwa\'sniewski [2] where this rule has been  formulated as an "`\textbf{of course}"'  Principle i.e. simultaneously  trivial and powerful statement.
\vspace{0.1cm}

\begin{quote}
\textbf{The Upside Down Notation Principle.} \\
\textbf{1. Let} the statement $s(F)$ depends only on the fact that $F$ is a natural numbers valued statement.\\
\textbf{2. Then} if one proves that $s(N)\equiv s(\left\langle n\right\rangle_{n \in N})$ is true  - the statement  $s(F)\equiv s(\left\langle n_F\right\rangle_{n \in N})$ is also true. Formally - use equivalence relation classes induced by co-images of $s : \left\{F\right\} \mapsto 2^{\left\{1\right\}}$ and proceed in a standard way.
\end{quote}

\vspace{0.1cm}

\noindent In order to proceed further let us now recall-rewrite  purposely  here Kwa\'sniewski $2003$ - formula for $\zeta$ function of \textbf{arbitrary}  cobweb poset in order  to see that its' algorithm rules automatically  make it valid for all $F$-cobweb posets where $F$ is any natural numbers valued sequence  i.e. with \textcolor{red}{\textbf{$F_0 > 0$}}. $I(\Pi,R)$ stays for the incidence algebra of the poset $\Pi$ over the commutative ring $R$  where  $x,y,k,s \in N \cup \left\{\textcolor{blue}{\textbf{0}}\right\}$.

$$\zeta(x,y) =\zeta_{1}(x,y) -\zeta_{0}(x,y)$$

$$\zeta_{\textcolor{green}{\textbf{1}}}(x,y)=\sum_{k=\textcolor{blue}{\textbf{0}}}^{\infty}\delta(x+k,y)$$
$$\zeta_{0}(x,y)=\sum_{k \geq \textcolor{red}{\textbf{1}}}\sum_{s \geq 0}\delta
(x,F_{s+1}+k)\sum_{r=1}^{F_{s}-k-1}\delta (k+F_{s+1}+r,y)$$ 
and naturally
$$ \delta (x,y)=\Big\{\begin{array}{l}1\;\;\;\;\; x=y\\0\;\;\;\;\;x\neq y\end{array}.$$
\vspace{0.1cm}

\noindent The above formula  for  $\zeta \in I(\Pi,R)$ rewritten in ($F_s \equiv s_F$) upside down notation equivalent form as below
is of course \textbf{\textcolor{red}{valid for all}} cobweb posets  ( $x,y,k,s \in N \cup \left\{\textcolor{blue}{\textbf{0}}\right\}$).

$$\zeta(x,y) =\zeta_{\textcolor{green}{\textbf{\textcolor{green}{\textbf{1}}}}}(x,y) -\zeta_{0}(x,y),$$

$$\zeta_{\textcolor{green}{\textbf{1}}}(x,y)=\sum_{k=\textcolor{blue}{\textbf{0}}}^{\infty}\delta(x+k,y),$$ 
$$\zeta_{0}(x,y)=\sum_{s \geq 1}\sum_{k \geq \textcolor{red}{\textbf{1}}}\delta (x,k+s_F)\sum_{r=1}^{(s-1)_F-k-1}\delta (x+r,y).$$

\noindent \textbf{Note.} $+\zeta_{\textcolor{green}{\textbf{1}}}$ "`produces the Pacific ocean of  \textcolor{green}{1's}"' in the whole upper triangle part of a would be incidence algebra $\sigma \in I(\Pi,R)$ matrix elements with then ($-\zeta_0$)  resulting zeros and ones multiplying arbitrary $\sigma$ choice fixed elements of  $R$],\\

\noindent \textbf{Note.}  $-\zeta_0$ cuts out \textcolor{red}{0's} i.e.  thus producing "`zeros' $F$-La Scala staircase"' in the \textcolor{green}{\textbf{1}'s}  delivered by $+\zeta_{\textcolor{green}{\textbf{1}}}$.\\

\vspace{0.1cm}

\noindent  This results exactly in forming \textcolor{red}{0's} rectangular triangles:  $s_F - 1$ of them at the start of subsequent stair and then  down to one \textcolor{red}{0} till - after  $s_F - 1$ rows passed by one reaches a half-line of $1's$ which is  running to the right- right  to infinity and thus marks the next in order stair of the $F$- La Scala.\\

\vspace{0.1cm}

\noindent The $\zeta$ matrix explicit formula was given for arbitrary graded posets with the finite set of minimal in terms of natural join of bipartite digraphs in [3]. 

\vspace{0.2cm}

\noindent\textbf{Recall 1.} \textcolor{blue}{\textbf{Recapitulation - the  La Scala Mantra.}}\\ 
What was said is equivalent to the fact that the cobweb poset coding La Scala is of the natural join operation origin  thus producing  $\zeta$ matrix [5,4,3] with  one down step of  La Scala being equivalent to $\os$ - adjoining the subsequent bipartite digraph and what results in: (quote from [3], see: Subsection 2.6.)

\vspace{0.1cm}

\noindent The explicit  expression for zeta matrix $\zeta_F$ of cobweb posets  via known blocks of zeros and ones for arbitrary natural numbers valued $F$- sequence  was given in [3]  due to more than  mnemonic  efficiency  of the up-side-down notation being applied (see [6] [v6] Feb 2009 and references therein). With this notation inspired by Gauss  and replacing  $k$ - natural numbers with   
"$k_F$"  numbers - elements of  the  $F $-sequence one gets 
$$
	\mathbf{A}_F = \left[\begin{array}{llllll}
	0_{1_F\times 1_F} & I(1_F \times 2_F) & 0_{1_F \times \infty} \\
	0_{2_F\times 1_F} & 0_{2_F\times 2_F} & I(2_F \times 3_F) & 0_{2_F \times \infty} \\
	0_{3_F\times 1_F} & 0_{3_F\times 2_F} & 0_{3_F\times 3_F} & I(3_F \times 4_F) & 0_{3_F \times \infty} \\
	0_{4_F\times 1_F} & 0_{4_F\times 2_F} & 0_{4_F\times 3_F} & 0_{4_F\times 4_F} & I(4_F \times 5_F) & 0_{4_F \times \infty} \\
	... & etc & ... & and\ so\ on & ...
	\end{array}\right]
$$
and

$$
	\zeta_F = exp_\copyright[\mathbf{A}_F] \equiv (1 - \mathbf{A}_F)^{-1\copyright} \equiv I_{\infty\times\infty} + \mathbf{A}_F + \mathbf{A}_F^{\copyright 2} + ... =
$$
$$
	= \left[\begin{array}{lllll}
	\textcolor{red}{\textbf{I}}_{1_F\times 1_F} & I(1_F\times\infty) \\
	O_{2_F\times 1_F} & \textcolor{red}{\textbf{I}}_{2_F\times 2_F} & I(2_F\times\infty) \\
	O_{3_F\times 1_F} & O_{3_F\times 2_F} & \textcolor{red}{\textbf{I}}_{3_F\times 3_F} & I(3_F\times\infty) \\
	O_{4_F\times 1_F} & O_{4_F\times 2_F} & O_{4_F\times 3_F} & \textcolor{red}{\textbf{I}}_{4_F\times 4_F} & I(4_F\times\infty) \\
	... & etc & ... & and\ so\ on & ...
	\end{array}\right]
$$
where  $I (s\times k)$  stays for $(s\times k)$  matrix  of  ones  i.e.  $[ I (s\times k) ]_{ij} = 1$;  $1 \leq i \leq  s,  1\leq j  \leq k.$  and  $n \in N \cup \{\infty\}$

\vspace{0.2cm}

\noindent In the $\zeta_F $ formula from [5,4,3]  $\copyright$ denotes the Boolean  product, hence - exactly this product is meant while  Boolean powers enter formulas. We readily recognize from its block structure that $F$-La Scala is formed by \textcolor{red}{\textbf{upper zeros}} of block-diagonal matrices $\textcolor{red}{\textbf{I}}_{k_F\times k_F}$ which sacrifice   these their  \textcolor{red}{\textbf{zeros}} to constitute the  $k$-th   subsequent stair in the $F$-La Scala descending and descending far away down
to infinity. Thus the cobweb poset coding La Scala is due to the natural join origin of $\zeta$ matrix. In general case  of any $F$-graded poset (with as in Remark 3.1 labeling fixed)  one naturally encounters - apart from  obligatory La Scala - zeros generated via \textbf{\textit{the ruling formula}} from (Remark.1.) those of    $B(A)$ which is biadjacency i.e cover relation $\prec\cdot$  matrix of the adjacency matrix $A$  of the $F$-graded poset.\\
\textcolor{red}{\textbf{Note}}:  biadjacency and  cover relation $\prec\cdot$  matrix for bipartite digraphs coincide. By extension - we call  \textbf{cover relation} $\prec\cdot$ matrix $\kappa$ the biadjacency matrix too  in order to keep reminiscent convocations going on.


\vspace{0.2cm}

\noindent  Note now that because of $\delta$'s  under summations in the former $\zeta$ formula  the following is obvious:  

$$ 1 \leq r = y - x \leq (s-1)_F - k - 1 \equiv  1 \leq r = y - k - s_F \leq s-1)_F - k -1  \equiv $$
$$ \equiv  1 \leq r = y \leq s_F - (s-1)_F  - 1 .  $$

\vspace{0.1cm}

\noindent Because of that the above last expression of the $\zeta$ expressed in terms of  $\delta \in I(\Pi,R)$  may be still simplified [for the sake of verification and portraying via computer simple program implementation]. Namely the following is true:

$$\zeta(x,y) =\zeta_{\textcolor{green}{\textbf{1}}}(x,y) -\zeta_{0}(x,y),$$

\vspace{0.1cm}

\noindent where

$$\zeta_{\textcolor{green}{1}}(x,y)=\sum_{\textcolor{blue}{k=\textbf{0}}}^{\infty}\delta(x+k,y),$$ 

\noindent [- note: $+\zeta_1$ "`produces the Pacific ocean of \textcolor{green}{1's}"' in the whole upper triangle part of a would be incidence algebra $\sigma \in I(\Pi,R)$ matrix elements with 
then ($-\zeta_0$) resulting zeros and ones multiplying arbitrary $\sigma$ choice fixed elements of  $R$],\\

\vspace{0.1cm}

\noindent and where (where  $x,y,k,s \in N \cup \left\{\textcolor{blue}{\textbf{0}}\right\}$)

$$\zeta_{0}(x,y)=\sum_{s \geq 1}\sum_{k \geq \textcolor{red}{\textbf{1}}}\delta (x,k+s_F)\sum_{r\geq 1}^{s_F + (s-1)_F-1} \delta (r,y),$$ 

\noindent [- note then again that $-\zeta_0$ cuts out "one's $F$-La Scala staircase"' in the \textcolor{green}{1's} provided by $+\zeta_{\textcolor{green}{\textbf{1}}}$].\\

\vspace{0.1cm}

\noindent Note, that for \textcolor{red}{\textbf{$F$ = Fibonacci}} this still more simplifies   as then  $$s_F + (s-1)_F-1 = (s+1)_F.$$

\vspace{0.3cm}

\noindent \textbf{Remark 3.2. ad Knuth notation} [23].\\

\vspace{0.1cm} 

\noindent In the  wise "`notationlogy"' Knuth's note [23]   one finds among others the notation just for the purpose here (see [6]  [v6] Fri, 20 Feb 2009)

$$ [ s ] =\left\{ \begin{array}{cl} 1&if \ s \ is \  true ,\\0&otherwise.

\end{array} \right. $$

\vspace{0.1cm} 

\noindent Consequently for any set or class

$$[x=y] \equiv \delta (x,y). $$

\vspace{0.1cm} 

\noindent Consequently for any set with addition (group, free group, semi-group, ring,...):

$$[x<y] \equiv \sum_{k\geq \textcolor{red}{\textbf{1}}} \delta (x+k,y),$$

$$[x\textcolor{red}{\leq} y] \equiv \sum_{k\geq \textcolor{red}{\textbf{0}}} \delta (x+k,y).$$

\noindent Using this makes my last above expression of the $\zeta$ in terms of  $\delta$ still more transparent and handy  if rewritten in Donald Ervin Knuth's  notation [23]. Namely:

$$\zeta(x,y) =\zeta_{\textcolor{green}{\textbf{1}}}(x,y) -\zeta_{0}(x,y)$$

$$\zeta_{\textcolor{green}{1}}(x,y)= [x \leq y]$$ 

$$\zeta_{0}(x,y)=\sum_{s \geq 1}\sum_{k \geq \textcolor{red}{\textbf{1}}}[x=k+s_F][1 \leq y \leq s_F + (s-1)_F\! - 1 ].$$

$$\zeta_{0}(x,y)=\sum_{s \geq 1}[x>s_F][1 \leq y \leq s_F + (s-1)_F\! - 1 ],$$ 
where, let us  recall: $x,y,k,s \in N \cup \left\{\textcolor{blue}{\textbf{0}}\right\}$.
\vspace{0.1cm}

\noindent Note, that for \textcolor{red}{\textbf{$F$ = Fibonacci}} this still more simplifies   as then  $$s_F + (s-1)_F-1 = (s+1)_F.$$

\vspace{0.3cm}



\noindent \textbf{Remark 3.3. Knuth notation [23] - and  Dziemia\'nczuk's ? guess}

\vspace{0.2cm}

\noindent It was remarked by my Gda\'nsk University Student Coworker Maciej Dziemia\'nczuk  - that my  $\zeta \in I(\Pi,R)$  (equivalent) expressions are valid according to him only for $F$ = Fibonacci sequence.  In view of the Upside Down Notation Principle if any of these  is proved valid for any particular natural numbers  valued sequence $F$ using no other particular properties of $F$ then it should be true for all of the kind.

\vspace{0.2cm}

\noindent His this being doubtful - has led him to invention of his own - in the course of our The Internet Gian Carlo Rota Polish Seminar discussions with me (see [6] 20 Feb 2009). 

\vspace{0.1cm}

\noindent  Here comes the  formula  postulated by him in the course the The Internet Seminar e-mail discussions (see then \textit{resulting now} Comment 5 referring to \textcolor{blue}{\textbf{Krot}}).

\vspace{0.1cm}

$$ \zeta(x,y) =  [x\leq y] - [x<y] \sum_{n\geq 0} [(x> S(n)] [y \leq S(n+1)],$$

\noindent where 

$$ S(n) = \sum_{k\geq 1}^n k_F $$

\vspace{0.2cm}

\noindent \textbf{Exercise.} My todays reply to his guess (compare [19] 20 Feb 2009)  is the following  Exercise.

\noindent  Let $x,y\in N\cup\left\{0\right\}$  be the labels of  vertices in their "`natural"' linear order as explained earlier. 

\noindent Prove the true claim:\\
\textit{Dziemia\'nczuk guess is equivalent to Kwa\'sniewski formulas}.\\ 

\noindent - What is for?  My "`for"' is the Socratic Method  question.  Why not use  the arguments in favor of 
$$\zeta_{0}(x,y)=\sum_{s \geq 1}\sum_{k \geq \textcolor{red}{\textbf{1}}}\delta (x,k+s_F)\sum_{r\geq 1}^{s_F + (s-1)_F-1} \delta (r,y),$$ 

\vspace{0.1cm}

\noindent  Hint. Use the same argumentation.  Hint. Then - contact Comment 5.

\vspace{0.3cm}

\noindent \textbf{Remark 4. Ewa Krot Choice.} While the above is established it is a matter of  simple observation by inspection  to find out how does  the the M{\"{o}}bius matrix $\mu =\zeta^{-1} $ looks like . Using in [24,25] this author example and expression for $\zeta$ matrix this  has been accomplished first (see also [27])   for Fibonacci sequence and then the same formula was declared to be valid for $F$ sequences as above  in [26,27]. Namely the author of [28] states that the M{\"{o}}bius function for the Fibonacci sequence designated cobweb poset can be easily extended to the whole family of all cobweb posets with indication to the reference [28] where one neither finds the proof except for declaration that the validity far all cobweb posets is OK. From the todays perspective the present author should say that it is not so automatic  if definition of cobweb poset via ordinal sum is not uncovered as in \textbf{1.2} above i.e.  $\Pi = \oplus_{k\geq 0} \Phi_k$. For that to see follow what follows.\\
By now here is her formula for  the cobweb posets' M{\"{o}}bius function (see: (6) in [24] then it is recommended to consult Comment 5).\\

\vspace{0.2cm}

\noindent Let $x= \left\langle s,t\right\rangle$ and  $y= \left\langle u,v\right\rangle$ where $1 \leq s \leq F_t$, $1 \leq u \leq F_v$ while $t,v \in \textcolor{red}{\textbf{\textit{N}}}$.\\ 

\vspace{0.1cm}

\noindent These are descriptive and extra external with respect to the Krot formula below conditions imposed in order to stay in accordance with the zeros' "`La Scala di Fibonacci"' structure of the present author "`discovered"' in 2003 [18,19]. This Ewa Krot brave independence declaration step formula was since now on presented by the author of [25-28]  in opposition (?) to the Kwa\'sniewski's choice which makes these conditions being automatically inherited from  $\zeta$ matrix formula by the present author (see 2003 [19] and consequently all relevant papers of Kwa\'sniewski later on till today).

\vspace{0.2cm}

\noindent \textbf{If} these external with respect to formula conditions  are assumed \textbf{then} $\mu$  M{\"{o}}bius function
for Fibonacci cobweb \textcolor{blue}{\textbf{Krot formula}}  reads  ([24])

$$\mu(x,y) = \mu (\left\langle s,t\right\rangle,\left\langle u,v\right\rangle) = \delta(s,u)\delta(t,v) - \delta(t+1,v) +  \sum_{k=2}^{\infty}\delta(t+k,v)(-1)^k\prod_{i= t+1}^{v-1}(F_i-1) $$

\vspace{0.2cm}

\noindent In particular for Fibonacci sequence either $F = \left\langle 0,\textcolor{red}{\textbf{1}},1,2,3,5,8,13,21,34,...\right\rangle$ or  $F = \left\langle \textcolor{red}{\textbf{1}},1,2,3,5,8,13,21,34,...\right\rangle$ we get the right number 

\vspace{0.2cm}

$ \mu (\left\langle 1,1\right\rangle,\left\langle 2,1\right\rangle)= -1,$ \ as \ $\prod_{i= 1+1}^{1-1}(F_i-1) = \prod_{i= 0}^{2}(F_i-1) = 0. $
\vspace{0.2cm}

\noindent The same is right for $F=N$. We shall see also by inspection via Examples below that this is a obviously decisive sensitive starting point in applying the recurrent definition of  M{\"{o}}bius function matrix $\mu$  and its descendant - the block structure of  M{\"{o}}bius function coding matrix $C(\mu)$ - with this latter recurrence for $C(\mu)$   allowing  simple solution simultaneously with combinatorial interpretation of  \textcolor{blue}{\textbf{Krot}}on matrix  $K =(K_s(r_F))$, \ where \  $ K_s(r_F) = \left|C(\mu)_{r,s}\right|$. 

\vspace{0.2cm}

\noindent Now bearing in mind the Upside Down Notation Principle let start to  prepare the formula  \textbf{for all} connected graded posets ($F$-cobweb posets included) with   \textcolor{red}{\textbf{$F_0>0$}} (as it should be for natural numbers valued sequences) and of course for  other natural numbers valued sequences $F$.\\

\vspace{0.1cm}
\noindent  Note the  condition resulting from  \textcolor{red}{\textbf{$F_0>0$}} unavoidable convention: Fibonacci means since now on $F = \left\langle \textcolor{red}{\textbf{1}},1,2,3,5,8,13,21,34,...\right\rangle$

\vspace{0.1cm}
\noindent At first the \textbf{First Step}. Let us formulate equivalent versions of the above Krot formula in coordinate grid  $Z \times Z$ adequately to to the task of verifying it
in the case of Fibonacci sequence $F = \left\langle \textcolor{red}{\textbf{1}},1,2,3,5,8,13,21,34,...\right\rangle$.  
This has been done we arrive at what follows.

\vspace{0.2cm}
\noindent \textbf{Let} $x= \left\langle s,t\right\rangle$ and  $y = \left\langle u,v\right\rangle$ where $1 \leq s \leq t_F$, $1 \leq u \leq v_F$ while $t,v \in \textcolor{red}{\textbf{\textit{N}}}$. \textbf{Then} (\textcolor{blue}{ versions equivalent to Krot formula})

$$\mu(x,y) = \mu (\left\langle s,t\right\rangle,\left\langle u,v\right\rangle) = [(s=u] [t=v]  - [t+1=v] +  \sum_{k=2}^{\infty}[t+k=v](-1)^k\prod_{i= t+1}^{v-1}(i_F-1)$$

$$\mu(x,y) = \mu (\left\langle s,t\right\rangle,\left\langle u,v\right\rangle) = [(s=u] [t=v]  - [t+1=v] +  [t+1<v](-1)^{v-t}\prod_{i= t+1}^{v-1}(i_F-1)$$

\noindent or with sine qua non  \textbf{conditions} being \textbf{implemented} in there:

$$\mu(x,y) = \mu (\left\langle s,t\right\rangle,\left\langle u,v\right\rangle) = [(s=u]\: [t=v]  - [t+1=v]+$$  
$$ +  [t+1<v]\:[1 \leq s \leq t_F ] \:[1 \leq u \leq v_F ]
(-1)^k\prod_{i= t+1}^{v-1}(i_F-1).$$

\vspace{0.2cm}

\noindent The above  M{\"{o}}bius function re-formulas \textcolor{red}{\textbf{if proved}} valid for $F=\textcolor{red}{\textbf{N}}$ thanks to no more than the assumption $n_N \in N$ \textbf{then} it should be literally valid for all natural numbers valued sequences $F$. 

\vspace{0.2cm}

\noindent These formulas for M{\"{o}}bius function appear suitable [check] for  Fibonacci sequence $F = \left\langle \textcolor{red}{\textbf{1}},1,2,3,5,8,13,21,34,...\right\rangle$  as well [check] 
as in the case of  $F = N$ (Example 9. ) and  as well as  in the case of Example 12  (see both below).\\
\vspace{0.2cm}

\noindent As a matter of  fact this might be already expected from the following simple check using any of the equivalent formulas:

$$
 \mu (\left\langle 1,1\right\rangle,\left\langle 2,1\right\rangle)= -1
$$
which is the right number for Fibonacci sequence (or see Example 11) as well as  

$$
 \mu (\left\langle 1,1\right\rangle,\left\langle 2,1\right\rangle)= (\left\langle 1,1\right\rangle,\left\langle 2,2\right\rangle)-1
$$
is right the number for $F=N$ natural number sequence (or see Example 12). The reason for that fact is at hand just by inspection
of Hasse digraphs of cobweb posets under consideration. And these checks are crucial at the start in view of recurrent form of  August Ferdinand M{\"{o}}bius matrix $\mu $ formula.    However...

\vspace{0.2cm}

\noindent \textcolor{blue}{\textbf{However we are in need of The Proof!}} of this \textbf{Krot Formula} for M{\"{o}}bius function for any one  - hence \textbf{for all} of the  relevant sequences $F$. Here let us call  back (hail) the mantra: The Upside Down Notation Principle.\\
Also let us recall here that  due to obvious observation of this article that the natural join  $P \os Q$ of graded posets $\left\langle P,\leq_P \right\rangle$ and $\left\langle Q,\leq_Q \right\rangle$ with correspondingly  maximal (in $P$) and minimal (in $Q$) sets  being identical is expressed by ordinal sum  $P\oplus Q$ (see [1]) one arrives in $\textbf{1.4 }$  at a very simple proof of the M{\"{o}}bius function formula for cobweb posets.

\vspace{0.2cm}

\noindent So - as now we see it - one is in need of  the \textbf{Second  Step}.
\noindent In the sequel this is to be done and we shall use formulas for M{\"{o}}bius function with the structure  inferred  from the  fact that incidence algebra $I(\Pi,R$ elements arise in the sequential natural join of di-bicliques or bipartite digraphs in the general case of $F$-graded posets as to be exemplified and derived below. 
Then  implementation  of the recurrent definition of  M{\"{o}}bius function matrix $\mu$  gives birth to daughter descendant of  $\mu$ i.e. the block structure of  M{\"{o}}bius function coding matrix $C(\mu)$ implying for $C(\mu)$  an  recurrence  allowing  simple solution simultaneously with combinatorial interpretation of  \textcolor{blue}{\textbf{Krot}}on matrix \ $K = ( K_s(r_F))$, \ where  $ K_s(r_F) = \left|C(\mu)_{r,s}\right|$.

\noindent And this is to be this Second Step.

\vspace{0.4cm}

\noindent Before doing this in the next section -  to this end - lets for now continue "`the Krot and Krot-Sieniawska contribution subject"'.  The author of  [24] introduces parallely also another form of  $\zeta$ function formula and since now on - except for [24,46] - in  subsequent papers [25,27,28] their author  uses the formula for  $\zeta$  function in this another form. Namely - this  other form formula for $\zeta$  function in the present authors' grid coordinate system description of the cobweb posets was given by  Krot in her  note on M{\"{o}}bius function and M{\"{o}}bius inversion formula for Fibonacci cobweb poset [25]  with $F$ designating the Fibonacci cobweb  posets. In [26] the  formula the \textbf{Krot formula} for the  M{\"{o}}bius function for Fibonacci sequence $F$ was declared as valid for all cobweb posets i.e. for all natural numbers valued sequences $F$ denominated cobweb posets. 
(Consult also so the  recent note "`On Characteristic Polynomials of the Family of Cobweb Posets"' [28] and see also Comment 5.).

\vspace{0.2cm}



\noindent \textbf{Comment 5.} ad $zeta$ and $\mu $.  Back to 03 Feb 2004 preprint [46] for to see the source of the past in the  future which is present in \textbf{1.4} subsection due to the definition of natural join of graded cobweb  posets via  ordinal sum of independent sets of these cobweb posets. 

\vspace{0.1cm}

\noindent  In [46]  the deliberate task was to consider just the case of \textbf{Fibonacci} sequence in order to  to find the inverse matrix $zeta^{-1}$ of  the $zeta$ from [17] (November 2003) using the present author $\zeta$  matrix expression in terms of the infinite Kronecker delta  matrix $\delta$ from [18] (November 2003) and [19] (December 2003). 
Why  \textbf{ Fibonacci}?  See the formula (5) page 9 in [46]. Applying (5) to the conditions on  the top of the page 9 above the relevant formula (5) one arrives at (6) which afterwords - in this author coordinate grid description reads:

$$x= \left\langle s,t\right\rangle,\   y= \left\langle u,v\right\rangle, \  where \ 1 \leq s \leq F_t, \ 1 \leq u \leq F_v,\  while \ t,v \in \textcolor{red}{\textbf{\textsl{N}}}.$$

\vspace{0.1cm}

\noindent Nevertheless already in  Ewa Krot preprint [46] (see the top of the page 9 above the relevant formulas (5) and (6) ...)  already there the \textbf{general case} conditions are stated which in notation of the present author labeling and upside down notation as well as due to Dziemia\'nczuk's observed Knuth notation  now simply read as follows:

$$[(x> S(n)] [y \leq S(n+1)]$$

\noindent where 

$$ S(n) = \sum_{k\geq 1}^n k_F $$

\noindent and accordingly we now infer \textbf{(}$x,y,k,s,n \in N \cup \left\{\textcolor{blue}{\textbf{0}}\right\}$  [Remark 2.1.]\textbf{)}

$$ \zeta(x,y) =  [x\leq y] - [x<y] \sum_{n\geq \textcolor{blue}{\textbf{0}}} [(x> S(n)] [y \leq S(n+1)].$$

\noindent Note. The author of [24-28,46] consequently avoids the upside down notation. However she had used this notation  then in her  Rota and cobweb posets related dissertation that she had defended with distinction on 30 September 2008 [8]. \textbf{The end of comment.}

\vspace{0.2cm}

\noindent \textbf{No doubt }the $\zeta$  function formulas - the former (Kwa\'sniewski)  and the latter (Krot) are valid for all natural numbers valued sequences $F$. \\
(Let us recall here that  due to obvious observation of this article that the natural join  $P \os Q$ of graded posets $\left\langle P,\leq_P \right\rangle$ and $\left\langle Q,\leq_Q \right\rangle$ with correspondingly  maximal (in $P$) and minimal (in $Q$) sets  being identical is expressed by ordinal sum  $P\oplus Q$ (see [1]) one arrives in $\textbf{1.4 }$  at a very simple proof of the M{\"{o}}bius function formula for cobweb posets).\\
Well, here is this other latter form of Krot formula for  $\zeta$  function (see: (7) in [24] or (1) in [26]).\\

\noindent Let $x= \left\langle s,t\right\rangle$ and  $y= \left\langle u,v\right\rangle$ where $1 \leq s \leq F_t$, $1 \leq u \leq F_v$ while $t,v \in \textcolor{red}{\textbf{\textit{N}}}$. Then

$$\zeta (x,y) = \zeta (\left\langle s,t\right\rangle,\left\langle u,v\right\rangle) = \delta(s,u)\delta(t,v) +  \sum_{k=1}^{\infty}\delta(t+k,v) $$

\noindent where here - recall  ($a,b \in Z$):

$$   \delta(a,b)=\left\{ \begin{array}{cl} 1&for\ a = b,\\0&otherwise.

\end{array} \right. $$

\vspace{0.2cm}
\noindent In February 2009 - in the course of The Internet Gian Carlo Rota Polish Seminar e-mail discussions with the present author - still another $\zeta$ - matrix formula was postulated by  Dziemia\'nczuk - in Knuth notation. See - below. We claim : all are - up to the equivalence of description - the same. See then \textbf{Comment 5.}

\vspace{0.2cm}



\noindent \textbf{Remark 4.1. again on $\zeta$ formulas.}

\vspace{0.1cm}

\noindent Let us compare the above  Krot formula for $\zeta$   with those by Kwa\'sniewski equivalent to the one from the Remark 3.1. ($x,y \in N $) i.e. with

$$\zeta(x,y) = [x \leq y] - \sum_{s \geq 1}\sum_{k \geq \textcolor{red}{\textbf{1}}}[x=k+s_F][1 \leq y \leq s_F + (s-1)_F\! - 1 ],$$ 

$$\zeta(x,y) = [x \leq y] - \sum_{s \geq 1}[x>s_F][1 \leq y \leq s_F + (s-1)_F\! - 1 ],$$ 
where, let us  recall: $k,s \in N \cup \left\{0\right\}$.

\vspace{0.2cm}

\noindent Let us rewrite the above  Krot formula in Knuth notation \textbf{keeping in mind the conditions}

$$1 \leq s \leq F_t, \; 1 \leq u \leq F_v , \ t,v \in \textcolor{red}{\textbf{\textit{N}}},$$
which should have been imposed altogether with:

$$ \zeta (<s,t>,<u,v>) = [s=u] [t=v] + [v>t].  $$

\noindent The above formula with sine qua non  conditions being implemented in there  reads:

$$ \zeta (<s,t>,<u,v>) = [s=u] [t=v] + [v>t][1 \leq s \leq t_F ] [1 \leq u \leq v_F ]  $$
and so,  if written with $\delta$'s it contains three subsequent summations as in the Kwa\'sniewski formula from 2003.



\section {The formula  of inverse zeta matrix for graded posets  with the finite set of minimal elements via natural join of matrices and digraphs technique.}

\vspace{0.1cm}

\noindent \textcolor{red}{\textbf{Training in relabeling - \textit{Exercise}.}}
 
\vspace{0.1cm}

\noindent As we were and are to compare formulas from papers using different labeling - write and/or learn to see formulas from the above and below Observations, definitions etc.  as for $x,y,k,s \in N \cup \left\{\textcolor{blue}{\textbf{0}}\right\}$  on one hand and as for $x,y,k,s \in \textcolor{red}{\textbf{\textsl{\textit{N}}}}$ on the other hand. Because of the comparisons reason 
we shall tolerate and use both being indicated explicitly.

\vspace{0.2cm}

\noindent  Let us start with picture Examples 9,10,11 of inverse zeta matrices subsequently corresponding to picture  Examples 1,2,5. For that to do it is enough for now to use
the recurrent definition of the M{\"{o}}bius function 
$$ \mu (x,y)=\Big\{\begin{array}{l}\;\;\;\;1\;\;\;\;\;\;\;\;\;\;\;\;\;\;\;\;\;\;\;\;\;\;\;\;\;x=y\\-\sum_{x\leq z <y} \mu (x,z),\;\; x<y\end{array}.$$

\vspace{0.2cm}

\noindent   Before doing that note that we deal with $F$-\textbf{graded} posets and contact Remark 1 for notation and typical relations relevant below.  

\vspace{0.2cm}

\noindent\textbf{Recall} \textcolor{blue}{\textbf{What form of the August Ferdinand M{\"{o}}bius matrix  we do expect by now.}} 

\vspace{0.1cm}

\noindent Recall:  (see Observation 3)  - in the case of M{\"{o}}bius  $\mu = \zeta^{-1}$ matrix as it is obligatory  \textcolor{red}{$\textbf{c}_{r,r+1}= \textbf{-1}$}.

\vspace{0.1cm}

\noindent Recall (Remark 1) Markov property and observe by inspection that - in the case of M{\"{o}}bius  $\mu = \zeta^{-1}$ matrix \textbf{\textit{for cobweb posets}} it is obligatory to put  
$$
M\left(r_F\times (r+2)_F\right)  = - [ I\left(r_F\times (r+1)_F\right) I\left( (r+1)_F\times (r+2)_F\right) - I\left(r_F\times (r+2)_F\right)] 
$$
i.e.
$$
M\left(r_F\times (r+2)_F\right) = - \textbf{[(r+1)}_F-\textbf{1]}I\left(r_F\times (r+2)_F\right),
$$
thereby :

$$
c_{r,r+2} = - \textbf{[(r+1)}_F-\textbf{1]}c_{r,r+1}, \;\; c_{r,r+1}= -1.
$$
- What about then with arbitrary $F$-graded posets $(P,\leq )$ ?

\vspace{0.2cm}

\noindent In what follows we consider (consult the Remark 1.) motivating examples and then representative Examples 9,10,11,12 of  M{\"{o}}bius matrix. After that the looked for 
\textbf{Theorem 4.} is stated for arbitrary $F$-graded posets $(P,\leq )$.

\vspace{0.2cm}

\noindent  \textbf{Motivating examples}.

\vspace{0.2cm}

\noindent Example 1. Let $i = 1,...,r_F$, $k = 1,...,(r+1)_F$,  $j = 1,...,(r+2)_F$  as now we consider (Remark 1.)  $x_{r,i} \prec\cdot x_{r+1,k}$ where $\left\{ x_{r,i}  \right\}= \Phi_r$ and  $\left\{ x_{r+1,k}  \right\}= \Phi_{r+1}$  are independent sets. Then  \\ 

$$
\mu (x_{r,i},x_{r+2,j}) =  -\sum_{x_{r,i}\leq z <x_{r+2,j}} \mu (x_{r,i},z)= -\left(1+ \sum_{k=1}^{(r+1)_F} \mu (x_{r,i},x_{r+1,k})\right), 
$$
i.e. 
$$
\mu (x_{r,i},x_{r+2,j}) = +[(r+1)_F - 1] = c_{r,r+2} = -[(r+1)_F - 1] c_{r,r+1}.
$$

\vspace{0.2cm}

\noindent Example 2. From Example 1 we infer that as  $\mu (x_{r,i},x_{r+2,j})= \mu (x_r,x_{r+2})$ then it is now enough to consider what follows   ($x_r,x_{r+3}$  any fixed):
$$
\mu(x_r,x_{r+3}) =  -\sum_{x_r\leq z <x_{r+3}} \mu(x_r,z)= -\left(1+ \sum_{x_{r+1}\leq z <x_{r+3}} \mu (x_r,z)\right)= 
$$

$$
= -\left(1+ (r+1)_F \mu(x_r,x_{r+1}) + \sum_{x_{r+2}\leq z <x_{r+3}} \mu (x_r,z)\right)= -\left(1- (r+1)_F  + (r+2)_F \mu (x_r,x_{r+2})\right),
$$
i.e. 
$$
\mu(x_r,x_{r+3}) = -[(r+ 2)_F - 1]c_{r,r+2}= -[(r+ 2)_F - 1][(r+ 1)_F - 1].
$$

\vspace{0.2cm}

\noindent  Via straightforward induction  we conclude that now for arbitrary $r,s \in N \cup \left\{0\right\} $ and for any \textbf{cobweb poset} - in accordance with subsection 1.4. induction - the following is true.

\vspace{0.2cm}

\noindent \textbf{Theorem 3 for cobweb posets.} ($ N \cup \left\{0\right\}$.)

$$ c_{r,s}= [s=r]  - [s=r+1] + [s>r+] (-1)^{s-r}\left((s-r-1)_F - 1)\right)...\left(3_F - 1) \right) \; (+1) =$$
$$ = [s=r] - [s=r+1] + [s>r+1](-1)^{s-r}\; \prod_{i=r+1}^{s-1}(i_F - 1).$$

\vspace{0.2cm}

\noindent Let us see now how it works and how this theorem may be extended to general case of arbitrary $F$-denominated poset. At first the representative Examples 9,10,11,12 of  M{\"{o}}bius matrix follow which might be derived right from the recurrent definition of M{\"{o}}bius  function without even referring to the above theorem .\\

$$ \left[\begin{array}{ccccccccccccccccc}
\textbf{1} & -1 & -1 & +1 & +1 & +1 & -2 & -2 & -2 & -2 & +6 & +6 & +6 & +6 & +6 & -24\cdots\\
0 & \textbf{1} & \textbf{\textcolor{blue}{0}} & -1 & -1 & -1 & +2 & +2 & +2 & +2 & -6 & -6 & -6 & -6 & -6 & +24\cdots\\
0 & \textbf{0} & \textbf{1} & -1 & -1 & -1 & +2 & +2 & +2 & +2 & -6 & -6 & -6 & -6 & -6 & +24\cdots\\
0 & 0 & 0 & \textbf{1} & \textbf{\textcolor{blue}{0}} & \textbf{\textcolor{blue}{0}} & -1 & -1 & -1 & -1 & +3 & +3 & +3 & +3 & +3 & -12\cdots\\
0 & 0 & 0 & \textbf{0} & \textbf{1} & \textbf{\textcolor{blue}{0}} & -1 & -1 & -1 & -1 & +3 & +3 & +3 & +3 & +3 & -12\cdots\\
0 & 0 & 0 & \textbf{0} & \textbf{0} &\textbf{1} & -1 & -1 & -1 & -1 & +3 & +3 & +3 & +3 & +3 & -12\cdots\\
0 & 0 & 0 & 0 & 0 & 0 & \textbf{1} & \textbf{\textcolor{blue}{0}}& \textbf{\textcolor{blue}{0}}&\textbf{\textcolor{blue}{0}} & -1 & -1 & -1 & -1 & -1 & +4\cdots\\
0 & 0 & 0 & 0 & 0 & 0 & \textbf{0} & \textbf{1} & \textbf{\textcolor{blue}{0}} & \textbf{\textcolor{blue}{0}} & -1 & -1 & -1 & -1 & -1 & +4\cdots\\
0 & 0 & 0 & 0 & 0 & 0 & \textbf{0} & \textbf{0} & \textbf{1 }& \textbf{\textcolor{blue}{0}}& -1 & -1 & -1 & -1 & -1 & +4\cdots\\
0 & 0 & 0 & 0 & 0 & 0 & \textbf{0} & \textbf{0} & \textbf{0} & \textbf{1} & -1 & -1 & -1 & -1 & -1 & +4\cdots\\
0 & 0 & 0 & 0 & 0 & 0 & 0 & 0 & 0 & 0 & \textbf{1 }& \textbf{\textcolor{blue}{0}} & \textbf{\textcolor{blue}{0}} & \textbf{\textcolor{blue}{0}}& \textbf{\textcolor{blue}{0}}& -1\cdots\\
0 & 0 & 0 & 0 & 0 & 0 & 0 & 0 & 0 & 0 & \textbf{0} & \textbf{1} & \textbf{\textcolor{blue}{0}}& \textbf{\textcolor{blue}{0}} & \textbf{\textcolor{blue}{0}} & -1\cdots\\
0 & 0 & 0 & 0 & 0 & 0 & 0 & 0 & 0 & 0 & \textbf{0} & \textbf{0} & \textbf{1} & \textbf{\textcolor{blue}{0}} & \textbf{\textcolor{blue}{0}} & -1\cdots\\
0 & 0 & 0 & 0 & 0 & 0 & 0 & 0 & 0 & 0 & \textbf{0} & \textbf{0} & \textbf{0} & \textbf{1} & \textbf{\textcolor{blue}{0}}& -1\cdots\\
0 & 0 & 0 & 0 & 0 & 0 & 0 & 0 & 0 & 0 & \textbf{0} & \textbf{0} & \textbf{0} & \textbf{0} & \textbf{1} & -1\cdots\\
0 & 0 & 0 & 0 & 0 & 0 & 0 & 0 & 0 & 0 & 0 & 0 & 0 & 0 & 0 &  \textbf{1} \cdots\\
. & . & . & . & . & . & . & . & . & . & . & . & . & . & . &  .\cdots\\
 \end{array}\right]$$

\begin{center}
\noindent \textbf{Example.9  $\zeta_N^{-1}$.  The  M{\"{o}}bius function matrix
$\mu = \zeta^{-1}$ for the  natural numbers  i.e. $N$ - cobweb poset.}
\end{center}

\vspace{0.2cm}
$$
	 \mu_N = \left[\begin{array}{lllll}
\textbf{I}_{1\times 1} & \textcolor{red}{-\textbf{I}}(1\times2) & +I(1\times 3) & -2I(1\times4) & +6I(1\times 5) \\
	\textcolor{red}{\textbf{O}}_{2\times 1} & \textbf{I}_{2\times 2} & \textcolor{blue}{-\textbf{I}}(2\times 3) & -2I(2\times 4) & -6I(2\times 5) \\
	O_{3\times 1} & \textcolor{blue}{\textbf{O}}_{3\times 2} & \textbf{I}_{3\times 3} & \textcolor{green}{-\textbf{I}}(3\times 4) & +3I(3\times 5) \\
	O_{4\times 1} & O_{4\times 2} & \textcolor{green}{\textbf{O}}_{4\times 3} & \textbf{I}_{4\times 4} & \textcolor{red}{-\textbf{I}}(4\times 5) \\
	O_{5\times 1} & O_{5\times 2} & O_{5\times 3} & \textcolor{red}{\textbf{O}}_{5\times 4} & \textbf{I}_{5\times 5} \\
	... & etc & ... & and\ so\ on & ...
	\end{array}\right]
$$

\vspace{1mm}
\noindent \textbf{Note.} $\mu$ has \textbf{of course}  natural join inherited structure, of course.
\vspace{1mm}

\begin{center}
\noindent \textbf{Fig.9a   $\mu_N = \zeta_N^{-1}$.  The \textit{block presentation} of the  M{\"{o}}bius function matrix
$\mu = \zeta^{-1}$ for the  natural numbers  i.e. $N$ - cobweb poset.}
\end{center}

\vspace{0.2cm}

\noindent The \textbf{ code} for this KoDAG is given by its KoDAG self-evident code-triangle of the \textcolor{blue}{\textbf{coding matrix}} $C(\mu_F)$( a starting part of it shown below):

\vspace{0.2cm}
$$
	 C(\mu_N) = \left[\begin{array}{llllll}
	+1  & -1 & +1 & -2 & +6  &  -24\\
	-0  & +1  & -1 & +2 & -6  &  +24\\
  +0  & -0  & +1 & -1 & +3  &  -12\\
  -0  & +0  & -0 & +1 & -1  &  +4 \\
  +0  & -0  & +0 & -0 & +1  &  -1\\
	. & . & . & . & .
	\end{array}\right]
$$

\vspace{2mm}

$$ \left[\begin{array}{ccccccccccccccccc}
\textbf{1} & -1 & +0 & -0 & +0 & -0 & +0 & -0 & +0 & -0 & +0 & -0 & +0 & -0 & +0 & -0 & \cdots\\
0 & \textbf{1} & -1 & +0 & -0 & +0 & -0 & +0 & -0 & +0 & -0 & +0 & -0 & +0 & -0 & +0 & \cdots\\
0 & 0 & \textbf{1} & -1 & -1 & +1 & +1 & +1 & -2 & -2 & -2 & -2 & -2 & +8 & +8 & +8 & \cdots\\
0 & 0 & 0 & \textbf{1} & \textbf{\textcolor{red}{0}} & -1 & -1 & -1 & +2 & +2 & +2 & +2 & +2 & -8 & -8 & -8 & \cdots\\
0 & 0 & 0 & \textbf{0} & \textbf{1} & -1 & -1 & -1 & +2 & +2 & +2 & +2 & +2 & -8 & -8 & -8 & \cdots\\
0 & 0 & 0 & 0 & 0 & \textbf{1} & \textbf{\textcolor{red}{0}} & \textbf{\textcolor{red}{0}} & -1 & -1 & -1 & -1 & -1 & +4 & +4 & +4 & \cdots\\
0 & 0 & 0 & 0 & 0 & \textbf{0} & \textbf{1} & \textbf{\textcolor{red}{0}} & -1 & -1 & -1 & -1 & -1 & +4 & +4 & +4 & \cdots\\
0 & 0 & 0 & 0 & 0 & \textbf{0} & \textbf{0} & \textbf{1} & -1 & -1 & -1 & -1 & -1 & +4 & +4 & +4 & \cdots\\
0 & 0 & 0 & 0 & 0 & 0 & 0 & 0 & \textbf{1 }& \textbf{\textcolor{red}{0}} & \textbf{\textcolor{red}{0}} & \textbf{\textcolor{red}{0}} & \textbf{\textcolor{red}{0}} & -1 & -1 & -1 & \cdots\\
0 & 0 & 0 & 0 & 0 & 0 & 0 & 0 & \textbf{0} & \textbf{1} & \textbf{\textcolor{red}{0}}& \textbf{\textcolor{red}{0}} & \textbf{\textcolor{red}{0}} & -1 & -1 & -1 & \cdots\\
0 & 0 & 0 & 0 & 0 & 0 & 0 & 0 & \textbf{0} & \textbf{0} & \textbf{1} & \textbf{\textcolor{red}{0}}& \textbf{\textcolor{red}{0}} & -1 & -1 & -1 & \cdots\\
0 & 0 & 0 & 0 & 0 & 0 & 0 & 0 & \textbf{0} & \textbf{0} & \textbf{0} & \textbf{1} & \textbf{\textcolor{red}{0}}& -1 & -1 & -1 & \cdots\\
0 & 0 & 0 & 0 & 0 & 0 & 0 & 0 & \textbf{0} & \textbf{0} & \textbf{0} & \textbf{0} & \textbf{1} & -1 & -1 & -1 & \cdots\\
0 & 0 & 0 & 0 & 0 & 0 & 0 & 0 & 0 & 0 & 0 & 0 & 0 & \textbf{1} & \textbf{\textcolor{red}{0}} & \textbf{\textcolor{red}{0}} & \cdots\\
0 & 0 & 0 & 0 & 0 & 0 & 0 & 0 & 0 & 0 & 0 & 0 & 0 & \textbf{0}& \textbf{1} & \textbf{\textcolor{red}{0}} & \cdots\\
0 & 0 & 0 & 0 & 0 & 0 & 0 & 0 & 0 & 0 & 0 & 0 & 0 & \textbf{0} & \textbf{0} & \textbf{1} & \cdots\\
. & . & . & . & . & . & . & . & . & . & . & . & . & . & . & . & . \cdots\\
\end{array}\right]$$
\begin{center}
\noindent \textbf{Example.10  $\zeta_F^{-1}$.  The  M{\"{o}}bius function matrix
$\mu = \zeta^{-1}$ for $F$=Fibonacci sequence.} 
\end{center}

\vspace{0.2cm}

$$
	 \mu_F = \left[\begin{array}{lllll}
	I_{1\times 1} & -I(1\times1) & 0I(1\times 1) & 0I(1\times 2) & 0I(1\times 3)\\
	O_{1\times 1} & I_{1\times 1} & -I(1\times 1) & 0I(1\times 2) & 0I(1\times 3)\\
	O_{1\times 1} & O_{1\times 1} & I_{1\times 1} & -I(1\times 2) & +I(1\times 3)\\
	O_{2\times 1} & O_{2\times 1} & O_{2\times 1} & I_{2\times 2} & -I(2\times 3)\\
	O_{3\times 1} & O_{3\times 1} & O_{3\times 1} & 0_{3\times 2} & I_{3\times 3}\\
	... & etc & ... & and\ so\ on & ...
	\end{array}\right]
$$\begin{center}
\noindent \textbf{Example.10a  $\zeta_F^{-1}$.  The \textit{block presentation} of the M{\"{o}}bius function matrix
$\mu = \zeta^{-1}$ for $F$=Fibonacci sequence.} 
\end{center}
Recall then and note here up and below the block structure.

$$
	\sigma  = \left[\begin{array}{llllll}
	I_{1_F\times 1_F} & B(1_F \times 2_F) & B(1_F \times 3_F) & B(1_F \times 4_F) & B(1_F \times 5_F) & B(1_F \times 6_F)\\
	0_{2_F\times 1_F} & I_{2_F\times 2_F} & B(2_F \times 3_F) & B(2_F \times 4_F) & B(2_F \times 5_F) & B(2_F \times 6_F)\\
	0_{3_F\times 1_F} & 0_{3_F\times 2_F} & I_{3_F\times 3_F} & B(3_F \times 4_F) & B(3_F \times 5_F) & B(3_F \times 6_F)\\
	0_{4_F\times 1_F} & 0_{4_F\times 2_F} & 0_{4_F\times 3_F} & I_{4_F\times 4_F} & B(4_F \times 5_F) & B(4_F \times 6_F) \\
	... & etc & ... & and\ so\ on & ...
	\end{array}\right]
$$
where  $B(k_F \times (k+1)_F)$  denote corresponding constant  $k_F \times (k+1)_F$ matrices in the case of $\zeta$ or $\zeta^{-1}$ matrices for example, with matrix elements from the ring  $R$= $2^{\left\{1 \right\}}$ , $Z_2=\left\{0,1\right\}$, $Z$ etc.

\vspace{1mm}
$$ \left[\begin{array}{ccccccccccccccccc}
\textbf{1} & -1 & +0 & -0 & +0 & -0 & +0 & -0 & +0 & -0 & +0 & - 0 & + 0 & - 0 & +0 & -0 & \cdots\\
0 & \textbf{1} & -1 & -1 & -1 & +2 & +2 & +2 & -4 & -4 & -4 & +8 & +8 & +8 & -16 & -16 & \cdots\\
0 & 0 & \textbf{1} &  \textbf{\textcolor{red}{0}} &  \textbf{\textcolor{red}{0}} & -1 & -1 & -1 & +2 & +2 & +2 & -4 & -4 & -4 & +8 & +8 & \cdots\\
0 & 0 & \textbf{0}& \textbf{1} & \textbf{\textcolor{red}{0}} & -1 & -1 & -1 & +2 & +2 & +2 & -4 & -4 & -4 & +8 & +8 & \cdots\\
0 & 0 & \textbf{0} & \textbf{0} & \textbf{1} & -1 & -1 & -1 & +2 & +2 & +2 & -4 & -4 & -4 & +8 & +8 & \cdots\\
0 & 0 & 0 & 0 & 0 & \textbf{1} & \textbf{\textcolor{red}{0}} & \textbf{\textcolor{red}{0}} & -1 & -1 & -1 & +2 & +2 & +2 & -4 & -4 & \cdots\\
0 & 0 & 0 & 0 & 0 & \textbf{0} & \textbf{1} & \textbf{\textcolor{red}{0}} & -1 & -1 & -1 & +2 & +2 & +2 & -4 & -4 & \cdots\\
0 & 0 & 0 & 0 & 0 & \textbf{0} & \textbf{0} & \textbf{1} & -1 & -1 & -1 & +2 & +2 & +2 & -4 & -4 & \cdots\\
0 & 0 & 0 & 0 & 0 & 0 & 0 & 0 & \textbf{1 }& \textbf{\textcolor{red}{0}} & \textbf{\textcolor{red}{0}} & -1 & -1 & -1 & +2 & +2 & \cdots\\
0 & 0 & 0 & 0 & 0 & 0 & 0 & 0 & \textbf{0} & \textbf{1} & \textbf{\textcolor{red}{0}}& -1 & -1 & -1 & +2 & +2 & \cdots\\
0 & 0 & 0 & 0 & 0 & 0 & 0 & 0 & \textbf{0} & \textbf{0} & \textbf{1} & -1 & -1 & -1 & +2 & +2 & \cdots\\
0 & 0 & 0 & 0 & 0 & 0 & 0 & 0 & 0 & 0 & 0 & \textbf{1} & \textbf{\textcolor{red}{0}}& \textbf{\textcolor{red}{0}} & -1 & -1 & \cdots\\
0 & 0 & 0 & 0 & 0 & 0 & 0 & 0 & 0 & 0 & 0 & \textbf{0} & \textbf{1} & \textbf{\textcolor{red}{0}} & -1 & -1 & \cdots\\
0 & 0 & 0 & 0 & 0 & 0 & 0 & 0 & 0 & 0 & 0 & \textbf{0} & \textbf{0} & \textbf{1} & -1 & -1 & \cdots\\
0 & 0 & 0 & 0 & 0 & 0 & 0 & 0 & 0 & 0 & 0 & 0 & 0 & 0 & \textbf{1} & \textbf{\textcolor{red}{0}} & \cdots\\
0 & 0 & 0 & 0 & 0 & 0 & 0 & 0 & 0 & 0 & 0 & 0 & 0 & 0 & 0 & \textbf{1} & \cdots\\
. & . & . & . & . & . & . & . & . & . & . & . & . & . & . & . & . \cdots\\
 \end{array}\right]$$

\begin{center}
\noindent \textbf{\textcolor{red}{Example.11}  $\zeta_F^{-1}$.  The  M{\"{o}}bius function matrix
$\mu = \zeta^{-1}$ for  ($1_F =2_F =1$ and $n_F=3$ for $n \geq 2$) the $F=Fibonacci$ relative special sequence $\textbf{\textcolor{red}{F}}$ constituting\textbf{ the} 
label sequence denominating cobweb poset associated to \textbf{$F$-KoDAG} Hasse digraph }
\end{center}

\vspace{0.3cm}

$$
	 \mu_F = \left[\begin{array}{lllll}
	I_{1\times 1} & -I(1\times1) & +0I(1\times 3) & -0I(1\times 3) & +0I(1\times 3)\\
	O_{1\times 1} & +I_{1\times 1} & -I(1\times 3) & +2I(1\times 3) & -4I(1\times 3)\\
	O_{3\times 1} & -O_{3\times 1} & +I_{3\times 3} & -I(3\times 3) & +2I(3\times 3)\\
	O_{3\times 1} & +O_{3\times 1} & -O_{3\times 3} & +I_{3\times 3} & -I(3\times 3)\\
	O_{3\times 1} & -O_{3\times 1} & +O_{3\times 3} & -0_{3\times 3} & +I_{3\times 3}\\
	... & etc & ... & and\ so\ on & ...
	\end{array}\right]
$$
\begin{center}
\noindent \textbf{\textcolor{red}{Example.11a} $\zeta_F^{-1}$.   The \textit{block presentation} of the  M{\"{o}}bius function matrix
$\mu = \zeta^{-1}$ for  ($1_F=2_F=1$ and $n_F=3$ for $n \geq 2$) the $F=Fibonacci$ relative special sequence $\textbf{\textcolor{red}{F}}$ constituting\textbf{ the} 
label sequence denominating cobweb poset associated to \textbf{$F$-KoDAG} Hasse digraph }
\end{center}

\vspace{0.3cm}

\noindent The  \textcolor{red}{\textbf{code}} for this KoDAG is given by its KoDAG self-evident code-triangle of the \textcolor{red}{\textbf{coding matrix}} $C(\mu_F)$ (a starting part of it shown below):

\vspace{0.2cm}
$$
	 C(\mu_F) = \left[\begin{array}{llllll}
	1 & -1 & +0 & -0 &  +0 & -0\\
	0 &  +1 & -1 & +2 &  -4 & +8\\
  0 &  -0 &  +1& -1 &  +2 & -4\\
	0 &  +0 &  -0 &  +1 & -1 & +2\\
	0 &  -0 &  +0 &  -0 & +1 & -1\\
	. & . & . & . & .
	\end{array}\right]
$$

\vspace{2mm}

$$ \left[\begin{array}{ccccccccccccccccc}
\textbf{1} & -1 & -1 & -1 & +2 & +2 & +2 & -4 & -4 & -4 & +8 & +8 & +8 & -16 & -16 & -16  \cdots\\
0 & +\textbf{1} &  +\textbf{\textcolor{blue}{0}} &  +\textbf{\textcolor{blue}{0}} & -1 & -1 & -1 & +2 & +2 & +2 & -4 & -4 & -4 & +8 & +8 & +8 \cdots\\
0 & -\textbf{0}& +\textbf{1} & +\textbf{\textcolor{blue}{0}} & -1 & -1 & -1 & +2 & +2 & +2 & -4 & -4 & -4 & +8 & +8 & +8\cdots\\
0 & +\textbf{0} & -\textbf{0} & +\textbf{1} & -1 & -1 & -1 & +2 & +2 & +2 & -4 & -4 & -4 & +8 & +8 & +8\cdots\\
0 & -0 & +0 & -0 & +\textbf{1} & +\textbf{\textcolor{blue}{0}} & +\textbf{\textcolor{blue}{0}} & -1 & -1 & -1 & +2 & +2 & +2 & -4 & -4 & -4\cdots\\
0 & +0 & -0 & +0 & -\textbf{0} & +\textbf{1} & +\textbf{\textcolor{blue}{0}} & -1 & -1 & -1 & +2 & +2 & +2 & -4 & -4 & -4\cdots\\
0 & -0 & +0 & -0 & +\textbf{0} & -\textbf{0} & +\textbf{1} & -1 & -1 & -1 & +2 & +2 & +2 & -4 & -4 & -4\cdots\\
0 & +0 & -0 & +0 & -0 & +0 & -0 & +\textbf{1 }& +\textbf{\textcolor{blue}{0}} & +\textbf{\textcolor{blue}{0}} & -1 & -1 & -1 & +2 & +2 & +2\cdots\\
0 & -0 & +0 & -0 & +0 & -0 & +0 & -\textbf{0} & +\textbf{1} & +\textbf{\textcolor{blue}{0}}& -1 & -1 & -1 & +2 & +2 & +2\cdots\\
0 & +0 & -0 & +0 & -0 & +0 & -0 & +\textbf{0} & -\textbf{0} & +\textbf{1} & -1 & -1 & -1 & +2 & +2 & +2\cdots\\
0 & -0 & +0 & -0 & +0 & -0 & +0 & -0 & +0 & -0 & +\textbf{1} & +\textbf{\textcolor{blue}{0}}& +\textbf{\textcolor{blue}{0}} & -1 & -1 & -1\cdots\\
0 & +0 & -0 & +0 & -0 & +0 & -0 & +0 & -0 & +0 & -\textbf{0} & +\textbf{1} & \textbf{\textcolor{blue}{0}} & -1 & -1 & -1\cdots\\
0 & -0 & +0 & -0 & +0 & -0 & +0 & -0 & +0 & -0 & +\textbf{0} & -\textbf{0} & +\textbf{1} & -1 & -1 & -1\cdots\\
0 & +0 & -0 & +0 & -0 & +0 & -0 & +0 & -0 & +0 & -0 & +0 & -0 & +\textbf{1} & +\textbf{\textcolor{blue}{0}} & +\textbf{\textcolor{blue}{0}}\cdots\\
0 & -0 & +0 & -0 & +0 & -0 & +0 & -0 & +0 & -0 & +0 & -0 & +0 & -0 & +\textbf{1} & +\textbf{\textcolor{blue}{0}}\cdots\\
0 & +0 & -0 & +0 & -0 & +0 & -0 & +0 & -0 & +0 & -0 & +0 & -0 & +0 & -0 & +\textbf{1}\cdots\\
. & . & . & . & . & . & . & . & . & . & . & . & . & . & . & . & . \cdots\\
 \end{array}\right]$$

\begin{center}
\noindent \textbf{\textcolor{blue}{Example.12}  $\zeta_F^{-1}$.  The  M{\"{o}}bius function matrix
$\mu = \zeta^{-1}$ for  ($1_F=1$ and $n_F=3$ for $n \geq 2$) the $N$ relative special sequence $\textbf{\textcolor{blue}{F}}$ constituting\textbf{ the} 
label sequence denominating cobweb poset associated to \textbf{$F$-KoDAG} Hasse digraph }
\end{center}

\vspace{0.3cm}

$$
	 \mu_F = \left[\begin{array}{lllll}
	I_{1\times 1} & -I(1\times1) & +2I(1\times 3) & -4I(1\times 3) & +8I(1\times 3)\\
	O_{1\times 1} & +I_{1\times 1} & -I(1\times 3) & +2I(1\times 3) & -4I(1\times 3)\\
	O_{3\times 1} & -O_{3\times 1} & +I_{3\times 3} & -I(3\times 3) & +2I(3\times 3)\\
	O_{3\times 1} & +O_{3\times 1} & -O_{3\times 3} & +I_{3\times 3} & -I(3\times 3)\\
	O_{3\times 1} & -O_{3\times 1} & +O_{3\times 3} & -0_{3\times 3} & +I_{3\times 3}\\
	... & etc & ... & and\ so\ on & ...
	\end{array}\right]
$$
\begin{center}
\noindent \textbf{\textcolor{blue}{Example.12a} $\zeta_F^{-1}$.   The \textit{block presentation} of the  M{\"{o}}bius function matrix
$\mu = \zeta^{-1}$ for  ($1_F=1$ and $n_F=3$ for $n \geq 2$) the $N$ relative special sequence $\textbf{\textcolor{blue}{F}}$ constituting\textbf{ the} 
label sequence denominating cobweb poset associated to \textbf{$F$-KoDAG} Hasse digraph }
\end{center}

\vspace{0.2cm}

\noindent The   \textcolor{blue}{\textbf{code}} for this KoDAG is given by its KoDAG self-evident code-triangle of the \textcolor{blue}{\textbf{coding matrix}} $C(\mu_F)$ (a starting part of it shown below):

\vspace{0.2cm}
$$
	 C(\mu_F) = \left[\begin{array}{llllll}
	1 & -1 & +2 & -4 &  +8 & -16\\
	0 &  +1 & -1 & +2 &  -4 & +8\\
  0 &  -0 &  +1 & -1 &  +2 & -4\\
	0 &  +0 &  -0 &  +1 & -1 & +2\\
	0 &  -0 &  +0 &  -0 & +1  & -1\\
	. & . & . & . & .
	\end{array}\right]
$$

\vspace{0.3cm}
\noindent From Observation 2  we infer what follows as obvious.

\begin{observen}
\noindent Compare with the Remark 1. The block structure of  $\zeta$  and consequently the block structure of $\mu$ \textbf{for any graded poset} with finite set of minimal elements (including cobwebs) is of the type:
$$
\zeta	= \left[ \begin{array}{lllll}
I_1,	B_1... \\
& I_2,	B_2... \\
& & I_3,	B_3... \\
	& & ... \\
	& & & & I_n, B_n...
	\end{array} \right]
$$

$$
\mu	= \left[ \begin{array}{lllll}
I_1,	-B_1 ...\\
&  I_2, -B_2...\\
& & I_3, -B_3... \\
	& & ... \\
	& & &  I_n, -B_n...
	\end{array} \right]
$$

\noindent $n \in N \cup \{\infty\}, \zeta,\mu \in I(\Pi;R)$  where  $I_r = I_{r_F\times r_F}$ and $B_r = B(r_F\times (r+1)_F)$ as introduced by Observation 2..
\end{observen}

\vspace{0.3cm}

\noindent \textcolor{blue}{\textbf{Recall}}  then and note here up and below the block structure $\zeta$  and consequently the block structure of $\mu$ \textbf{for any graded poset $P$} with finite set of minimal elements (including cobwebs)  which is proprietary characteristic for any $\sigma \in I(P;R)$ where the ring  $R$= $2^{\left\{1 \right\}}$ , $Z_2=\left\{0,1\right\}$, $Z$ etc.

$$
	\sigma  = \left[\begin{array}{llllll}
	I_{1_F\times 1_F} & M(1_F \times 2_F) & M(1_F \times 3_F) & M(1_F \times 4_F) & M(1_F \times 5_F) & M(1_F \times 6_F)\\
	0_{2_F\times 1_F} & I_{2_F\times 2_F} & M(2_F \times 3_F) & M(2_F \times 4_F) & M(2_F \times 5_F) & M(2_F \times 6_F)\\
	0_{3_F\times 1_F} & 0_{3_F\times 2_F} & I_{3_F\times 3_F} & M(3_F \times 4_F) & M(3_F \times 5_F) & M(3_F \times 6_F)\\
	0_{4_F\times 1_F} & 0_{4_F\times 2_F} & 0_{4_F\times 3_F} & I_{4_F\times 4_F} & M(4_F \times 5_F) & M(4_F \times 6_F) \\
	... & etc & ... & and\ so\ on & ...
	\end{array}\right]
$$
where  in the case of  $\os$-natural $\zeta$ or $\zeta^{-1}$ matrices , with matrix elements from the ring  $R$= $2^{\left\{1 \right\}}$ , $Z_2=\left\{0,1\right\}$, $Z$ etc the rectangle \textbf{\textit{non-zero}} block matrices  $M(k_F \times (k+1)_F)$  denote corresponding connected \textbf{graded poset characteristic} $k_F \times (k+1)_F$ matrices.\\ 
Note then that $M(k_F \times (k+1)_F)_{r,s} =  c_{i,j,k}B(k_F \times (k+1)_F)_{i,j}$, $i=1,...,k_F$ and  $i=1,...,(k+1)_F$ where the rectangular "`zero-one"' $B(k_F \times (k+1)_F)$ matrices were introduced by the Observation 2. Consult Remark 1. - apart from the  motivating examples - for $i=1,...,k_F$ and  $i=1,...,(k+1)_F$ as the layer $\left\langle \Phi_k \longrightarrow \Phi_{k+1} \right\rangle$ variables.

\vspace{0.2cm}

\noindent \textbf{Note} now the \textcolor{red}{\textbf{important}} fact. The relation

$$M(k_F \times (k+1)_F)_{i,j} =  c_{i,j,k}B(k_F \times (k+1)_F)_{i,j},$$
where 
$$i=1,...,k_F ,\;\; i=1,...,(k+1)_F$$

\noindent does not fix uniquely \textbf{the  layer}   $\left\langle \Phi_k \longrightarrow \Phi_{k+1} \right\rangle$  \textbf{coding matrix} 
$C_{k,k+1}= (c_{i,j,k})$,  $i=1,...,k_F$, $i=1,...,(k+1)_F$ for $F$-denominated \textbf{arbitrary graded poset} - \textit{ except for} cobweb posets for which $B(k_F \times (k+1)_F)= I(k_F \times (k+1)_F)$. In order to delimit this layer coding matrix \textbf{\textit{uniquely}} we  define $en \:bloc$ the coding matrix $\textcolor{blue}{\textbf{C}}(\mu_F)$ for all layers.

\begin{defn}  $F$-graded poset $\left\langle \Phi,\mu_F \right\rangle$  \textcolor{blue}{\textbf{coding matrix}}  $\textcolor{blue}{\textbf{C}}(\mu_F)$.\\ 
Let $k,r,s \in N \cup \left\{\textcolor{blue}{\textbf{0}}\right\}$. Then we define $\textcolor{blue}{\textbf{C}}(\mu_F)$ via $\os$ originated blocks as follows:

$$\textcolor{blue}{\textbf{C}}(\mu_F) =  \left( \textcolor{red}{\textbf{c}}_{r,s} \right)$$
where $\textcolor{red}{\textbf{c}}_{r,s}$ are coding matrix elements for $F$-denominated cobweb poset, hence

$$\mu_F =\left([r=s]I_{r_F,r_F} + [s>r]\textcolor{red}{\textbf{c}}_{r,s}B(r_F \times s_F)\right),$$
and where
$$ c_{i,j,k} \equiv M(k_F \times (k+1)_F)_{i,j} = c_{i,j}B(k_F \times (k+1)_F)_{i,j} .$$
\end{defn}

\noindent thus the following identifications are self-evident: 
$$\left\langle \Phi,\mu_F \right\rangle   \equiv \left\langle \Phi,\zeta_F \right\rangle  \equiv  \left\langle \Phi,\leq \right\rangle \equiv  \left\langle \Phi,\textcolor{blue}{\textbf{C}}(\mu_F) \right\rangle  .$$

\vspace{0.2cm}

\noindent \textbf{Result:} $\textcolor{blue}{\textbf{C}}(\mu_F)$ as well as block sub-matrices $M(k_F \times (k+1)_F)= (c_{i,j,k})$ where $k \in N \cup \left\{\textcolor{blue}{\textbf{0}}\right\}$ are defined i.e are given unambiguously.

\vspace{0.2cm}

\noindent Specifically, in \textbf{cobweb posets case}: for $\zeta$ function (matrix) we have   $M(k_F \times (k+1)_F) = I(k_F \times (k+1)_F)$, while for $\zeta^{-1}= \mu$ M{\"{o}}bius function (matrix) - from already considered  examples' prompt we have already deduced these unambiguous $\textcolor{red}{\textbf{c}}_{r,s}$ ( see Theorem 2 for cobweb posets - above). Namely :

$$M(r_F \times (r+1)_F) = \textcolor{red}{\textbf{c}}_{r,r+1}I(r_F \times (r+1)_F).$$

\noindent What about any $F$-denominated  graded posets then? The answer \textbf{\textit{now}} is of course secured now to be the same as for $F$-cobweb posets. The answer is automatically secured by the Definitions 7,8 . Just replace in the above Theorem 3 for cobweb posets $I(r_F \times (r+1)_F)$  by $B(r_F \times (r+1)_F)$ and-or see the Theorem 4 below for the corresponding recurrence - the recurrence equivalent to the recurrence relation definition of $\textcolor{red}{\textbf{c}}_{r,s}$.

\vspace{0.2cm}

\noindent In order to be complete also with the next section content  another important example - the example of cover relation $\kappa_{\Pi} \in I(\Pi,R)$ matrix follows. Recall for that purpose now Observation 1 and the Remark 1 as to conclude what follows.

\vspace{0.2cm}

\begin{observen}($n \in N \cup \{\infty\}$)
\noindent The block structure of  cover relation $\kappa_{\Pi} \in I(\Pi,R)$  ($\chi\left( \prec\cdot_{\Pi} \right) \equiv \kappa_{\Pi},$) is the following

$$\kappa_{\Pi} = \os_{k=1}^n \kappa_k = $$

$$
	= \left[ \begin{array}{lllll}
	0_{1_F\times 1_F} & I(1_F\times 2_F) & 0_{1_F\times \infty}\\
	0_{2_F\times 1_F} & 0_{2_F\times 2_F} & I(2_F\times 3_F) & 0_{2_F\times \infty}\\
	0_{3_F\times 1_F} & 0_{3_F\times 2_F} & 0_{3_F\times 3_F} & I(3_F\times 4_F) & 0_{3_F\times \infty}\\
	& & ... \\
	0_{n_F\times 1_F} & ...& 0_{n_F\times n_F} & I(n_F\times (n+1)_F) & 0_{n_F\times \infty}
	\end{array} \right]
$$
\vspace{0.1cm}
\noindent where $\kappa_k$ is a cover relation of di-biclique $\langle \Phi_k\to\Phi_{k+1}\rangle$, $I_k \equiv I(k_F \times (k+1)_F)$, $k = 1,...,n$ and where  - recall - $I (s\times k)$  stays for $(s\times k)$  matrix  of  ones  i.e.  $[ I (s\times k) ]_{ij} = 1$;  $1 \leq i \leq  s,  1\leq j  \leq k.$  while  $n \in N \cup \{\infty\}$.  

\vspace{0.1cm}

\noindent and consequently the block structure of \textbf{reflexive cover} relation $\eta_{\Pi} \in I(\Pi,R)$  ($ \chi\left( \leq \cdot_{\Pi} \right)= \prec\cdot_{\Pi} +\delta \equiv \eta_{\Pi}$) is given by
$$
	= \left[ \begin{array}{lllll}
	I_{1_F\times 1_F} & I(1_F\times 2_F) & 0_{1_F\times \infty}\\
	0_{2_F\times 1_F} & I_{2_F\times 2_F} & I(2_F\times 3_F) & 0_{2_F\times \infty}\\
	0_{3_F\times 1_F} & 0_{3_F\times 2_F} & I_{3_F\times 3_F} & I(3_F\times 4_F) & 0_{3_F\times \infty}\\
	& & ... \\
	0_{n_F\times 1_F} & ...& I_{n_F\times n_F} & I(n_F\times (n+1)_F) & 0_{n_F\times \infty}
	\end{array} \right]
$$
\end{observen}

\vspace{0.2cm}

\noindent Specifically, \textbf{if} restricting to   \textbf{cobweb posets}: for $\zeta$ function (matrix) we have   $B(k_F \times (k+1)_F) = I(k_F \times (k+1)_F)$, while for $\zeta^{-1}= \mu$ M{\"{o}}bius function (matrix) we would expect 

$$B(r_F \times (r+1)_F) = c_{r,r+1}I(r_F \times (r+1)_F)$$
\noindent where  $c_{k,k+1}= [C(\mu_F)]_{k,(k+1)}$.

\vspace{0.2cm}

\noindent What is then the explicit formula for $c_{k,k+1} ?$ It is of course equivalent to the question: what is then the explicit formula for $c_{r,s}?$  Let us recapitulate our experience till now in order to infer the closing answer Theorem 4.  and its equivalent proof method.

\vspace{0.2cm}

\noindent \textcolor{red}{\textbf{Training in relabeling - \textit{Exercise}.}}
 
\vspace{0.1cm}

\noindent As we were and are to compare formulas from papers using different labeling - write and learn to see formulas from the above and below Observations as for $x,y,k,s \in N \cup \left\{\textcolor{blue}{\textbf{0}}\right\}$  on one hand and as for $x,y,k,s \in \textcolor{red}{\textbf{\textsl{N}}}$ on the other hand. Because of the comparisons  reason 
we shall tolerate and use both being indicated explicitly if needed.

\vspace{0.4cm}

\noindent \textbf{Recapitulation 4.1. ; notation and \textcolor{blue}{The Formula.}} The code $C(\mu_F)$ matrix no more secret.

\vspace{0.2cm}

\noindent \textbf{Notation}. Upside down notation development continuation. 

\vspace{0.2cm}

\noindent Recall: $$n^{\overline{k}} = n(n+1)(n+2)...(n+k-1) ,$$
Denote:
$$n_F^{\overline{k}} \equiv n_F(n+1)_F(n+2)_F...(n+k-1)_F $$

\noindent Denote (valid whenever defined for corresponding functions $f$ of the natural number argument or of an argument from any chosen ring ):

$$
f(r_F)^{\overline{k}} = f(r_F)f([r+1]_F)...f([r + k -1]_F), \ n^{\overline{0}} \equiv 1,\: n \in N \cup \ \left\{ \textit{\textbf{0}} \right\},Z,R,etc. , 
$$

$$
f(r_F)^{\underline{k}} = f(r_F)f([r-1]_F)...f([r - k +1]_F) , \ n^{\underline{0}} \equiv 1,\: n \in N \cup \ \left\{ \textit{\textbf{0}} \right\},Z,R,etc. .
$$

\vspace{0.2cm}

\noindent Define Krot-on-shift-functions  $K_s ,\:s,r,i \in N \cup \left\{0\right\}$ or \textcolor{blue}{\textbf{Kroton}} functions in brief -(Kroton = Croton = Codiaeum).

\begin{defn} ( $\textcolor{red}{\textbf{N }}\cup \left\{\textcolor{blue}{\textbf{0}}\right\}$ \ \textbf{labels})
$$
K_s(r_F) = [s>r] [(r+1)_F - 1]^{\overline{s-r}}
$$
\end{defn}

\noindent  These of course constitute an upper triangle matrix with zeros on the diagonal for $s,r \in N \cup \left\{0\right\}$, (\textcolor{red}{\textbf{r}} = labels \textcolor{red}{\textbf{r}}ows). 

\vspace{0.2cm}

Note  two cases:\\

Let $s-r-1 \neq 0$. Then

$$K_s(r_F) = [s>r]\prod_{i=r+1}^{s-1}(i_F - 1)$$

Let $s-r-1 = 0$. Then

$$K_s(r_F) = [s>r] .$$

\vspace{0.2cm}

\noindent Now - with this $\textcolor{red}{\textbf{N }}\cup \left\{\textcolor{blue}{\textbf{0}}\right\}$ labeling as established in this note (Remark2.1.) - perform simple calculations. \textcolor{red}{\textbf{Fibonacci}} sequence $F = \left\langle \textcolor{red}{\textbf{1}},1,2,3,5,8,13,21,34,...\right\rangle$ case \textcolor{red}{\textbf{Example}}.

\vspace{0.1cm}

\noindent  $K_2(1_F) = 1 $ , $K_s(1_F) = 0 $ for $s>2$;

\noindent  $K_3(2_F) = 1 $ , $K_s(2_F) = 0 $ for $s>3$;  

\noindent  $K_4(3_F) = 1 $ , $K_5(3_F) = 1 $ ,  $K_6(3_F) = 2$, $K_7(3_F) = 2 \cdot 4 = 8$, $K_8(3_F) = 8 \cdot 12 = 96$, $K_9(3_F) = 96 \cdot 20 = 1920$ , and so on, 

\noindent  $K_5(4_F) = 1 $ , $K_6(4_F) = 1\cdot 4 $ ,  $K_7(4_F) = 4 \cdot 7 = 14 $, $K_8(4_F) = 14 \cdot 12 = 168$, $K_9(4_F) = 168 \cdot [F_8 -1] = ?$, $K_{10}(4_F) = 3360 \cdot [9_F - 1] = ?$ , and so on. Note that in the course of the above the following was used  ( $N \cup \left\{0\right\}$ - labeling).

\vspace{0.3cm}

\noindent \textbf{Lemma 4.1} ($r,s \in N \cup \left\{0\right\}$ . Obvious)

$$K_{s+1}(r_F) = K_s(r_F)\bullet [s_F -1], \ \ K_{r+1}(r_F) = 1,$$

\vspace{0.3cm}

\noindent \textcolor{blue}{\textbf{N}} sequence  case \textcolor{blue}{\textbf{Example}}. This exercise has obvious outcomes in view of the Lemma 2.1. For the just check results see absolute values of 
\textcolor{blue}{\textbf{coding matrix}} matrix elements from the Example 9. .


\vspace{0.2cm}

\noindent The next fact we mark as Lemma because of its importance.

\vspace{0.2cm}

\noindent \textbf{Lemma 4.2} (Obvious - recapitulation.)\\
Let  $R=N,Z$,...,any commutative ring. For any graded $F$-denominated poset (hence connected) i.e for any chain of subsequent natural joins of bipartite digraphs (di-bicliques for KoDAGs) and with the linear labeling of nodes fixed ( $s,r \in N \cup \left\{0\right\}$  as in Remark 2.1.  or $s,r \in N$) : 

$$\mu = \left( \delta_{r,s}I_{r_F\times r_F} + [s>r]C(\mu_F)_{r,s}B(r_F \times s_F) \right)$$
where $C(\mu_F)_{r,s} \in R$  are given by Definition 6. while $B(r_F \times s_F)$ are nonzero matrices introduced in the Observation 2.


\vspace{0.2cm}

\noindent Bearing in mind Definitions 7 and 8 and the the above Lemma 4.2. we see that the Theorem 3 for cobweb posets extends to be true for all $F$-denominated posets.

\vspace{0.2cm}

\noindent \textbf{Theorem 4} (Kwa\'sniewski)

\vspace{0.1cm}

\noindent Let $F$ be \textbf{any} natural numbers valued sequence. Then \textcolor{red}{\textbf{for arbitrary}} $F$-denominated graded poset (cobweb posets included)    

$$ C(\mu_F)_{r,s}= c_{r,s} = [r=s] + K_s(r_F)(-1)^{s-r} =[r=s] + [s>r] (-1)^{s-r}[(r+1)_F - 1]^{\overline{s-r}} $$ 
i.e.

$$C(\mu_F)_{r,s}= c_{r,s} = [s\geq r] (-1)^{s-r}[(r+1)_F - 1]^{\overline{s-r}}, $$
i.e. 
$$C(\mu_F)_{r,s}= c_{r,s} = [s\geq r] (-1)^{s-r}K_s(r_F), $$

\noindent with matrix elements from $N$ or the ring  $R$= $2^{\left\{1 \right\}}$ , $Z_2=\left\{0,1\right\}$, $Z$ etc.\\
i.e. for cobweb posets

$$\mu = \left( \delta_{r,s}I_{r_F\times r_F} + (-1)^{s-r}K_s(r_F)I(r_F \times s_F)   \right)$$ 
i.e.

$$
	\mu  = \left[\begin{array}{llllll}
	I_{1_F\times 1_F} & c_{1,2}I(1_F \times 2_F) & c_{1,3}I(1_F \times 3_F) & c_{1,4}I(1_F \times 4_F) & c_{1,5}I(1_F \times 5_F) & c_{1,6}I(1_F \times 6_F)\\
	0_{2_F\times 1_F} & I_{2_F\times 2_F} & c_{2,3}I(2_F \times 3_F) & c_{2,4}I(2_F \times 4_F) & c_{2,5}I(2_F \times 5_F) & c_{2,6}I(2_F \times 6_F)\\
	0_{3_F\times 1_F} & 0_{3_F\times 2_F} & I_{3_F\times 3_F} & c_{3,4}I(3_F \times 4_F) & c_{3,5}I(3_F \times 5_F) & c_{3,6}I(3_F \times 6_F)\\
	0_{4_F\times 1_F} & 0_{4_F\times 2_F} & 0_{4_F\times 3_F} & I_{4_F\times 4_F} & c_{4,5}I(4_F \times 5_F) & c_{4,6}I(4_F \times 6_F)\\
	... & etc & ... & and\ so\ on & ...
	\end{array}\right]
$$
where  $I(k_F \times (k+1)_F)$  denotes (recall)  $k_F \times (k+1)_F$ matrix  of all entries equal to one. \textbf{For any} $F$-\textbf{denominated poset} \textcolor{red}{\textbf{replace}} $I(k_F \times (k+1)_F)$ \textcolor{red}{\textbf{by}} $B(k_F \times (k+1)_F)$ obtained from  $I(k_F \times (k+1)_F)$ via replacing adequately (in accordance with Hasse  digraph) corresponding \textit{ones} by \textbf{zeros}.

\vspace{0.2cm}

\noindent \textit{Another Proof }:  One may prove the above also as follows.

\vspace{0.1cm}

\noindent  From  motivating examples we know that $\mu (x_{r,i},x_{s,j})= \mu (x_r,x_s)$. Observe then how the recurrent definition of  M{\"{o}}bius function matrix $\mu$  gives birth to daughter descendant of  $\mu$ i.e. the block structure of  M{\"{o}}bius function coding matrix $C(\mu)$ implying for $C(\mu)$  a recurrence  allowing  simple solution simultaneously with combinatorial interpretation of  \textcolor{blue}{\textbf{Krot}}on matrix  $K =(K_s(r_F))\equiv (K_{r,s})$ , where  
$ K_s(r_F) = \left|C(\mu)_{r,s}\right|$.

\vspace{0.1cm}

\noindent For that to do call back the recurrent definition of the M{\"{o}}bius function where $x,y \in \Phi$  for $\Pi = (\Phi, \leq)$ and where - note: $\mu (x,y) = -1 $ for $x\prec\cdot y $ :

$$ \mu (x,y)=\Big\{\begin{array}{l}\;\;\;\;1\;\;\;\;\;\;\;\;\;\;\;\;\;\;\;\;\;\;\;\;\;\;\;\;\;x=y\\-\sum_{x\leq z <y} \mu (x,z),\;\; x<y\end{array}.$$

\noindent  The above recurrent definition M{\"{o}}bius function becomes - after \textbf{linear} order labeling has been applied - either  $r,s,i \in N \cup \left\{0\right\}$ - as fixed-stated in this note, Remark 2. and/or  fact that $r,s \in N$ - whereby $r,s$ are block-row and block-column indexes correspondingly -  say it again - the above recurrent definition M{\"{o}}bius function 
in the case of $F$-denominated graded posets becomes ( $c_{r,r+1} = -1$ )

$$ c_{r,s}=\Big\{\begin{array}{l}\;\;\;\;1\;\;\;\;\;\;\;\;\;\;\;\;\;\;\;\;\;\;\;\;\;s=r\\-\sum_{r\leq i < s} c_{r,i} ,\;\; r<s\end{array}.$$

\vspace{0.1cm}

\noindent For that to see \textbf{ note that} $\forall x,y,z \in \Phi$,  $\exists \; r,s,i \in N$ such that  $x_r \in \Phi_r$ , $y_s \in \Phi_s $,  $z_i \in \Phi_i $, hence  for $x_r < y_s \:\equiv\: r<s$  where (\textcolor{red}{\textbf{Important}}!)  $r,s,i$ stay now for \textbf{labels of independent sets }(levels)  $\left\{\Phi_k\right\}$ i.e. label steps of La Scala i.e. label blocks. Thereby 
$$  c_{r,s} = \mu (x_r,y_s) = - \sum_{x_r\leq z <y_s} \mu (x_r,z)  = - \sum_{x_r\leq z_i <y_s} \mu (x_r,z_i) = \sum_{r\leq i < s}\ c_{r,i}. $$

\noindent (Bear  in mind  Lemma 2.2. in order to get back to $\mu$ matrix unblocked appearance if needed.)  From this recurrence the thesis follows.


\vspace{0.2cm}

\noindent How does this happens?  \textbf{1)} Let us put  $r=1$ just for the moment in order to make an inspection via example  ($r$ stays for \textbf{\textit{block - row}} label and  $k>1$) and  \textbf{2)} use the Russian babushka in Babushka inspection i.e. apply the recurrent relation above subsequently till the end - till the smallest of size 1  babushka is encountered which is here $c_{r,r+1}=-1$ . Use then trivial induction to state the validity of what follows  below for all relevant values of variables $r,s \in N.$ 

$$ c_{1,k}= -\sum_{1\leq i <k} c_{1,i} =  \left(-\sum_{1\leq i <k_F}\right)\; \left(-\sum_{1\leq i <(k-1)_F}\right)...\left(-\sum_{1\leq i <3_F} \right) \; c_{1,2},$$ 
i.e.
$$ c_{1,k}= (-1)^{k-1}\left(\sum_{1 \leq i < k_F}\right)...\left(\sum_{1 \leq i < 4_F} \right) \;\left(\sum_{1 \leq i < 3_F}\right) (+1),$$
i.e.
$$ c_{1,k}= - [1+1=k] + [k>2] (-1)^{k-1}  \left(k_F - 1)\right)...\left(3_F - 1) \right) \; (+1) = $$
$$ =  - [1+1=k] + [k>2](-1)^{k-1}\; \prod_{i=2+1}^k(i_F - 1).$$
Similarly we conclude that now for arbitrary $r,s \in N$

$$ c_{r,s}= [s=r]  - [s=r+1] + [s>r+1] (-1)^{s-r}\left(s_F - 1)\right)...\left(3_F - 1) \right) \; (+1) =$$
$$ = [s=r] - [s=r+1] + [s>r+1](-1)^{s-r}\; \prod_{i=r+2}^s(i_F - 1),$$
Equivalently  we conclude that now for arbitrary $r,s \in N \cup \left\{0\right\} $

$$ c_{r,s}= [s=r]  - [s=r+1] + [s>r+] (-1)^{s-r}\left((s-r-1)_F - 1)\right)...\left(3_F - 1) \right) \; (+1) =$$
$$ = [s=r] - [s=r+1] + [s>r+1](-1)^{s-r}\; \prod_{i=r+1}^{s-1}(i_F - 1),$$

\vspace{0.2cm}

\noindent \textbf{To colligate and to imagine hint.} Starting from the left upper corner of La Scala of  $\zeta$, $\mu$,...,$\sigma \in I(\Pi,R)$  \textcolor{red}{\textbf{down}} $\Downarrow$ is biunivoquely starting from the "`bottom"' or "`root"' minimal elements level $\Phi_0$  \textcolor{green}{\textbf{up}}  $\Uparrow$ the Hasse digraph $(\Pi,\prec\cdot)$  uniquely representing  the "`much, much more cobwebbed tree'"' - the digraph  $(\Pi,\leq )$ 

\vspace{0.2cm}

\noindent \textsl{\textbf{Descriptive - combinatorial interpretation}}:  Once the formula has been observed-derived as above the following turns out perceptible. Namely note  that \\ 
\textbf{1.}  for  $F=N$,  $[s\neq r  ]$, the Kroton matrix element  $\left|\textcolor{blue}{\textbf{C}}(\mu_N)_{r,s}  \right|$, where 

$$ \textcolor{blue}{\textbf{C}}(\mu_N)_{r,s}= \textcolor{blue}{\textbf{c}}_{r,s} = [s>r] (-1)^{s-r}[(r+1)_N - 1]^{\overline{s-r}} $$

\noindent is equal to the number of heads' dispositions of maximal chains tailed at \textbf{one} vertex of the  $r-th$  level  and headed up  at \textbf{one} vertex of the $s$-th level. This biunivoquely corresponds to the number of summands $=\left|\textcolor{blue}{\textbf{C}}(\mu_N)_{r,s} \right|$  entering the recurrence calculation of the $\textcolor{blue}{\textbf{C}}(\mu_N)$ matrix ("`the Russian babushka in Babushka introspection"' with interchangeable signs) being in one to one correspondence with climbing up Hasse digraph i.e. descending down the matrix $\mu$ La Scala along the way uniquely encoded by  the subjected to their \textcolor{red}{\textbf{heads}} disposition maximal chains

$$c=<x_{\textcolor{blue}{\textbf{r}}},x_{r+1},...,x_{s-1},x_{\textcolor{red}{\textbf{s}}}>, \: x_i \in \Phi_i, \:i=\textcolor{blue}{\textbf{r}},r+1,...,s-1,\textcolor{red}{\textbf{s}}$$ 

\noindent  with the tail  $\textcolor{blue}{\textbf{r}}$  and  the head $\textcolor{red}{\textbf{s}}$ fixed as start and the end points of the descending down the La Scala blocks trip
 ( $\equiv $  climbing up the levels of the  graded Hasse  digraph  $\left\langle \Phi,\prec\cdot\right\rangle$).

\vspace{0.2cm}

\noindent \textbf{2.} For the same interpretation in the general $F$-case \textbf{ apply} the Upside Down Notation Principle.

\vspace{0.2cm}

\noindent According to and from  the above one extracts the obvious now property of \textcolor{blue}{\textbf{Krot}}on functions i.e. matrix elements of \textcolor{blue}{\textbf{Krot}}on matrix  $K =(K_s(r_F))\equiv (K_{r,s})$ 

\vspace{0.2cm}

\noindent \textbf{Lemma 4.3} ($r,s \in N \cup \left\{0\right\}$ . 

$$K_{s+1}(r_F) = K_s(r_F)\bullet [s_F -1], \ \ K_{r+1}(r_F) = 1$$
is equivalent to 

$$ K_{r,s}= -\sum_{r\leq i < s} (-1)^{s-i} K_{r,i} \ \ K_{r+1}(r_F) = 1, \ \ s>r.$$

\noindent \textbf{Exercise.} Deliver the descriptive combinatorial interpretation of \textcolor{blue}{\textbf{Krot}}on matrix in the language of  hyper-boxes from [11].

\vspace{0.2cm}

\noindent \textcolor{blue}{\textbf{Compare.}}  All the above may be now compared with   \textbf{1.4. Cobweb posets' M\"{o}bius function} where it has been proved  that  for  $x \in \Phi_r$ ,   $z \in \Phi_s$  and  for     
$$s>r , \ [x,z] = x\oplus \Phi_{r+1}\oplus...\oplus \Phi_{s-1}\oplus  z$$    

  $$\mu(x,y) = (-1)^{s-r}\prod_{k=r+1}^{s-1}[k_F-1]. $$

\vspace{0.1cm}



\noindent \textcolor{blue}{\textbf{Whitney numbers.}}  Let us remind the notation: $r_F = \left|\Phi_r\right|$. Let us then recall (see [49]) definitions of  Whitney numbers of the first kind $w_r(P)$ and Whitney numbers of second kind $W_r(P)$ where $P$ is any given graded poset with bounded independence sets - i.e.  $\left|\Phi_r\right| \in N$ for $r \in N \cup \left\{0\right\}$  and we assume that $\left|\Phi_0\right|= 1 = 0_F$  hence  $\Phi_0 = \left\{0\right\}$; $0 \in \P$ stays for minimal element of the poset $P$. Now here are these definitions. 

$$w_{r}(P)=\sum_{x\in P,\,r(x)=r}\mu(0,x),$$
$$W_r(P)=\sum_{x\in P,\,r(x)=r}1=|\{x\in P\,:\,r(x)=r\}|.$$

\vspace{0.2cm}

\noindent It is obvious just by notation  that 

$$W_r(P)= r_F .$$
Of course in the general case of finite posets $P_n$ the number $\left|\Phi_k\right|$ might depend on $n$ and then we end up with an array $\left(W_k(P_n)\right)$ of Whitney numbers as it is the case with binomials or Gaussian binomials for example - (see further classical examples in [49]).

\vspace{0.2cm}

\noindent According to Theorem 4 we have for \textcolor{blue}{\textbf{cobweb posets}} i.e. for $\Pi$'s

$$C(\mu_F)_{0,s}= c_{0,s} = [s\geq 0] (-1)^s(1_F - 1)^{\overline s}, $$
or equivalently - as the values of $\mu(x,y)$ depend only on the rank of its arguments  
$$\mu_F(0,x)=  [r(x)\geq 0] (-1)^{r(x)}(1_F - 1)^{\overline {r(x)}}, $$
or equivalently (compare all this with (3) in [28])
$$\mu_F(0,x)=  [x=0] - [r(x)=1](1_F - 1) +  [r(x)>1](-1)^{r(x)}\prod_{k=1}^{r(x)-1}(k_F-1), $$
or equivalently  just 
$$\mu_F(0,x)=  [r(x)\geq 0] (-1)^{r(x)}K_r(0_F). $$

\vspace{0.1cm}

\noindent Consequently - as the values of $\mu(x,y)$ depend only on the rank of its arguments - the Whitney numbers of the first kind for the denominated by $F$ 
\textcolor{blue}{\textbf{cobweb poset}} $\Pi$ may be calculated along the formula

  $$w_r(\Pi) =  \sum_{\{x\in\Pi\,:\,r(x)=r\}}\mu_F(0,x)= r_F \cdot\mu_F(0,x)$$
 i.e.
   
   $$w_r(\Pi) = r_F\cdot(-1)^r K_r(0_F).$$
Naturally $w_0(\Pi)=1$. Compare  the above with (4) in [28].  

\vspace{0.2cm}

\noindent Of course in the \textit{general case} of finite posets $P_n = \bigcup_{k=0}^n \Phi_k(n) $ the number $\left|\Phi_k(n)\right|$ might depend on $n$ and then we end up with an array $\left(w_k(P_n)\right)$ of Whitney numbers of the first kind as it is the case with binomials or Gaussian binomials for example - (see further classical examples in [49]).

\vspace{0.2cm}

\noindent To this end - consequently - let us consider characteristic polynomials  $\chi_{P_n}(t)$, $n\geq 0$ 
defined as ([50,51],[28])

$$\chi_{P_n}(t)=\sum_{x\in P_n}\mu(0,x)t^{n-r(x)}=\sum_{k = 0}^n w_k(P_n)t^{n-k}.$$
The formula for characteristic polynomials - here for \textbf{specific} $F$-denominated  finite cobweb sub-posets    $P_n = \oplus_{k=0}^n \Phi_k$   (i.e. $\left|\Phi_k(n)\right|$ does not depend on $n$)- namely - this  formula for characteristic polynomials  obviously is of the form 

$$ \chi_{P_n}(t)=\sum_{k=0}^n (-1)^k \cdot k_F \cdot x^{n-k} \cdot K_k(0_F) $$
or equivalently  (compare with Theorem 3.1 in [28])

$$ \chi_{P_n}(t)=x^n - x^{n-1}1_F(1_F-1)  +  \sum_{k=2}^n (-1)^k \cdot k_F \cdot x^{n-k} \cdot \prod_{r=1}^{k-1}(r_F-1).$$

\vspace{0.1cm}

\noindent \textbf{Recapitulation 4.2.} \textcolor{red}{natural join}.

\vspace{0.2cm}

\noindent \textbf{Recall} that  both $\leq$ partial order  and  $\prec\cdot$ cover relations are \textcolor{red}{\textbf{natural join}} of their bipartite correspondent chains, and this is 
exactly the reason and the very source of the Theorem 2  validity and shape. This is also the obvious  clue statement for what follows. Note also that all on structure of any $P$ poset's information is 
coded by the  $\zeta$ matrix - a characteristic function of $\leq \in \  P = \left\langle \Phi,\leq\right\rangle$. In short: $\zeta$ and equivalently $\mu = \zeta^{-1}$ are the Incidence algebra of $P$ 
coding elements. In brief - recall -  the following identifications are self-evident:

$$\left\langle \Phi,\mu_F \right\rangle   \equiv \left\langle \Phi,\zeta_F \right\rangle  \equiv  \left\langle \Phi,\leq \right\rangle \equiv  \left\langle \Phi,\textcolor{blue}{\textbf{C}}(\mu_F) \right\rangle .$$



\section{$F$-nomial coefficients and  $[Max]$ matrix of the $N$ weighted reflexive reachability relation }

\vspace{0.2cm}

\noindent  Call back now the \textcolor{red}{\textbf{Remark 1}}. Then consider the  incidence algebra  of the \textbf{cobweb poset} $\Pi$  as the algebra over (simultaneously) the ring $R$ and the Boolean algebra $2^{\left\{1\right\}}$. Denote this incidence algebra by  $I(\Pi,R, 2^{\left\{1\right\}})$).

\vspace{0.2cm}
\noindent In the case $R= 2^{\left\{1\right\}}$  denote it by  $I(\Pi,2^{\left\{1\right\}})\equiv I(\Pi,2^{\left\{1\right\}}, 2^{\left\{1\right\}}).$ Then for $\zeta \in I(\Pi,2^{\left\{1\right\}})$ we have of course $\zeta^{-1}= \zeta$ ("`reflexive reachability"'),\ $\zeta_{\leq \cdot}^{-1}  =  \zeta_{\leq \cdot}$ \ (reflexive "`cover"') and so on. This is of course true for any poset relevant algebra i.e. for $I(P,2^{\left\{1\right\}})$  - graded posets with finite set of minimal elements - included.

\vspace{0.2cm}

\noindent  Consider now the algebra $I(\Pi, \textcolor{red}{\textbf{Z}}, 2^{\left\{1\right\}}))$. We shall define now another characteristic matrix  $[Max]$ as the matrix of the "` $N$ weighted"' reflexive reachability relation. For that to do  recall that in case of    $I(\Pi,2^{\left\{1\right\}})$ 

\vspace{0.2cm}
 
\begin{center}
$ \leq \:= \:\prec\cdot^* $ = reflexive reachability of $\prec\cdot$ \\
\vspace{0.2cm} 
$ \prec\cdot^* \equiv {(I - \prec\cdot)}^{-1} = \prec\cdot^{0@} + \prec\cdot^{1@} + \prec\cdot^{2@} +...+ \prec\cdot^{k@} + ...  \ \equiv \ \bigcup_{k\geq 0}\prec\cdot^k,$  
\end{center}

\vspace{0.1cm}

\noindent where binary relations $\leq \  \subset \ \Phi \times \Phi $ and $\prec\cdot \ \subset \  \Phi \times \Phi $  etc.  as subsets  are identified with their matrices (see [3,2]), 
for example  $\prec\cdot \equiv \kappa $. In the above the Boolean powers of $\kappa$ were in action while here below this are to be powers over the  $R = N, Z, 2^{\left\{1\right\}}$, etc.

\begin{defn}
 The $[Max]$ matrix of the $N$ weighted reflexive reachability relation is defined by the over the ring $Z$ power series formula

$$[Max] = {(I - \prec\cdot)}^{-1}= \prec\cdot^0 + \prec\cdot^1 + \prec\cdot^2 +...+ \prec\cdot^k + ... = \sum_{k\geq 0} \kappa^k = {(I - \kappa)}^{-1}$$

\end{defn}

\vspace{0.2cm}

\noindent  Naturally  

$$[Max]^{-1} = \delta - \kappa = 
 = \left[ \begin{array}{llllll}
I_1 &	-B_1 & zeros \\
& I_2 &	-B_2 & zeros \\
& & I_3 &	-B_3  & zeros \\
	& & ... & \\
	 & & & I_n & -B_n & zeros
	\end{array} \right]
$$
where (recall from Section I. 1.5 )
$$
	[Max]_F = \mathbf{A}_F^0 + \mathbf{A}_F^1 + \mathbf{A}_F^2 + ...= (1 - \mathbf{A}_F)^{-1}=
$$
$$
= \left[\begin{array}{llllll}
	I_{1_F\times 1_F} & B(1_F \times 2_F) & B(1_F \times 3_F) & B(1_F \times 4_F) & B(1_F \times 5_F) & ... \\
	0_{2_F\times 1_F} & I_{2_F\times 2_F} & B(2_F \times 3_F) & B(2_F \times 4_F) & B(2_F \times 5_F) & ... \\
	0_{3_F\times 1_F} & 0_{3_F\times 2_F} & I_{3_F\times 3_F} & B(3_F \times 4_F) & B(3_F \times 5_F)  & ...\\
	0_{4_F\times 1_F} & 0_{4_F\times 2_F} & 0_{4_F\times 3_F} & I_{4_F\times 4_F} & B(4_F \times 5_F) & ... \\
	... & etc & ... & and\ so\ on & ...
	\end{array}\right].
$$

\vspace{0.2cm}

\noindent \textbf{Comment 6.} Combinatorial interpretation of  $[Max]$.

\begin{center}
$[Max]_{s,t}$ = the number of all maximal chains in the poset interval $[x_{s,i},x_{t,j}] = [x_s,x_t] \equiv  [s,t].$
\end{center}

\vspace{0.1cm}

\noindent where $x_{s,i},x_s \in \Phi_s$  and $x_{t,j},x_t \in \Phi_t$ for , say , $s \leq t$ with the reflexivity (loop) convention adopted i.e. $[Max]_{t,t}=1$.

\noindent The above obvious statement being taken into the account, in view and in conformity with the environment of the Theorem 1 we arrive at the trivial and  powerful Theorem 5.

\vspace{0.2cm}

\noindent \textbf{Theorem 5.}

\vspace{0.1cm}

\noindent Consider any $F$-cobweb poset with $F$ being a natural numbers valued sequence. Let  $x_k \equiv k \in \Phi_k$ and $x_t \equiv t \in \Phi_n$. Then 

$$  \sum_{i \in \Phi_n}[Max]_{k,i} \equiv   \sum_{i = 1}^{n_F}[Max]_{k,i} =  \left|C_{max}\langle\Phi_{k+1} \to \Phi_n \rangle \right| = n^{\underline{m}}_F,$$

\noindent where  $m = n-k$.

\vspace{0.2cm}

\noindent \textcolor{red}{\textbf{Note}} that  $k,m,n$ are level labels (vertical) while $i = 1,...,n_F $ stays for  horizontal - along the fixed level - label. With this in mind fixed we observe what follows.

\vspace{0.2cm}

\noindent \textbf{Corollary 5.1.}.

\vspace{0.1cm}

\noindent Consider any $F$-cobweb poset with $F$ being a \textbf{cobweb admissible} sequence.

\noindent Let  $x_k \equiv k \in \Phi_k$ and $x_n \equiv n \in \Phi_n$. Let $n \geq k \equiv(n-m) \geq 2$. Then 

$$  [Max]_{k,n} \left|\Phi_n\right|= n^{\underline{m}}_F $$
i.e.

$$  [Max]_{k,n} = \fnomial{n-1}{k-2}(n-k+1)_F! $$

\vspace{0.2cm}

\noindent \textbf{Corollary 5.2.}  \textit{ colligate with heads dispositions allied to the Theorem 3}. 

\vspace{0.1cm}

\noindent Consider any $F$-cobweb poset with $F$ being a \textbf{cobweb admissible} sequence.

\noindent Let  $x_k \equiv k \in \Phi_k$ and $x_m \equiv n \in \Phi_n$. Let $l+1= n \geq k \equiv(n-m) \geq 2$. Then 

$$  [Max]_{k,n} \left|\Phi_n\right|= n^{\underline{m}}_F $$
i.e.

$$  \fnomial{n-1}{n-1-k}(n-1-k)_F! = \fnomial{n-1}{k}(n-1-k)_F! =  [Max]_{k-2,n} $$

$$  \fnomial{n-1}{n-1-k} = \frac {[Max]_{k-2,n}}{(n-1-k)_F!} $$
i.e. ($n-1=l$)

$$  \fnomial{l}{k} = \fnomial{l}{l-k} = \frac {[Max]_{k-2,l+1}}{(l-k)_F!} $$

\vspace{0.2cm}

\noindent \textcolor{red}{\textbf{Note}} that  $k,m,n,l$ are level labels (vertical) and this is convention to be kept till the end of this note. 
\vspace{0.2cm}

\noindent The above obvious statement being taken into the account, in view and in conformity with the environment  of  Theorems 1 and 3 we are prompt to extract the trivial and powerful statement as the Theorem 6.

\vspace{0.2cm}

\noindent \textbf{Theorem 6.}   

\vspace{0.1cm}

\noindent Consider any $F$-cobweb poset with $F$ being a \textbf{cobweb admissible} sequence.  Let  $x_k \equiv k \in \Phi_k$ and $x_m \equiv n \in \Phi_n$. Let $(l+1) \geq k  \geq 2$. Then 

$$  \fnomial{l}{k} = \fnomial{l}{l-k} = \frac {[Max]_{k-2,l+1}}{(l-k)_F!} $$
i.e. 

\begin{center}
$\fnomial{l}{k}$ = $(l-k)_F!$'th \  fraction  of the number of all maximal chains  in the poset interval\  $[x_{k-2},x_{l+1}],$
\end{center}

\vspace{0.1cm}

\noindent where $x_l \in \Phi_l$  and $x_k \in \Phi_k$ with the reflexivity (loop) convention adopted i.e. $[Max]_{n,n}=1$.

\vspace{0.2cm}

\noindent \textbf{Farewell Exercises.}

\vspace{0.2cm}

\noindent \textbf{Problem-Exercise 5.1.} Rewrite Markov property in $F$-nomials language.

\vspace{0.2cm}

\noindent \textbf{Problem-Exercise 5.2.} Find the inverse of $\fnomial{l}{k}$ using the Theorem 4 and the knowledge of $[Max]^{-1}$. Compare  with [11].

\vspace{0.3cm}

\noindent \textbf{Acknowledgments}
\vspace{0.1cm}
\noindent Thanks are expressed here to the Student of Gda\'nsk University Maciej Dziemia\'nczuk for applying his skillful   TeX-nology with respect most of my articles since three years as well as for his general assistance and cooperation on KoDAGs  investigation. 

\vspace{0.2cm}

\noindent The author expresses his gratitude  also   Dr Ewa Krot-Sieniawska for her several years' cooperation and vivid application  of the alike  material deserving  Students' admiration for her being such a comprehensible and reliable  Teacher before she as Independent Person was fired by local Bialystok University local authorities exactly on the day she had defended  Rota and cobweb posets related dissertation with distinction.

\end{document}